\documentclass[11pt,a4paper]{article}

\usepackage[labelfont=bf, margin=1em]{caption}

\usepackage[T1]{fontenc}
\usepackage[american]{babel}
\usepackage{lmodern} 
\usepackage{microtype} 

\usepackage{amsmath, amssymb, amsfonts, amsthm}
\usepackage{mathtools} 
\usepackage{mathrsfs}  
\usepackage{textcomp}  
\usepackage{siunitx}   

\usepackage{graphicx}
\usepackage{float}
\usepackage{subcaption}
\usepackage{fullpage}  
\usepackage{appendix}  
\usepackage{setspace}  
\usepackage{tikz}
\usepackage[section]{placeins} 

\usepackage{booktabs}  
\usepackage{array}
\usepackage{tabularx}  
\usepackage[table]{xcolor} 

\usepackage{enumitem}  
\usepackage[skip=0.6em, indent=0pt]{parskip} 

\usepackage[numbers,sort&compress]{natbib}
\usepackage{url}
\usepackage{hyperref}
\hypersetup{
  colorlinks=true,
  linkcolor=blue!50!black,
  citecolor=blue!50!black,
  urlcolor=blue!50!black
}

\usepackage{authblk}   

\numberwithin{equation}{section}
\theoremstyle{definition}

\newcommand{\R}{\mathbb{R}}

\newcolumntype{Y}{>{\raggedright\arraybackslash}X}

\title{\LARGE \textbf{Early warning of critical transitions: Distinguishing supercritical and subcritical bifurcations from noisy observations}}

\author[1,*]{Paul A. Sanders}
\author[2]{Christian Kuehn}
\author[2]{Bastian Hilder}
\author[1]{Robbin Bastiaansen}

\affil[1]{\textit{Mathematical Institute (MI), Institute for Marine and Atmospheric Research Utrecht (IMAU) \&
Centre for Complex Systems Studies (CCSS), Utrecht University, 3584 CD Utrecht, The Netherlands}}
\affil[2]{\textit{Chair of Multiscale and Stochastic Dynamics, Department of Mathematics,
TUM School of Computation, Information and Technology,
Technical University of Munich, 85748 Garching, Germany}}

\affil[*]{\textit{Corresponding author: p.a.sanders@uu.nl}}

\date{}

\begin{document}
\maketitle

\begin{abstract}
Early-warning methods can indicate that a system is approaching a loss of stability, but they do not usually reveal whether the approaching bifurcation is supercritical or subcritical. We investigate whether this distinction can be made from noisy pre-transition observations. For this investigation, we extend a regression-based stability method by including the leading cubic nonlinear term, whose sign is used to predict the bifurcation type. We test the method on a stochastic pitchfork model and a spatial stochastic Swift--Hohenberg equation.
In both models, classification improves with longer observation horizons, closer proximity to the bifurcation, stronger noise, and larger nonlinear effects. However, these conditions are also linked to a higher risk of noise-induced tipping in the subcritical case. The useful early-warning window is therefore determined by a balance between obtaining enough information to identify the bifurcation type and doing so before noise-induced tipping becomes too likely. In spatially extended systems, we find that observing the full spatial field allows classification over shorter observation horizons, and performance improves when a broader range of spatial modes contributes to the observed dynamics. Altogether, we identify the conditions under which the method can reliably distinguish between supercritical and subcritical bifurcations.
\end{abstract}

\begin{quotation}
\noindent\small\bfseries
Many natural and engineered systems, from ice sheets and ocean currents to
ecosystems and power grids, can change abruptly when slowly changing
conditions push them past a critical threshold. Existing early-warning
indicators can signal that a system is losing stability, but they cannot
tell whether the coming change will be gradual or abrupt. This difference
depends on whether the underlying bifurcation is supercritical, where a new
state grows smoothly out of the old one, or subcritical, where the system
can jump to a distant state and random fluctuations can trigger that jump
before the threshold itself is reached. We show that noisy observations
recorded before the transition can contain enough information to tell these
two cases apart. Fitting a simple model that includes the leading nonlinear
term to simulated stochastic data, we find that the distinction becomes
reliable for longer records, closer to the transition, and for stronger noise
or nonlinearity, but these conditions also raise the risk of an early
noise-induced jump. For spatially extended systems, observing the full
spatial field shortens the record needed, especially when fluctuations are
spread over many spatial scales.
\end{quotation}

\section{Introduction}

An important motivation for studying bifurcations is the understanding of
tipping points in complex systems. A tipping point is a critical
threshold at which a system undergoes a qualitative transition to a substantially
different dynamical regime \cite{Scheffer2009, van2016you, Lenton2008}. Such transitions occur in
many natural and engineered systems, including climate dynamics, ecology,
financial markets, infrastructure networks, and social systems
\cite{Scheffer2009,Artime2024,Vespignani2012}. Examples include shifts in ocean circulation \cite{Lenton2008},
tipping of the Greenland and Antarctic ice sheets
\cite{Pattyn2018}, abrupt shifts in arid ecosystems
\cite{Scheffer2001}, cascading failures in power grids
\cite{Artime2024}, financial market crashes
\cite{Sornette2003}, and abrupt changes in social conventions or the
spreading of information through social systems
\cite{Centola2018,Vespignani2012}. In many systems,
the main question is therefore not only when a transition occurs, but also
whether it develops gradually or through an abrupt shift between states.

\begin{figure}[t]
\centering
\includegraphics[width=0.9\textwidth]{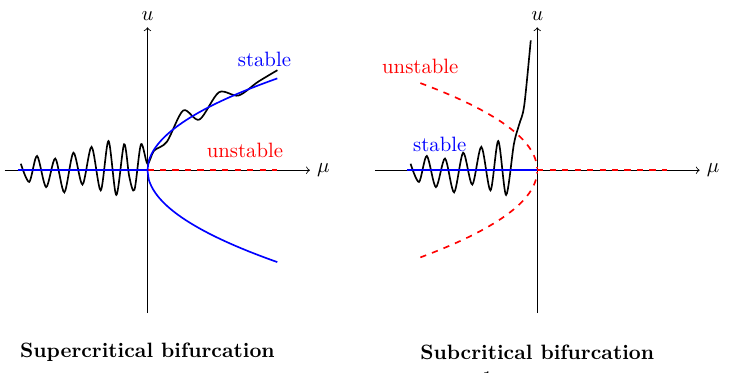}
\caption{
Schematic comparison of noisy trajectories near supercritical and subcritical pitchfork bifurcations. Blue solid lines denote stable equilibrium branches, red dashed lines denote unstable branches, and the black curves illustrate noisy trajectories. In the supercritical case, stable finite-amplitude branches emerge after the loss of stability, whereas in the subcritical case unstable finite-amplitude branches are already present for $\mu<0$.
}
\label{fig:subsupercritical}
\end{figure}

Specifically, in the case that the system undergoes a pitchfork bifurcation with normal form \eqref{eq:normal_form}, the distinction between gradual and abrupt transitions is closely
related to whether the underlying bifurcation is supercritical
or subcritical; the qualitative
difference is illustrated in Figure~\ref{fig:subsupercritical}, where $u = u(t) \in \R$ denotes the unknown variable and $\mu \in \R$ is a parameter in the system that encodes the distance to the bifurcation point. In a
supercritical bifurcation, stable finite-amplitude states emerge
continuously when the critical threshold is crossed. In a subcritical
bifurcation, unstable finite-amplitude states already exist before the
bifurcation point \cite[Sec.~3.4]{Strogatz2015}. These states form a threshold that can be crossed by
sufficiently large perturbations. Stochastic fluctuations can trigger such abrupt shifts before the deterministic bifurcation point is reached, a process known as noise-induced tipping \cite{Ashwin2012}. The distinction between supercritical and subcritical
bifurcations is therefore important both for the nature of the transition
and for the sensitivity of the system to stochastic perturbations.

The distinction between supercritical and subcritical bifurcations also applies to spatially extended systems. In these systems, instabilities can give rise to spatial patterns, as in a Turing bifurcation \cite{Turing1952Morphogenesis,CrossGreenside2009,Rietkerk2021}. The dynamics close to such an instability are determined by a universal amplitude equation, the (real) Ginzburg–Landau equation \cite{CrossGreenside2009}, whose structure closely resembles the normal form of the pitchfork bifurcation. The formal derivation and rigorous justification of these models in pattern-forming systems is by now a well-established theory; see, for example, \cite{schneider2017book} for the deterministic case, and \cite{Blomker2007} for the stochastic case.
Therefore, as in the pitchfork case depicted in Figure \ref{fig:subsupercritical}, the pattern emerges continuously from the homogeneous state for a supercritical bifurcation, and it typically appears through an abrupt finite-amplitude transition for a subcritical bifurcation.

Traditional early-warning indicators, such as those based on the theory of critical slowing down,
identify changes in the linearized dynamics near a stable equilibrium
\cite{Scheffer2009, van2007slow, wissel1984universal, dakos2011slowing}. These changes can signal a loss of stability and thereby provide an early warning of an approaching bifurcation. However, these indicators do not distinguish between supercritical and subcritical
bifurcations because this distinction is determined by the nonlinear
dynamics. In fact, leading nonlinear terms can already influence a system's behavior and predictability before the bifurcation is crossed \cite{RomanoKuehn2018}, suggesting that pre-transition observations may contain information beyond the linearized dynamics.
This raises the question under which observational conditions these
nonlinear contributions can be recovered reliably, thus enabling an early-warning indicator that is capable of distinguishing imminent supercritical bifurcations from subcritical bifurcations.

In this work, we investigate whether noisy observations of a system approaching a bifurcation can be used to determine its criticality. We extend the
regression method introduced in
\cite{SANDERS2026135300}, which estimated the linearized dynamics, by also
estimating the leading cubic drift term from temporal and spatio-temporal
data (Section \ref{sec:theory_method}). We consider synthetic data generated from a stochastic variant of the normal form for a pitchfork bifurcation \eqref{eq:normal-form-stochastic-pitchfork}, and a
stochastic Swift--Hohenberg equation \eqref{eq:SH_dimensional}. The latter is a standard model for a pattern-forming system close to a Turing instability, and was first introduced as a model for temperature-driven convection in fluids \cite{swift1977hydrodynamic}. For both models, we study how the
prediction depends on the distance to the bifurcation and the observation
horizon (Sections~\ref{sec:Stochastic_pitchfork}--\ref{sec:Swift-Hohenberg}). For the Swift--Hohenberg equation, we additionally investigate how the strength of the
spatial operator and the critical spatial scale affect the prediction of the
bifurcation criticality (Section \ref{sec:Swift-Hohenberg}). Finally, in Section \ref{sec:discussion} we discuss the
feasibility and limitations of applying the method to observational data from complex systems.

\section{Distinguishing supercritical and subcritical bifurcations}
\label{sec:theory_method}
To predict the bifurcation criticality from data, we approximate the local
dynamics by a minimal normal form containing the linear term and the
leading nonlinear contribution. For the bifurcations considered in this paper, the
leading nonlinear term is cubic. For non-spatial models, we thus consider
\begin{align}
\frac{du}{dt} = \mu u + s u^3,
\label{eq:normal_form}
\end{align}
where $u = u(t) \in \R$ denotes the deviation from the equilibrium, $\mu \in \R$ describes its
linear stability, and $s \in \R \setminus \{0\}$ is the leading nonlinear coefficient. We
use $\mu$ and $s$ throughout to emphasize their roles as the stability
parameter and the coefficient determining the bifurcation criticality,
respectively. For the equilibrium $u=0$, $\mu<0$ indicates linear stability, and
the bifurcation occurs as $\mu$ crosses zero. The sign of $s$ determines
the bifurcation criticality: $s<0$ corresponds to a supercritical
bifurcation and $s>0$ to a subcritical bifurcation, as illustrated
in Figure~\ref{fig:subsupercritical}. 

Equation~\eqref{eq:normal_form} is the normal form of a pitchfork bifurcation and describes the leading-order dynamics close to the bifurcation after reduction to the critical mode \cite{wiggins2003book}. A similar cubic structure arises in spatially extended systems close to pattern-forming instabilities.

For spatially extended systems, we consider the cubic Swift--Hohenberg equation \cite{swift1977hydrodynamic}
\begin{equation}
\label{eq:Swift-Hohenberg}
\partial_t u = \mu u - \Gamma (k_\ast^2+\partial_x^2)^2u + s u^3,
\end{equation}
where $u=u(x,t)\in\R$, $x\in\R$, $t\geq0$, and $\mu\in\R$, $\Gamma>0$, and $s\in\R\setminus\{0\}$. At $\mu=0$, the homogeneous equilibrium $u=0$ loses stability to spatial modes with critical wavenumber $k_\ast$, resulting in a Turing instability. Close to this instability, the dynamics are dominated by these critical modes, whose amplitude varies slowly in space and time. At leading order, this slow modulation is described by a cubic Ginzburg--Landau-type amplitude equation \cite{CrossGreenside2009,schneider2017book}. For a spatially uniform amplitude with fixed phase, the spatial derivative terms in the amplitude equation vanish, leaving a cubic ordinary differential equation with the same normal-form structure as Equation~\eqref{eq:normal_form} \cite{CrossGreenside2009}. A non-zero amplitude corresponds to a spatially periodic pattern in the original Swift--Hohenberg system, and the cubic coefficient in the amplitude equation is proportional to $s$. The sign of $s$ therefore determines whether the patterned branch emerges supercritically or subcritically, just as in the non-spatial normal form. Since such amplitude equations arise generically close to stationary finite-wavenumber instabilities, the Swift--Hohenberg equation provides a representative test case for a broader class of pattern-forming systems \cite{CrossGreenside2009}.

Throughout this paper, we use \emph{bifurcation criticality classification} to
refer to distinguishing between supercritical and subcritical bifurcations
based on the estimated cubic coefficient $s$. The classification performance
therefore measures how well the two bifurcation types can be separated
from the fitted values of this coefficient $s$, which we denote by $\hat{s}$.

We estimate these coefficients by writing the governing dynamics as a regression problem, similar to SINDy methods \cite{brunton2016discovering}
\begin{align}
\frac{du}{dt} = \Theta F(u),
\label{eq:regression_general}
\end{align}
where $F(u)$ contains the selected candidate terms and $\Theta$ contains
their coefficients. Note that, although the observations are noisy, we do not
explicitly include a stochastic forcing in the regression model.
Instead, we fit the deterministic drift of the dynamics. For zero-mean
stochastic forcing, the noise contribution is expected to average out over
sufficiently many observations or realizations and is therefore treated as
a residual in the regression. This allows us to focus on the deterministic
coefficients that determine the qualitative properties of the transition,
while keeping the fitted model as simple as possible.

We first describe the method for ordinary differential equations and then extend it to spatially extended systems.

\subsection{Method for ordinary differential equations}
\label{sec:method_ode}

We assume that we have scalar observations $Y(t_i)\in\mathbb{R}$ of a system at times
$\{t_1,\ldots,t_{N_t}\}$ and seek to describe the dynamics using
the scalar normal form in Equation~\eqref{eq:normal_form}.

The data are divided into observation windows that are analyzed separately.
Within each window, the dynamics are assumed to be approximately
stationary. The linear and cubic coefficients are estimated using linear
regression, after which the stability and the predicted bifurcation
criticality are determined.

\subsubsection*{Step 1: Prepare the data and define the observation windows}

The fluctuations around the equilibrium are defined as
\begin{align}
u(t_i)=Y(t_i)-Y^*(t_i),
\label{eq:ode_fluctuations}
\end{align}
where $Y^*(t_i)$ is the equilibrium state. For the models considered in this paper, the equilibrium is known to be located at $u=0$. For observational data with an unknown equilibrium, it can instead be estimated from the data.

From the full dataset, we define $N_{\mathrm{w}}$ possibly overlapping observation windows,
\begin{align}
\mathcal{W}_{\ell}
=
\left\{
t_1^{(\ell)},\ldots,t_{N_t^{(\ell)}}^{(\ell)}
\right\}
\subseteq
\left\{
t_1,\ldots,t_{N_t}
\right\},
\qquad
\ell=1,\ldots,N_{\mathrm{w}},
\label{eq:ode_windows}
\end{align}
where $N_t^{(\ell)}$ denotes the number of temporal observations in window $\ell$. The number of observations may therefore differ between windows, and each window is analyzed independently.

\subsubsection*{Step 2: Construct the regression model}

The minimal library contains the linear and cubic terms, i.e.
\begin{align}
\mathbf{F}(u)=
\begin{bmatrix}
u\\
u^3
\end{bmatrix},
\qquad
\Theta=\begin{bmatrix}\mu & s\end{bmatrix}.
\label{eq:ode_library}
\end{align}
Within observation window $\ell$, the time derivative is estimated using forward differences,
\begin{align}
u_t\left(t_i^{(\ell)}\right)
\approx
\frac{
u\left(t_{i+1}^{(\ell)}\right)-u\left(t_i^{(\ell)}\right)
}{
t_{i+1}^{(\ell)}-t_i^{(\ell)}
},
\qquad
i=1,\ldots,N_t^{(\ell)}-1.
\label{eq:ode_time_derivative}
\end{align}

\subsubsection*{Step 3: Estimate the coefficients}

For each observation window, the coefficients are estimated by minimizing
\begin{align}
\left(\hat{\mu}_{\ell},\hat{s}_{\ell}\right)
=
\operatorname*{argmin}_{\tilde{\mu},\tilde{s}}
\sum_{i=1}^{N_t^{(\ell)}-1}
\left|
u_t\left(t_i^{(\ell)}\right)
-
\tilde{\mu}u\left(t_i^{(\ell)}\right)
-
\tilde{s}u^3\left(t_i^{(\ell)}\right)
\right|^2.
\label{eq:ode_least_squares}
\end{align}
This gives one estimate of $\mu$ and $s$ for each observation window.

\subsubsection*{Step 4: Determine the bifurcation type and linear stability}

The estimate $\hat{\mu}_{\ell}$ is an indicator for the linear stability, similar to classic early-warning signs \cite{Scheffer2009, van2007slow, wissel1984universal, dakos2011slowing}, with $\hat{\mu}_{\ell}<0$ corresponding to a stable equilibrium and values approaching zero from below indicating an approaching loss of stability. The estimated cubic coefficient $\hat{s}_{\ell}$ is an indicator for the bifurcation criticality. In an ideal situation, a negative value, $\hat{s}_{\ell}<0$, indicates a supercritical bifurcation, and a positive value, $\hat{s}_{\ell}>0$, a subcritical bifurcation. However, due to neglected effects (e.g. stochasticity, temporal variations of forcing within a time window), the best criticality classification threshold might be located at a different value for $\hat{s}_{\ell}$.

\subsection{Method for spatially extended systems}
\label{sec:method_pde}

For spatially extended systems, we assume scalar observations $Y(x_j,t_i)\in\mathbb{R}$ at times
$\{t_1,\ldots,t_{N_t}\}$ and spatial locations
$\{x_1,\ldots,x_{N_x}\}$. The formulation presented here considers a single spatial variable and one
spatial dimension, but can be extended to multiple variables and higher
spatial dimensions.

As in the ODE case, the data are divided into temporal observation windows
that are analyzed separately, and the coefficients are assumed to remain
approximately constant within each window. The method builds on the spatial stability framework introduced in
\cite{SANDERS2026135300}. In the present formulation, spatial derivatives are evaluated
spectrally using Fourier transforms, while the fitted dynamics are extended
to include the leading cubic nonlinear term. The cubic coefficient is used
to predict the bifurcation type, while the linear spatial coefficients are
used to construct the dispersion relation.

\subsubsection*{Step 1: Prepare the data and define the observation windows}

The fluctuations around the homogeneous equilibrium are defined as
\begin{align}
u(x_j,t_i)
=
Y(x_j,t_i)-Y^*(t_i).
\label{eq:pde_fluctuations}
\end{align}
If the equilibrium is unknown, it can be estimated, e.g. from the spatial
average,
\begin{align}
Y^*(t_i)
=
\frac{1}{N_x}
\sum_{j=1}^{N_x}
Y(x_j,t_i).
\label{eq:pde_spatial_mean}
\end{align}

We define $N_{\mathrm{w}}$ possibly overlapping temporal observation
windows,
\begin{align}
\mathcal{W}_{\ell}
=
\left\{
t_1^{(\ell)},\ldots,t_{N_t^{(\ell)}}^{(\ell)}
\right\}
\subseteq
\left\{
t_1,\ldots,t_{N_t}
\right\},
\qquad
\ell=1,\ldots,N_{\mathrm{w}},
\label{eq:pde_windows}
\end{align}
where $N_t^{(\ell)}$ denotes the number of temporal observations in
window $\ell$. Each temporal observation in the window contains
measurements at all available spatial locations.

\subsubsection*{Step 2: Estimate temporal and spatial derivatives}

The temporal derivative is approximated from consecutive observations
using a forward difference,
\begin{align}
u_t(x_j,t_i)
\approx
\frac{
u(x_j,t_{i+1})-u(x_j,t_i)
}{
t_{i+1}-t_i
}.
\label{eq:pde_time_derivative}
\end{align}

Spatial derivatives are evaluated spectrally. Denoting the spatial discretized Fourier transform by $\hat{u}(k,t)=\mathcal{F}_x\{u(x,t)\}$, derivatives
of order $m$ are calculated according to
\begin{align}
\frac{\partial^m u(\cdot,t_i)}{\partial x^m}
=
\mathcal{F}_x^{-1}
\left\{
(ik)^m\hat{u}(k,t_i)
\right\}.
\label{eq:fourier_derivative}
\end{align}
This avoids additional finite-difference approximation errors in the
spatial derivatives, which is particularly useful when fitting dynamics
from noisy observations.

\subsubsection*{Step 3: Construct and fit the spatial regression model}

The candidate library contains the linear spatial terms expected in the
governing equation together with the leading cubic nonlinear term. In one
spatial dimension, it can be written as
\begin{align}
\mathbf{F}(u)
=
\begin{bmatrix}
u\\
u_x\\
u_{xx}\\
\vdots\\
u^3
\end{bmatrix},
\qquad
\Theta
=
\begin{bmatrix}
a_0 & a_1 & a_2 & \cdots & s
\end{bmatrix},
\label{eq:pde_spatial_library}
\end{align}
where dots represent higher order derivatives when required; only the spatial derivative terms corresponding to the assumed governing equation are included in the regression.

The coefficients are estimated separately for each observation window
using all spatial locations and time steps within that window. The
least-squares estimate is given by
\begin{align}
\hat{\Theta}_{\ell}
=
\operatorname*{argmin}_{\Theta}
\sum_{i=1}^{N_t^{(\ell)}-1}
\sum_{j=1}^{N_x}
\left|
u_t
\left(
x_j,t_i^{(\ell)}
\right)
-
\Theta
\mathbf{F}\left(u\left(
x_j,t_i^{(\ell)}
\right)\right)
\right|^2.
\label{eq:pde_least_squares}
\end{align}
This gives one estimate of the linear spatial coefficients $a_0, a_1, \ldots$ and the cubic
coefficient $s$ for each observation window.

\subsubsection*{Step 4: Determine the bifurcation type and construct the dispersion relation}

As before the estimated cubic coefficient $\hat{s}_{\ell}$ is an indicator for the bifurcation criticality, with negative values suggesting supercritical bifurcations and positive values subcritical bifurcations.

For spatial systems, the linear stability can be determined via the dispersion relation \cite{SANDERS2026135300}. For example, for the Swift--Hohenberg model (Section \ref{sec:Swift-Hohenberg}), the fitted linear operator is
\begin{align}
\hat{a}_{0,\ell}
+
\hat{a}_{2,\ell}\partial_x^2
+
\hat{a}_{4,\ell}\partial_x^4,
\end{align}
which gives the estimated dispersion relation
\begin{align}
\hat{\lambda}_{\ell}(k)
=
\hat{a}_{0,\ell}
-
\hat{a}_{2,\ell}k^2
+
\hat{a}_{4,\ell}k^4.
\label{eq:estimated_dispersion_sh}
\end{align}
The dominant wavenumber and corresponding growth rate are then obtained
from
\begin{align}
\hat{k}_{*,\ell}
&=
\operatorname*{argmax}_{k}
\operatorname{Re}
\left\{
\hat{\lambda}_{\ell}(k)
\right\},
&
\hat{\lambda}_{*,\ell}
&=
\operatorname{Re}
\left\{
\hat{\lambda}_{\ell}
\left(
\hat{k}_{*,\ell}
\right)
\right\}.
\end{align}
The value of $\hat{\lambda}_{*,\ell}$ describes the stability of the
dominant spatial mode, while $\hat{k}_{*,\ell}$ identifies the dominant
wavenumber and hence the corresponding characteristic spatial scale.

\section{Criticality of pitchfork bifurcation}
\label{sec:Stochastic_pitchfork}

We first apply and test the bifurcation criticality classification method to the stochastic pitchfork model. This provides
a controlled setting in which the linear coefficient and the sign of the
cubic term are known exactly. This section is structured as follows. First, we briefly describe the model. Second, we give details for the numerical experimental setup. Third, we give the results of these experiments.

\subsection{Model}
\label{sec:PF_model}

We consider the stochastic pitchfork normal form with additive noise, i.e.
\begin{align}\label{eq:normal-form-stochastic-pitchfork}
du=\left(ru+\alpha u^3\right)dt+\sigma\,dW_t,
\end{align}
where $W_t$ is a standard Wiener process. We introduce the rescaled state
and time variables
\begin{align}
u=\left(\frac{\sigma^2}{|\alpha|}\right)^{1/4}q,
\qquad
t=\frac{1}{\sigma\sqrt{|\alpha|}}\tau.
\label{eq:rescaled_variables}
\end{align}
to obtain the rescaled model
\begin{align}
dq=\left(\mu q+s q^3\right)d\tau+dW_\tau,
\label{eq:reduced_sde}
\end{align}
with
\begin{align}
\mu=\frac{r}{\sigma\sqrt{|\alpha|}},
\qquad
s=\operatorname{sgn}(\alpha)\in\{-1,1\}.
\label{eq:rescaled_parameters}
\end{align}
This scaling normalizes both the noise strength and the magnitude of the
cubic coefficient, leaving $\mu$ and $s$ as the reduced dynamical
parameters. The parameter $\mu$ gives the linear stability relative to
the noise and nonlinear scales. The bifurcation remains located at
$\mu=0$. The sign $s$ determines the bifurcation type, with $s=-1$
corresponding to the supercritical case and $s=1$ to the subcritical
case.

\begin{figure}[tbp]
\centering
\includegraphics[width=0.9\textwidth]{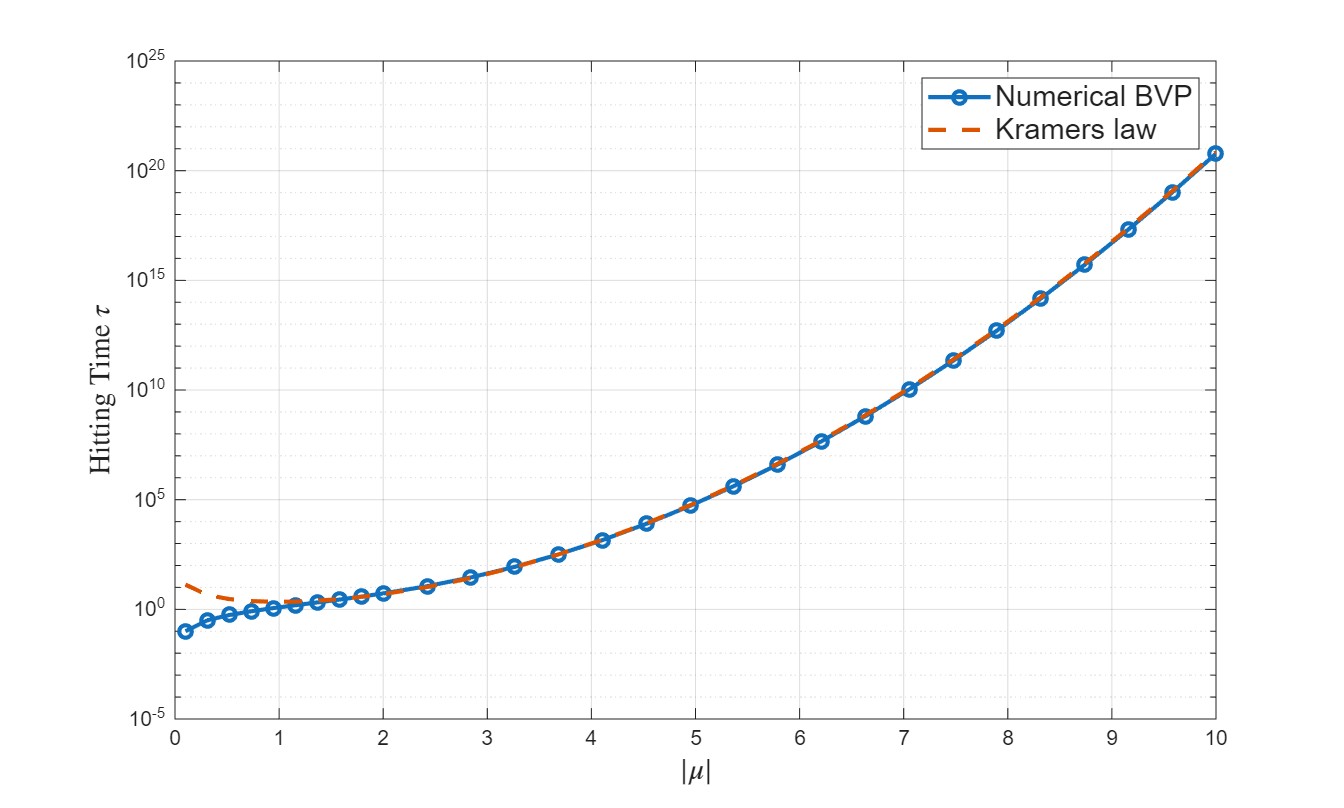}
\caption{Mean first passage time for the subcritical reduced pitchfork model in Equation~\eqref{eq:reduced_sde}. Starting from the stable equilibrium $q=0$, the first passage time measures the average time required for a trajectory to reach either unstable threshold $q=\pm\sqrt{-\mu}$ and leave the local stable region. The solid curve shows the mean first passage time obtained by numerically solving the backward boundary-value problem derived in Appendix~\ref{app:hitting_time}. The dashed curve shows the Eyring--Kramers approximation \cite{Kramers1940,berglund2013kramers}, multiplied by an empirical factor $C=0.6$ to improve the comparison over the finite parameter range shown.}
\label{fig:mfpt_kramers}
\end{figure}

For the subcritical model, the stable equilibrium $q=0$ is bounded by the
unstable states $q=\pm\sqrt{-\mu}$ for $\mu<0$. A trajectory that crosses either threshold leaves the region in which the cubic normal form provides a bounded local description. Appendix~\ref{app:hitting_time} derives the mean first passage time to these thresholds using the standard backward equation \cite{Gardiner1985} and compares it with the Eyring--Kramers approximation \cite{Kramers1940,berglund2013kramers}. Based on these calculations, Figure~\ref{fig:mfpt_kramers} shows how the escape time depends on $\mu$. This provides the basis for the treatment of escaped trajectories in the numerical experiments described in Section~\ref{sec:PF_setup}.

\subsection{Experimental setup}
\label{sec:PF_setup}

We generate numerical data from the reduced stochastic pitchfork model in
Equation~\eqref{eq:reduced_sde} using the Euler--Maruyama method,
\begin{align}
q_{i+1}
=
q_i
+
\left(
\mu q_i+s q_i^3
\right)\Delta\tau
+
\sqrt{\Delta\tau}\,\xi_i,
\label{eq:pitchfork_euler_maruyama}
\end{align}
where $\xi_i\sim\mathcal{N}(0,1)$ are sampled independently at every time
step. Each simulation starts at the stable equilibrium, i.e. $q(0)=0$.

To study how the fitting performance depends on the distance to the
bifurcation point ($\mu$), the bifurcation criticality $(s)$, and the available observation
time horizon ($\tau_{\max}$), we vary these parameters as
summarized in Table~\ref{tab:pitchfork_parameters}.

\begin{table}[tbp]
\centering
\caption{
Parameters varied in the stochastic pitchfork experiments. The fixed
numerical settings are the time step $\Delta\tau=0.01$, initial condition
$q(0)=0$, unit noise strength in the reduced model, and
$N_{\mathrm{sim}}=100$ independent realizations for each parameter
combination. Simulations are performed using the Euler--Maruyama method.
}
\begin{tabularx}{\textwidth}{ll c}
\hline
Parameter & Description & Values\\
\hline
$\mu$ & scaled linear stability
& $50$ log-spaced in $[-10^2,-10^{-1}]$\\
$s$ & cubic coefficient sign
& $\{-1,1\}$\\
$\tau_{\max}$ & observation horizon
& $\{1,2,5,10,20,50,100,200,500,1000,2000,5000\}$\\
\hline
\end{tabularx}
\label{tab:pitchfork_parameters}
\end{table}

For each combination of $\mu$ and $s$, the different observation horizons
are compared using the same underlying stochastic trajectories. A common
ensemble of $100$ trajectories is used for all included observation
horizons, with the shorter datasets obtained by truncating the trajectories
at the corresponding value of $\tau_{\max}$. This choice keeps the
underlying noise realizations fixed, so that differences between
observation horizons result from the amount of available data rather than
from changes in the sampled stochastic trajectories.

In the subcritical case, some trajectories eventually leave the local
regime described by the cubic normal form. Since the fitting procedure
itself does not assume prior knowledge of the unstable thresholds
$q=\pm\sqrt{-\mu}$, departures are identified using a numerical criterion:
a trajectory is excluded from the analysis if non-finite values occur or if
\begin{align}
\max_i |q_i|>10^{10}.
\label{eq:pitchfork_escape_criterion}
\end{align}
For each parameter setting, an observation horizon is included in the
analysis only if at least $10\%$ of the trajectories generated from the
first $100$ random seeds are not excluded by this criterion before reaching
that horizon. For each included observation horizon, a common ensemble of
$100$ trajectories is selected from a larger finite pool of simulated
trajectories, excluding those that leave the local regime before the end of
the observation window (with shorter observation horizons obtained by truncation). As $\mu$ approaches the bifurcation point at $\mu=0$ from
below, satisfying this criterion becomes harder for long
observation horizons, consistent with the decrease in the mean first
passage time shown in Figure~\ref{fig:mfpt_kramers}.

For each trajectory included in the analysis, the
coefficients are estimated using the ODE fitting procedure described in
Section~\ref{sec:method_ode}. 
This gives one pair of estimates
$(\hat{\mu},\hat{s})$ for each realization and observation
horizon.

The estimates $\hat{s}$ are used to distinguish the two bifurcation types.
The subcritical case is predicted when $\hat{s}$ exceeds a classification
threshold. To test this procedure and threshold, for each value of $\mu$ and observation horizon $\tau_{\max}$ we consider an ROC classification ensemble that contains the 100 simulations for both $s = -1$ and $s = +1$ (as described above) and test how well the bifurcation criticality classification can distinguish the super- and subcritical cases. In particular, we vary the classification threshold to obtain receiver operating
characteristic (ROC) curves \cite{Fawcett2006,zhang2015predictability}. The corresponding area under the curve (AUC)
measures how well the distributions of $\hat{s}$ for the supercritical
and subcritical ensembles can be separated. An AUC of $0.5$ corresponds
to chance-level discrimination between the two bifurcation types, while an
AUC of $1$ corresponds to complete separation.

\subsection{Results}
\label{sec:PF_results}

In this section we focus on the criticality classification; full details of the distributions of the fitted coefficients $\hat{\mu}$ and $\hat{s}$
for representative parameter settings are shown in
Appendix~\ref{app:pitchfork}. In short, the linear coefficient $\hat{\mu}$ is
generally recovered well, with improved estimates for
longer observation horizons and closer to the bifurcation point.

\begin{figure}[tbp]
    \centering    \includegraphics[width=0.8\textwidth]{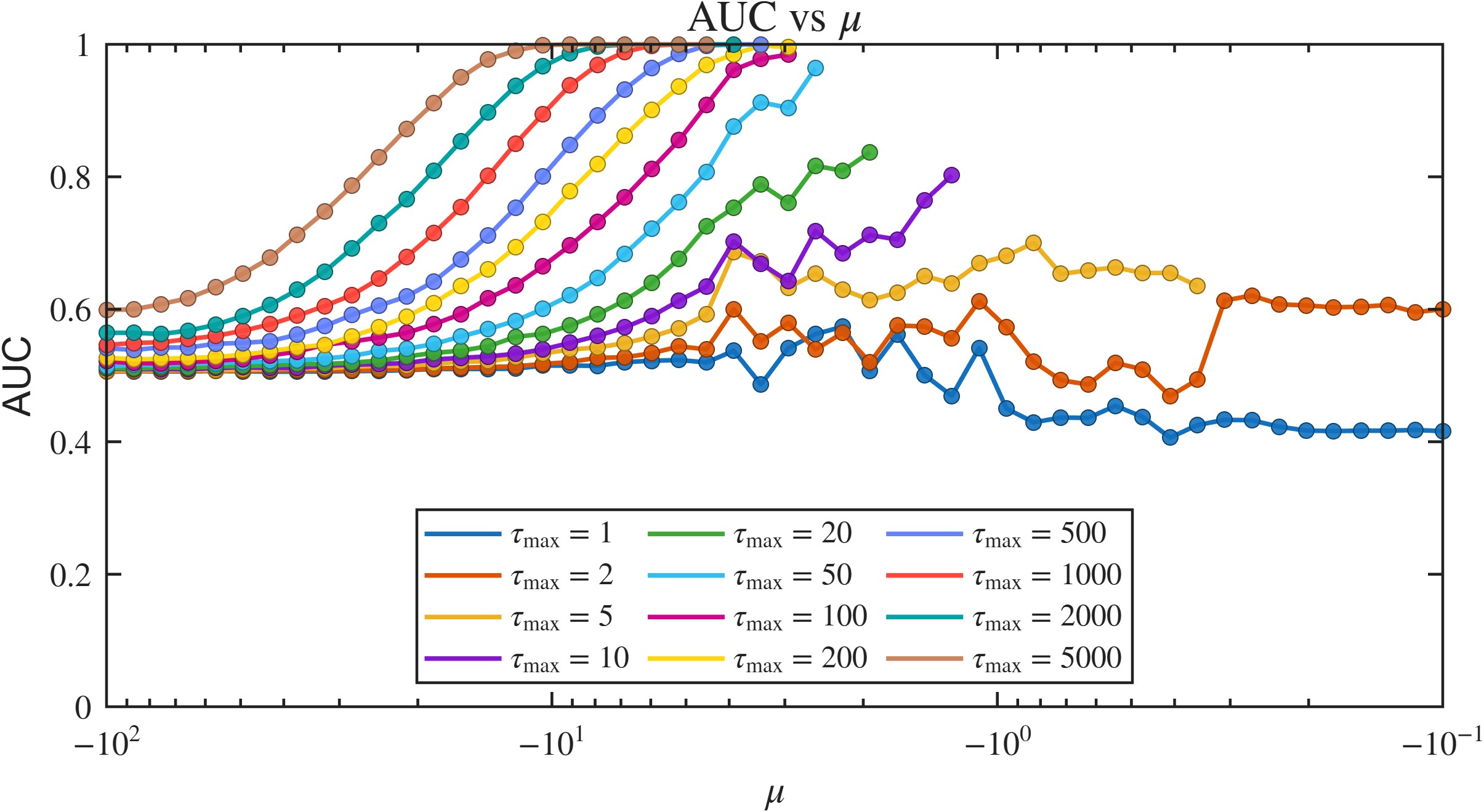}
    \caption{Area under the ROC curve (AUC) as a function of parameter $\mu$ for the different observation horizons $\tau_{\max}$. An AUC value of $0.5$ indicates the classification method is as good as flipping a coin; an AUC value of $1$ indicates potential for a perfect classification. $\mu$ indicates closeness to the bifurcation point located at $\mu = 0$. A value of $\mu$ is shown only if at least $10\%$ of the subcritical realizations remain within the local regime for the full observation horizon (Section~\ref{sec:PF_setup}); curves for longer horizons therefore end farther from the bifurcation. Some corresponding ROC-curves are shown in Figure~\ref{fig:pitchfork_combined_roc}.
}
    \label{fig:pitchfork_auc_vs_r}
\end{figure}

Figure~\ref{fig:pitchfork_auc_vs_r} shows the AUC for bifurcation criticality
classification as a function of the distance to the bifurcation, $\mu$,
and the observation horizon, $\tau_{\max}$. For a fixed observation
horizon, the AUC generally increases as $\mu$ approaches zero, showing that
the nonlinear coefficient becomes easier to distinguish closer to the
bifurcation. At the same time, for a fixed value of $\mu$, longer
observation horizons lead to higher AUC values. Hence, a longer observation
time allows the bifurcation criticality to be predicted farther away from the
bifurcation point. Some curves, particularly for the longer observation horizons, terminate
before reaching $\mu=0$ because too few subcritical trajectories remain
within the local regime for the full observation horizon.

The differences between the observation horizons are substantial. For the longest
observation horizon ($\tau_{\max}=5000$), the AUC already exceeds
approximately $0.6$ at $|\mu|\approx100$ and approaches perfect classification around
$|\mu|\approx10$. This means that accurate prediction is possible well before $|\mu|\approx3.5$ where too few trajectories remain within the local regime for the analysis, and noise-induced tipping is very probable. For intermediate observation horizons ($\tau_{\max}=100$--$500$), perfect classification ($\mathrm{AUC}=1$) can also be approached, but only relatively
close to the onset of frequent noise-induced tipping.
For shorter horizons ($\tau_{\max}=10$--$50$), the criticality classification becomes
weaker, and deteriorates further for even shorter observation windows, where the AUC
stays below approximately $0.7$ over most of the parameter range. For the shortest observation horizon ($\tau_{\max}=1$), the AUC even falls
below $0.5$, to approximately $0.42$, for $|\mu|\lesssim1$. At this horizon the fitted $\hat{s}$ is strongly negative for both bifurcation types (Appendix~\ref{app:pitchfork}), which indicates that the bifurcation criticality classification is unreliable when used with too little data.

\begin{figure}[tbp]
    \centering
    \includegraphics[width=0.48\textwidth]{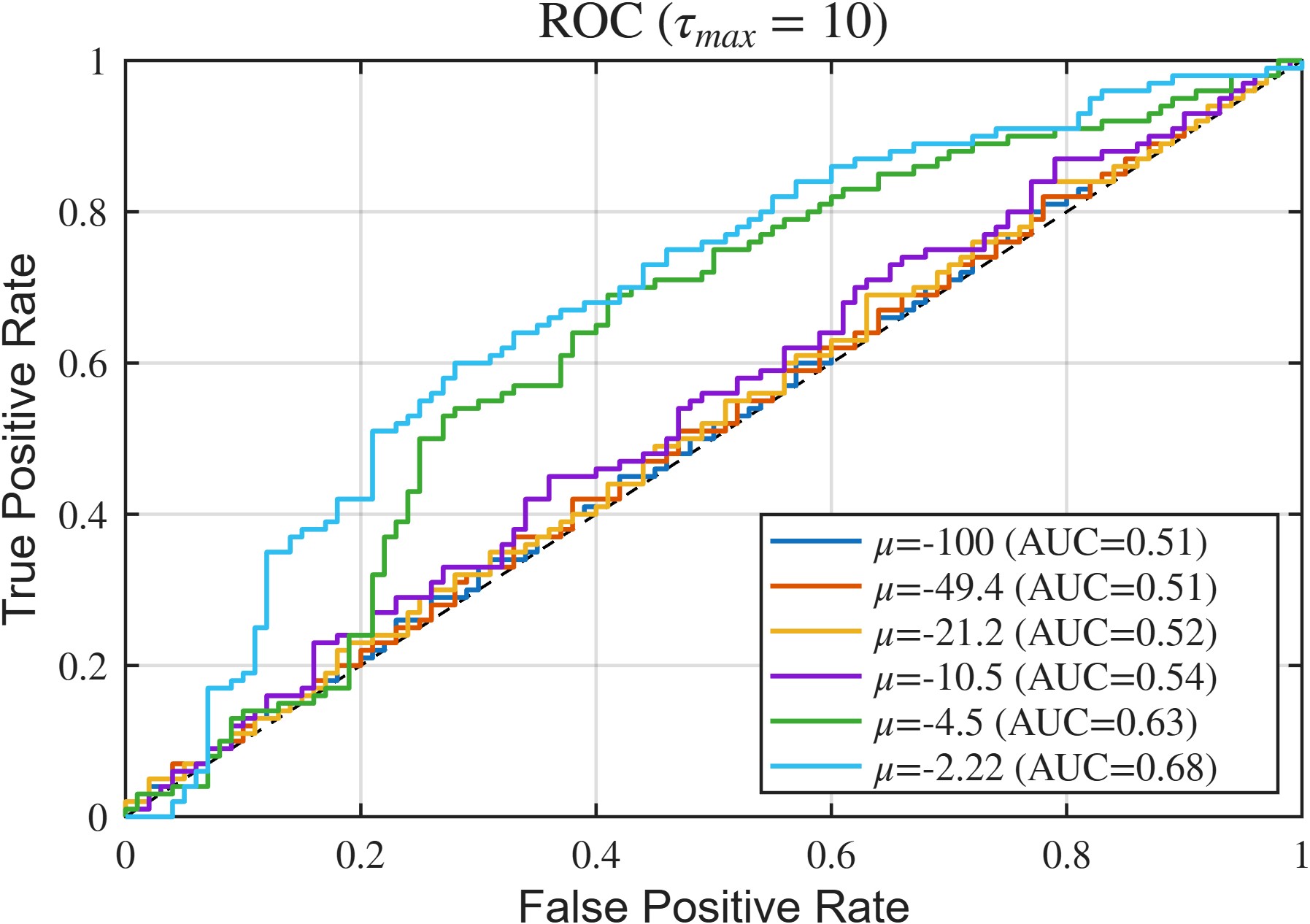}
    \includegraphics[width=0.48\textwidth]{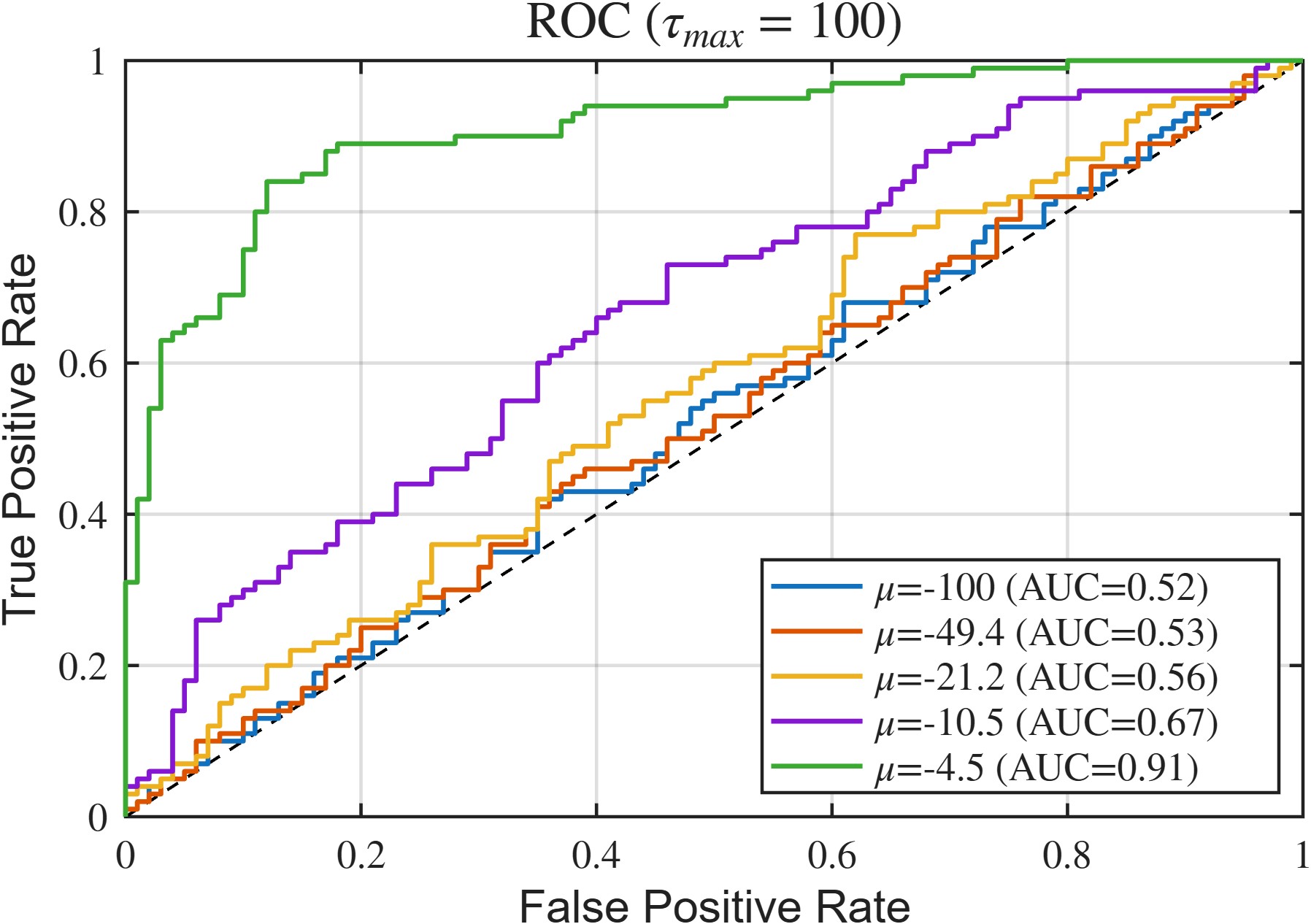}

    \caption{
Receiver operating characteristic (ROC) curves for the stochastic
pitchfork model in Equation~\eqref{eq:reduced_sde} at observation
horizons $\tau_{\max}=10$ (left) and $\tau_{\max}=100$ (right). Each
curve corresponds to a different value of $\mu$ and is constructed from
the fitted values of $\hat{s}$. The horizontal axis gives the
false-positive rate and the vertical axis the true-positive rate, with the
subcritical case ($s=1$) taken as the positive class. The diagonal line
denotes chance-level discrimination (i.e. flipping a coin).
}
    \label{fig:pitchfork_combined_roc}
\end{figure}

This dependence on the observation horizon is also visible in the
corresponding ROC curves in Figure~\ref{fig:pitchfork_combined_roc},
shown for two representative observation horizons ($\tau_{\max} = 10$ and $\tau_{\max}=100$) and a selection of $\mu$ values. For
$\tau_{\max}=10$, most curves remain relatively close to the diagonal over
much of the parameter range. For $\tau_{\max}=100$, a clearer separation
develops as $\mu$ approaches zero. The parameter distributions for $\hat{s}$ shown in
Appendix~\ref{app:pitchfork} give the same picture: the distributions of
$\hat{s}$ overlap strongly when little information is available and become
increasingly separated as the observation horizon increases or $\mu$ approaches zero, which makes it possible to choose the classification threshold to separate sub- and supercritical cases.

\section{Criticality of bifurcation in the Swift--Hohenberg equation}
\label{sec:Swift-Hohenberg}

Next, we apply and test the bifurcation criticality classification method to the stochastic Swift--Hohenberg model. In contrast to the stochastic pitchfork model, each temporal observation now contains the complete spatial field. This allows us to investigate how the additional spatial information affects the estimation of the nonlinear coefficient and the predicted bifurcation criticality.

\subsection{Model}
\label{sec:SH_model}

We consider the stochastic Swift--Hohenberg equation \cite{swift1977hydrodynamic}
\begin{align}
du =
\left[
ru
-
\Gamma
\left(k_\ast^2+\nabla_x^2\right)^2u
+
\alpha u^3
\right]dt
+
\sigma\,dW_t,
\label{eq:SH_dimensional}
\end{align}
on a bounded, periodic domain, where $u=u(x,t)$ denotes the spatially varying system state. Here, $r$ determines the linear
stability of the homogeneous equilibrium, $\alpha$ is the leading cubic
nonlinear coefficient, and $\sigma$ determines the noise strength. The
additional spatial parameters are $\Gamma$, which sets the strength of
the spatial operator, and $k_\ast$, the critical wavenumber at
which the linear growth rate is maximal. Here,
$\nabla_x^2=\partial_x^2$ denotes the spatial Laplacian in one dimension,
and $W_t$ is a space--time Wiener process.

We first rescale space by the critical physical wavenumber,
\begin{align}
X=k_\ast x.
\label{eq:SH_spatial_scaling}
\end{align}
The space--time white noise also changes under this spatial transformation.
In one spatial dimension, the transformation $X=k_\ast x$ introduces an
additional factor $\sqrt{k_\ast}$ in the noise amplitude
\cite[Chapter~1]{Blomker2007}. Accounting for this additional spatial
scaling of the noise, we choose the amplitude and time scales as
\begin{align}
u &=
\left(
\frac{\sigma^2k_\ast}{|\alpha|}
\right)^{1/4}q,
&
t &=
\frac{1}{\sigma\sqrt{k_\ast|\alpha|}}\tau.
\label{eq:SH_scaling}
\end{align}
This choice normalizes both the magnitude of the cubic coefficient and
the strength of the space--time white noise. The rescaled Swift--Hohenberg model then becomes
\begin{align}
dq =
\left[
\mu q
-
\gamma
\left(1+\nabla_X^2\right)^2q
+
sq^3
\right]d\tau
+
dW_\tau,
\label{eq:reduced_sh}
\end{align}
with
\begin{align}
\mu &=
\frac{r}{\sigma\sqrt{k_\ast|\alpha|}},
&
\gamma &=
\frac{\Gamma k_\ast^4}
{\sigma\sqrt{k_\ast|\alpha|}},
&
s &= \operatorname{sgn}(\alpha).
\label{eq:SH_scaled_parameters}
\end{align}
Here, $W_\tau$ denotes the correspondingly rescaled space--time Wiener
process. The parameters $\mu$ and $s$ have the same interpretation as in
the stochastic pitchfork model: $\mu$ describes the linear stability and $s$ determines the bifurcation criticality. The additional
parameter $\gamma$ describes the strength of the spatial operator relative
to the noise and nonlinear time scales.

The effect of $\gamma$ can be seen from the dispersion relation of the
linearized reduced model,
\begin{align}
\lambda(k)
=
\mu-\gamma(1-k^2)^2.
\label{eq:SH_dispersion}
\end{align}
Through the spatial scaling $X=k_\ast x$, the physical critical wavenumber
$k_\ast$ is mapped to the rescaled wavenumber $k=1$, where
$\lambda(1)=\mu$. Hence, for fixed $\mu$, changing $\gamma$ does not
change the stability of the critical mode. Instead, it determines how
strongly the modes around $k=1$ are damped. Larger $\gamma$ gives a
narrower band of weakly damped modes, while smaller $\gamma$ allows a
broader range of spatial modes to fluctuate, as illustrated in
Figure~\ref{fig:SH_gamma_spectra}.

Expanding the spatial operator in Equation~\eqref{eq:reduced_sh} gives
\begin{align}
dq
=
\left[(\mu-\gamma)q
-2\gamma\nabla_X^2q
-\gamma\nabla_X^4q
+sq^3\right]d\tau
+dW_\tau.
\label{eq:SH_expanded}
\end{align}
The corresponding coefficients of the linear spatial terms are therefore
\begin{align}
a_0=\mu-\gamma,
\qquad
a_2=-2\gamma,
\qquad
a_4=-\gamma,
\label{eq:SH_coefficients}
\end{align}
while the cubic coefficient remains $s$. These coefficients provide the
theoretical values against which the estimates from the regression
are compared.

\subsection{Experimental setup}
\label{sec:SH_setup}

We generate numerical data from the reduced Swift--Hohenberg model in
Equation~\eqref{eq:reduced_sh}. The simulations are performed on a
one-dimensional periodic domain in the reduced spatial coordinate $X$, the
spatial derivatives are evaluated spectrally, and the system is integrated
in time using the Euler--Maruyama method. Unit space--time white noise is
discretized such that the stochastic increment at each space--time grid
point is $\sqrt{\Delta\tau/\Delta X}\,\xi_{j,n}$, where
$\xi_{j,n}\sim\mathcal{N}(0,1)$ are independent.

We perform two main experiments. First, we vary the observation horizon
$\tau_{\max}$ at fixed $\gamma=1$ to study how the available observation
time affects the ability to distinguish the two bifurcation types. Second,
we vary $\gamma$ at fixed $\tau_{\max}=2$ to study how the spatial dynamics
affect this distinction. The parameter values and numerical settings used in
the two experiments are summarized in Table~\ref{tab:SH_parameters}. In
addition, Appendix~\ref{app:SH_spatial_scaling} investigates the spatial
effect associated with changes in the unscaled critical wavenumber $k_\ast$.
This is done by varying the reduced domain length and grid spacing
according to the transformation $X=k_\ast x$, while keeping the reduced
dynamical parameters and the number of spatial observations fixed.

\begin{table}[tbp]
\centering
\small
\caption{
Parameters used in the two stochastic Swift--Hohenberg experiments.
Both experiments use a one-dimensional periodic domain of length
$L=20\pi$ with grid spacing $\Delta X=0.2\pi$, initial
condition $q(X,0)=0$, unit spatially and temporally white noise, and
$N_{\mathrm{sim}}=100$ independent realizations for each parameter
combination. Spatial derivatives are evaluated spectrally and the
simulations are integrated using the Euler--Maruyama method. The
observation-horizon experiment uses
$\Delta\tau=2.5\times10^{-4}$, while the $\gamma$-sweep experiment
uses $\Delta\tau=5\times10^{-5}$ to retain numerical stability
for the largest values of $\gamma$. The value $\gamma=0$ is included as a
test case without the spatial operator, for which Eq.~\eqref{eq:reduced_sh}
reduces to uncoupled pitchfork dynamics $dq=(\mu q+sq^3)\,d\tau+dW_\tau$ at
each spatial location, driven by independent noise, i.e., a spatial version of
the pitchfork model.
}
\begin{tabularx}{\textwidth}{lccc}
\hline
Parameter
& Description
& Varying $\tau_{\max}$
& Varying $\gamma$\\
\hline
$\mu$
& scaled stability
& $20$ log-spaced in $[-10^2,-10^{-1}]$
& $20$ log-spaced in $[-10^2,-10^{-1}]$\\

$s$
& cubic coefficient sign
& $\{-1,1\}$
& $\{-1,1\}$\\

$\gamma$
& spatial strength
& $1$
& $\{0,0.1,0.2,0.5,1,2,5,10,20,50\}$\\

$\tau_{\max}$
& observation horizon
& $\{0.1,0.2,0.5,1,2,5,10,20\}$
& $2$\\
\hline
\end{tabularx}
\label{tab:SH_parameters}
\end{table}
For the observation-horizon experiment, the different observation horizons
are compared using the same underlying stochastic realizations. Each
realization is simulated up to the largest required value of
$\tau_{\max}$, with the shorter datasets obtained by truncating the
realization at the corresponding observation horizon, as before
(Section~\ref{sec:PF_setup}).

In the subcritical case, some realizations eventually leave the local stable regime described by the cubic Swift--Hohenberg model. A realization is
excluded from the analysis if the simulated field becomes non-finite or if
\begin{align}
\max_{X,\tau}|q(X,\tau)|>10^{10}.
\label{eq:SH_escape_criterion}
\end{align}
As in the pitchfork case (Section~\ref{sec:PF_setup}), an observation horizon is included in the
analysis only if at least $10\%$ of the realizations generated from the
first $100$ random seeds are not excluded by this criterion before reaching
that horizon. For each included observation horizon, an ensemble of $100$
realizations is selected from a larger finite pool, excluding those that
leave the local regime before the end of the observation window.

In the second experiment, we vary $\gamma$ to investigate the role of the
spatial dynamics in the fitting procedure. Since $\gamma$ controls the
range of spatial modes that contribute strongly to the fluctuations,
this experiment tests how the available spatial information affects the
recovery of the nonlinear coefficient and the predicted bifurcation type.

Each realization is analyzed using the spatial fitting procedure described
in Section~\ref{sec:method_pde}. For each observation window
$\mathcal{W}_{\ell}$, the spatio-temporal data are used to estimate
\begin{align}
\hat{\Theta}_{\ell}
=
\begin{bmatrix}
\hat{a}_{0,\ell} &
\hat{a}_{2,\ell} &
\hat{a}_{4,\ell} &
\hat{s}_{\ell}
\end{bmatrix}.
\end{align}
The linear coefficients define the estimated dispersion relation
\begin{align}
\hat{\lambda}_{\ell}(k)
=
\hat{a}_{0,\ell}
-
\hat{a}_{2,\ell}k^2
+
\hat{a}_{4,\ell}k^4,
\label{eq:SH_estimated_dispersion}
\end{align}
while $\hat{s}_{\ell}$ is used to predict the bifurcation type. The ROC
curves and AUC values are constructed from the distributions of
$\hat{s}_{\ell}$ in the same way as in Section~\ref{sec:PF_setup}. Since each fitted dataset is considered separately in the following results,
we omit the observation-window index $\ell$ from the fitted quantities for
notational simplicity.

\subsection{Results}
\label{sec:SH_results}

\begin{figure}[tbp]
    \centering
    \includegraphics[width=0.9\textwidth]{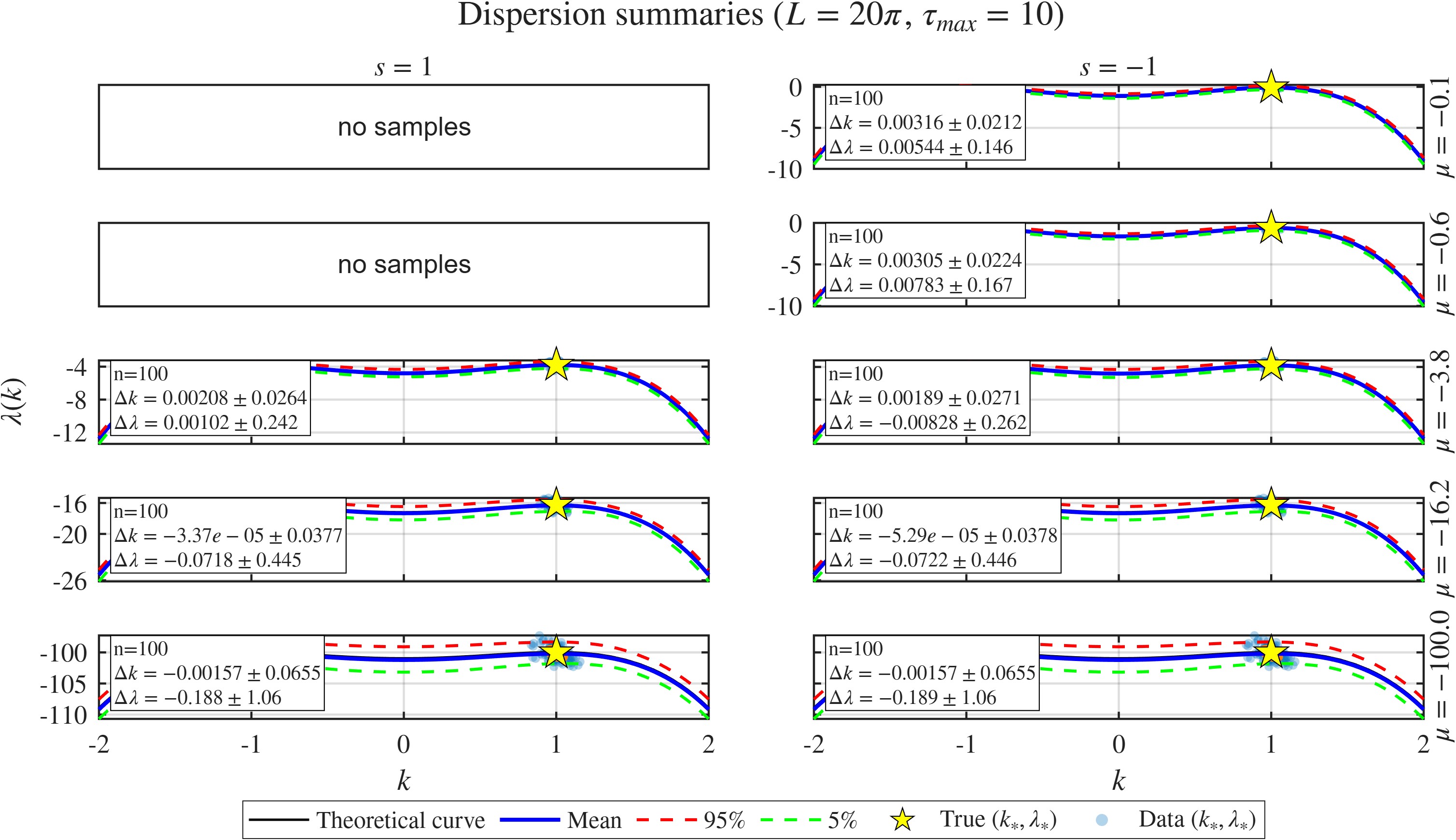}
    \caption{
Representative fitted dispersion relations for the stochastic Swift--Hohenberg model at $\gamma=1$ and observation horizon $\tau_{\max}=10$. Each panel corresponds to a fixed combination of $\mu$ and $s$. The black curve shows the theoretical dispersion relation, the blue curve the ensemble mean of the fitted relations, and the dashed curves the 5th and 95th percentiles. Blue points denote the fitted maxima $(\hat{k}_*,\hat{\lambda}_*)$ of the individual realizations, while the yellow star denotes the theoretical maximum $(k_*,\lambda_*)=(1,\mu)$. The text box reports the number of available realizations and summary statistics for $\Delta k=\hat{k}_*-k_*$ and $\Delta\lambda=\hat{\lambda}_*-\lambda_*$, following \cite{SANDERS2026135300}.
}
\label{fig:SH_dispersion_representative}
\end{figure}

The distributions of the fitted coefficients for representative parameter
settings are shown in Appendix~\ref{app:SH_ridgefits}. As in the pitchfork
model, the linear part of the dynamics is generally recovered well, with
improved estimates for longer observation horizons and closer to the
bifurcation point. Since the linear spatial dynamics are described by
multiple fitted coefficients, we additionally show a representative fitted
dispersion relation in Figure~\ref{fig:SH_dispersion_representative}. The
fitted dispersion relation closely follows the theoretical curve and
reproduces both the dominant wavenumber $k_*$ and its corresponding growth
rate $\lambda_*$. Additional fitted dispersion relations for varying
$\gamma$, together with the errors in the estimated maxima, are shown in
Appendix~\ref{app:SH_dispersion}, indicating that the results of \cite{SANDERS2026135300} remain accurate when the regression is extended to include the cubic term. We therefore focus below on the additional question of when the fitted
nonlinear coefficient contains sufficient information to predict the
bifurcation criticality, which we assess from the distribution of $\hat{s}$ using
ROC curves and AUC values.

\subsubsection*{Effect of observation horizon}

\begin{figure}[tbp]
\centering
\includegraphics[width=0.8\textwidth]{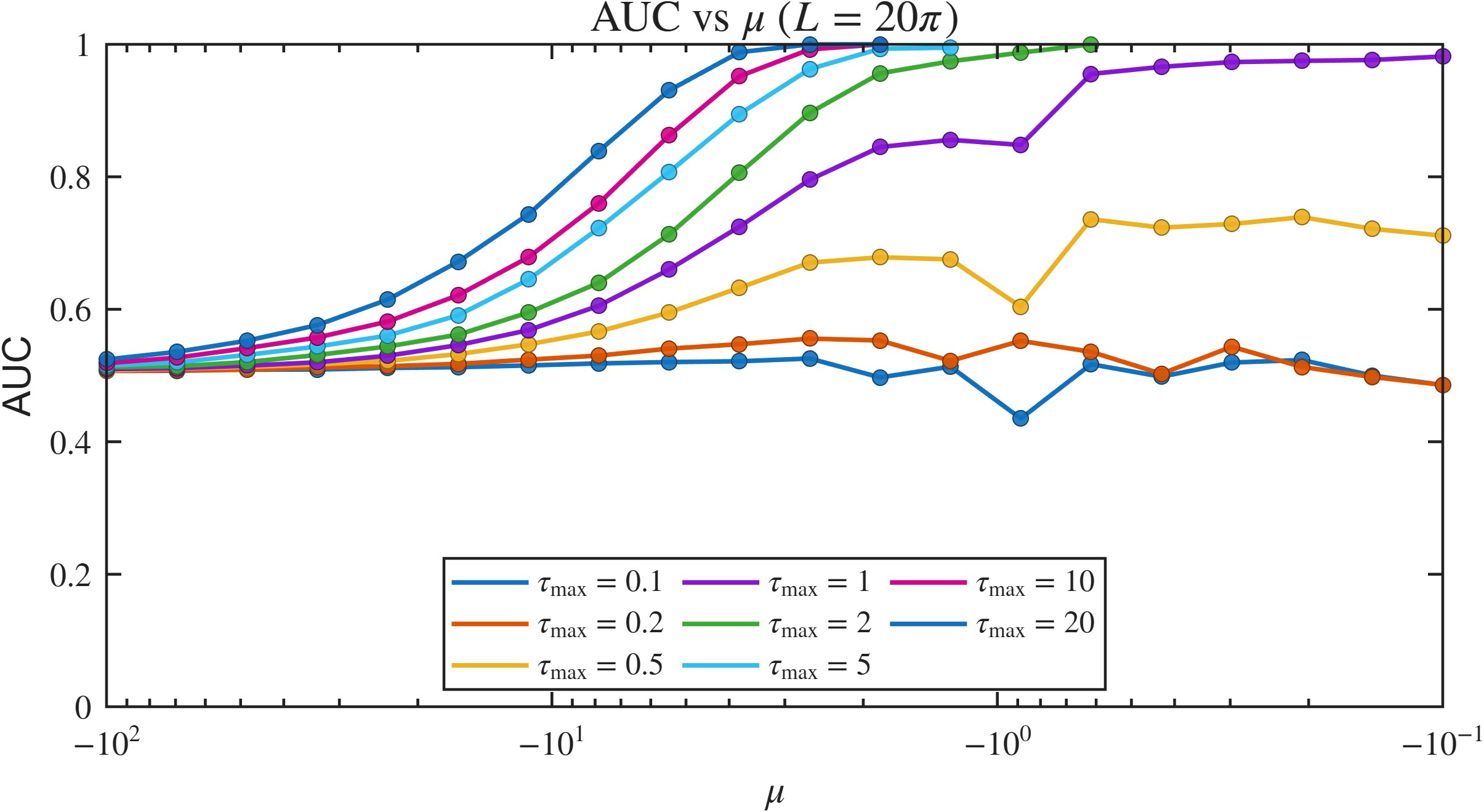}
\caption{
Area under the ROC curve (AUC) as a function of $\mu$ for different observation horizons $\tau_{\max}$ at fixed $\gamma=1$. Each curve corresponds to one observation horizon. The AUC and class conventions are as in Figure~\ref{fig:pitchfork_auc_vs_r}.
}
\label{fig:SH_auc_tau}
\end{figure}

Figure~\ref{fig:SH_auc_tau} shows the AUC for bifurcation criticality
classification as a function of the distance to the bifurcation, $\mu$, and
the observation horizon, $\tau_{\max}$, at fixed $\gamma=1$. For a fixed
observation horizon, the AUC generally increases as $\mu$ approaches zero,
showing that the nonlinear coefficient becomes easier to distinguish closer
to the bifurcation. At the same time, for a fixed value of $\mu$, longer
observation horizons lead to higher AUC values. Hence, as in the stochastic
pitchfork model, a longer observation time allows the bifurcation type to be
predicted farther away from the bifurcation point.

\begin{figure}[tbp]
\centering
\includegraphics[width=0.48\textwidth]{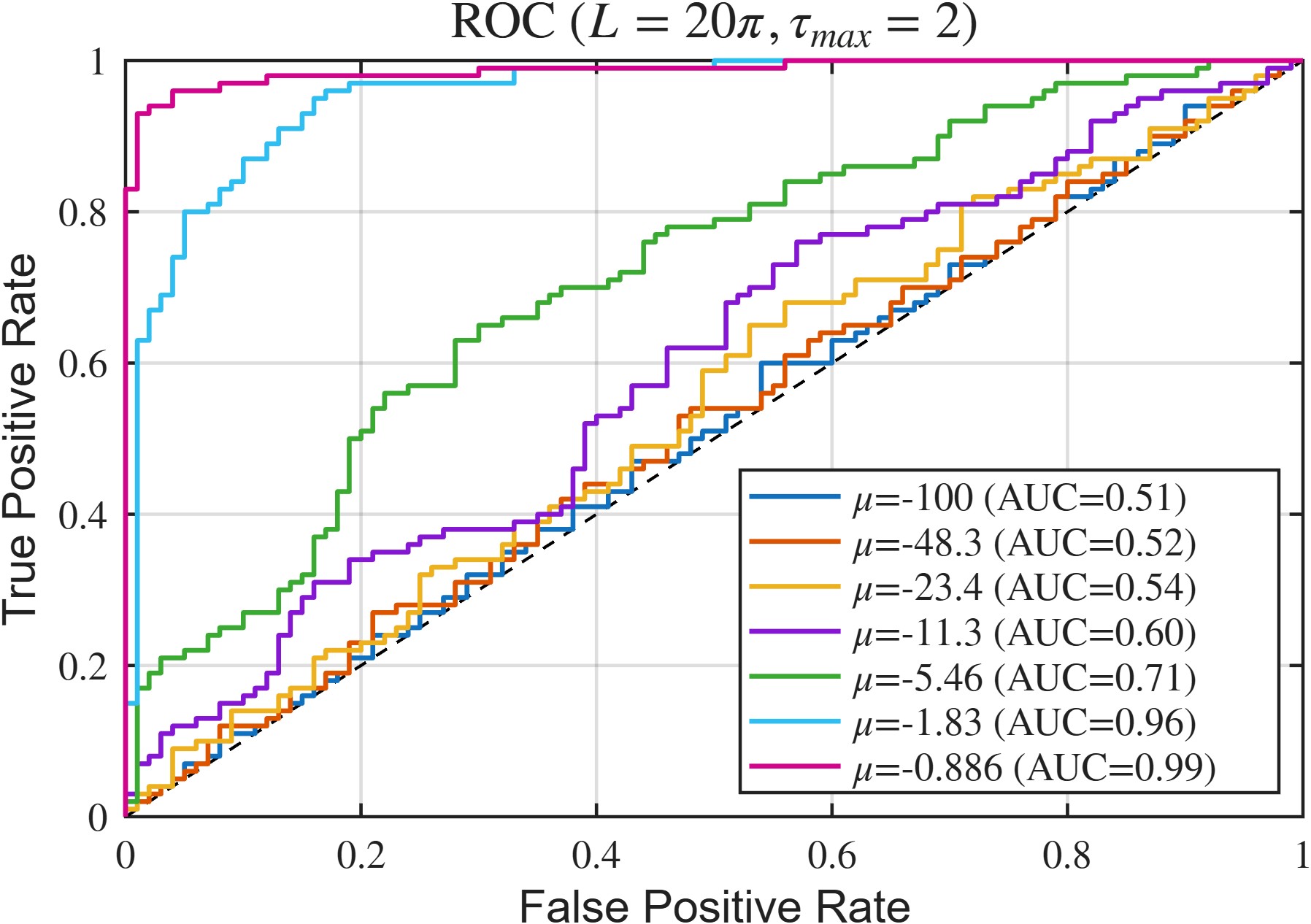}
\includegraphics[width=0.48\textwidth]{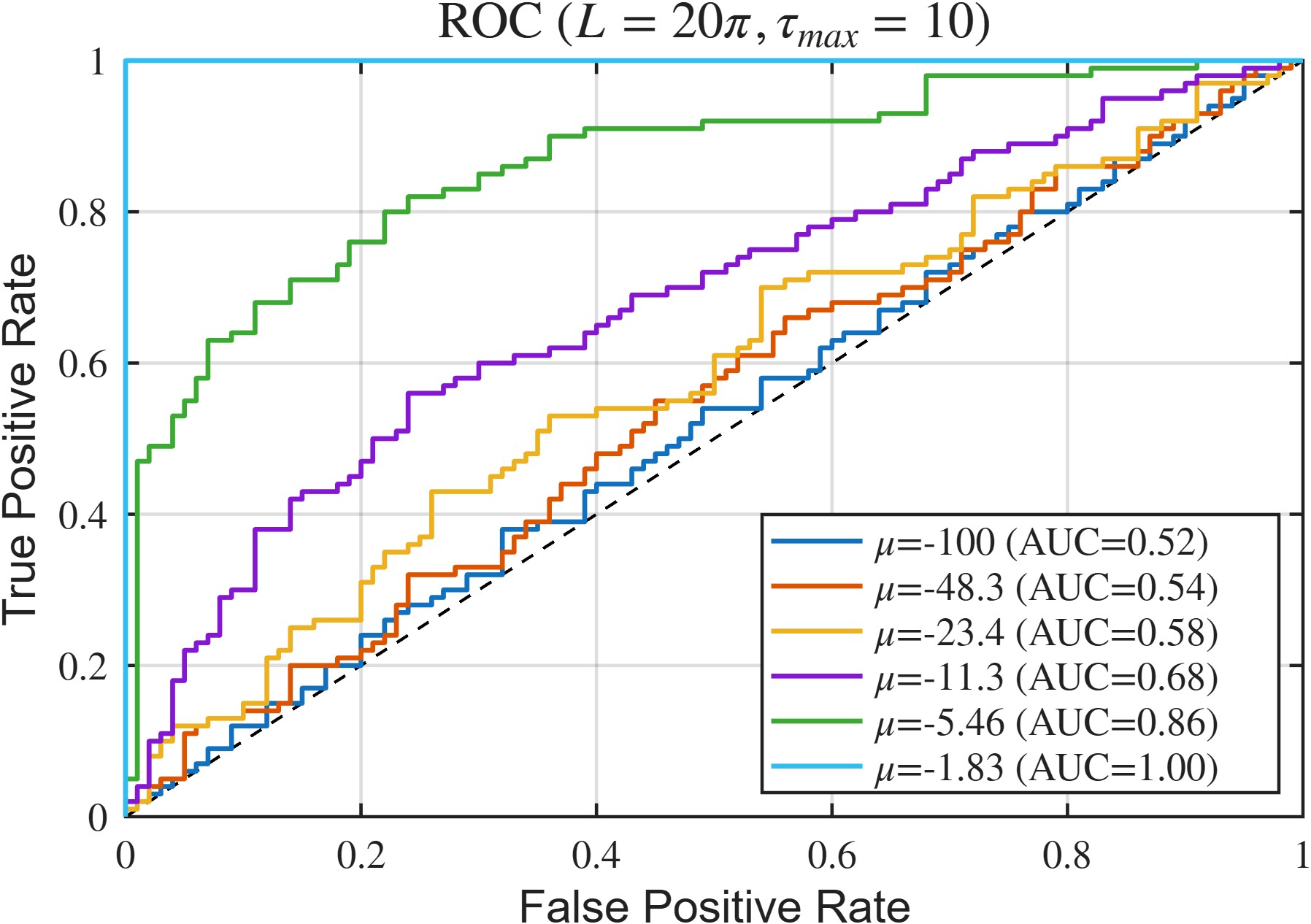}
\caption{
Representative ROC curves for the reduced stochastic Swift--Hohenberg model at $\tau_{\max}=2$ (left) and $\tau_{\max}=10$ (right), with $\gamma=1$. Each curve corresponds to a different value of $\mu$ and is constructed from the fitted values of $\hat{s}$. The axis and class conventions are as in Figure~\ref{fig:pitchfork_combined_roc}. 
}
\label{fig:SH_combined_roc_tau}
\end{figure}

This dependence on the observation horizon is also visible in the
corresponding ROC curves in Figure~\ref{fig:SH_combined_roc_tau}, shown for
two representative observation horizons ($\tau_{\max} = 2$ and $\tau_{\max} = 10$). The longer observation horizon
gives a clearer separation between the supercritical and subcritical
estimates over a larger range of $\mu$. The parameter distributions shown
in Appendix~\ref{app:SH_ridgefits} give the same picture: the distributions
of $\hat{s}$ become increasingly separated as the observation horizon
increases or $\mu$ approaches zero.

Compared with the non-spatial pitchfork case, high AUC values are obtained over much shorter
scaled observation horizons. Each temporal observation contains the complete
spatial field, so the regression uses variation in both space and time rather
than a single scalar observation at each time step. The additional spatial
information therefore allows the nonlinear dynamics to be constrained over
shorter observation horizons.

This improvement at shorter observation horizons is accompanied by a smaller
range of $\mu$ over which subcritical realizations can be included for some
of these horizons. Comparing Figures~\ref{fig:pitchfork_auc_vs_r} and
\ref{fig:SH_auc_tau}, we find that for the shared shorter observation
horizons $\tau_{\max}=2$ and $5$, the Swift--Hohenberg curves
terminate farther from the bifurcation, while for $\tau_{\max}=1$ both
models can be analyzed up to $\mu=-10^{-1}$. Thus, although accurate
classification can be achieved over shorter observation times, noise-induced tipping also seems more likely. However, for $\tau_{\max}=10$ and $20$, the 
termination points occur at comparable values of $|\mu|$ in both models.

\subsubsection*{Effect of the spatial parameter \texorpdfstring{$\gamma$}{gamma}}

\begin{figure}[tbp]
\centering
\includegraphics[width=0.8\textwidth]{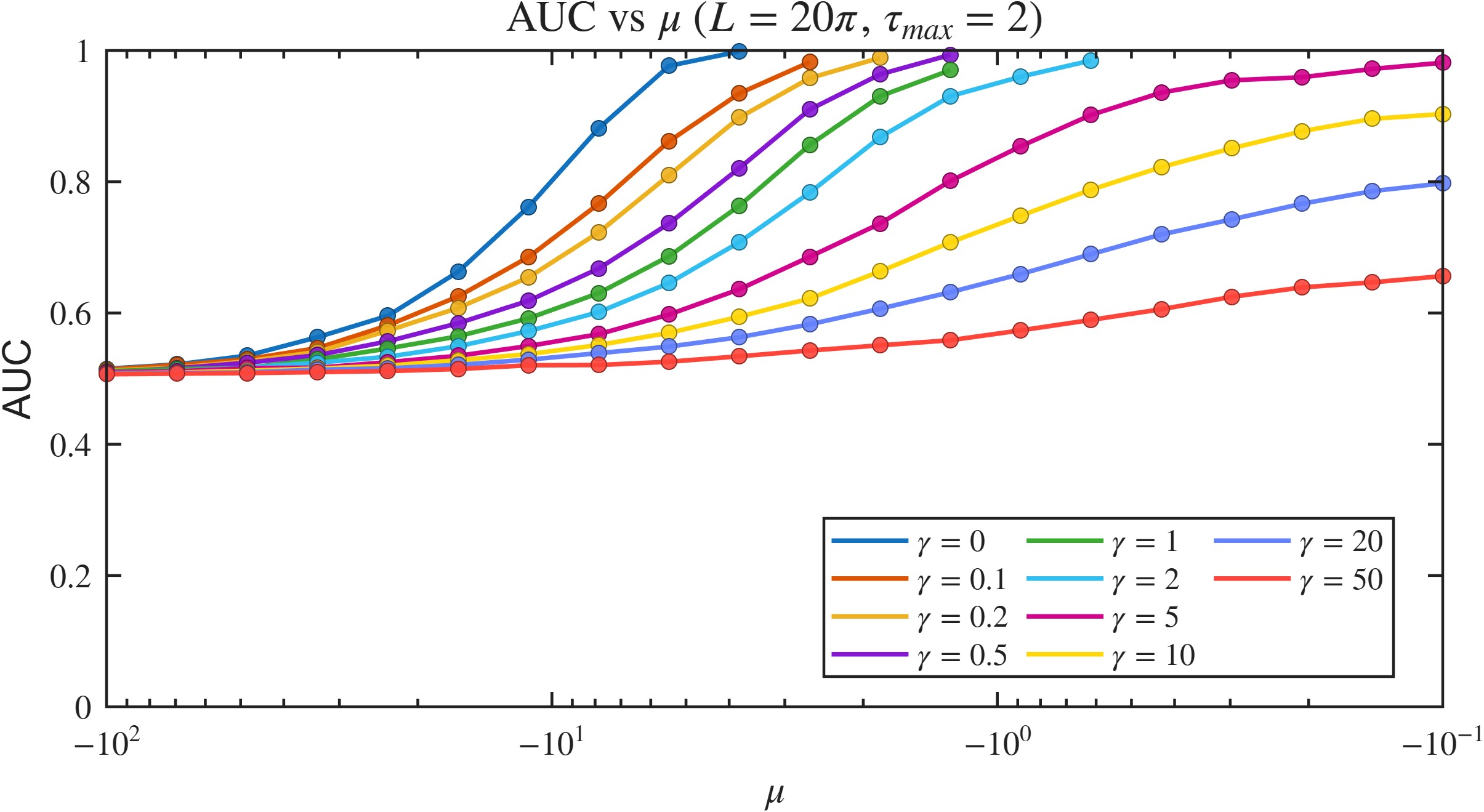}
\caption{
Area under the ROC curve (AUC) as a function of $\mu$ for different values of $\gamma$ at fixed observation horizon $\tau_{\max}=2$. Each curve corresponds to one value of $\gamma$. The AUC and class conventions are as in Figure~\ref{fig:pitchfork_auc_vs_r}. The $\gamma=1$, $\tau_{\max}=2$ curve differs slightly from the corresponding curve in the observation-horizon experiment because the two experiments use separate simulation ensembles.
}
\label{fig:SH_auc_gamma}
\end{figure}

Figure~\ref{fig:SH_auc_gamma} shows the AUC for
bifurcation criticality classification for different values of $\gamma$ at fixed
$\tau_{\max}=2$. The two criticality types are most clearly distinguished
for small $\gamma$, while the classification performance generally decreases
as $\gamma$ increases. At the same time, the curves for larger $\gamma$
extend closer to the bifurcation point. Stronger damping of the modes away
from the critical wavenumber reduces the fluctuations, so that fewer
subcritical realizations leave the local regime within the observation
horizon.

\begin{figure}[tbp]
\centering
\includegraphics[width=0.48\textwidth]{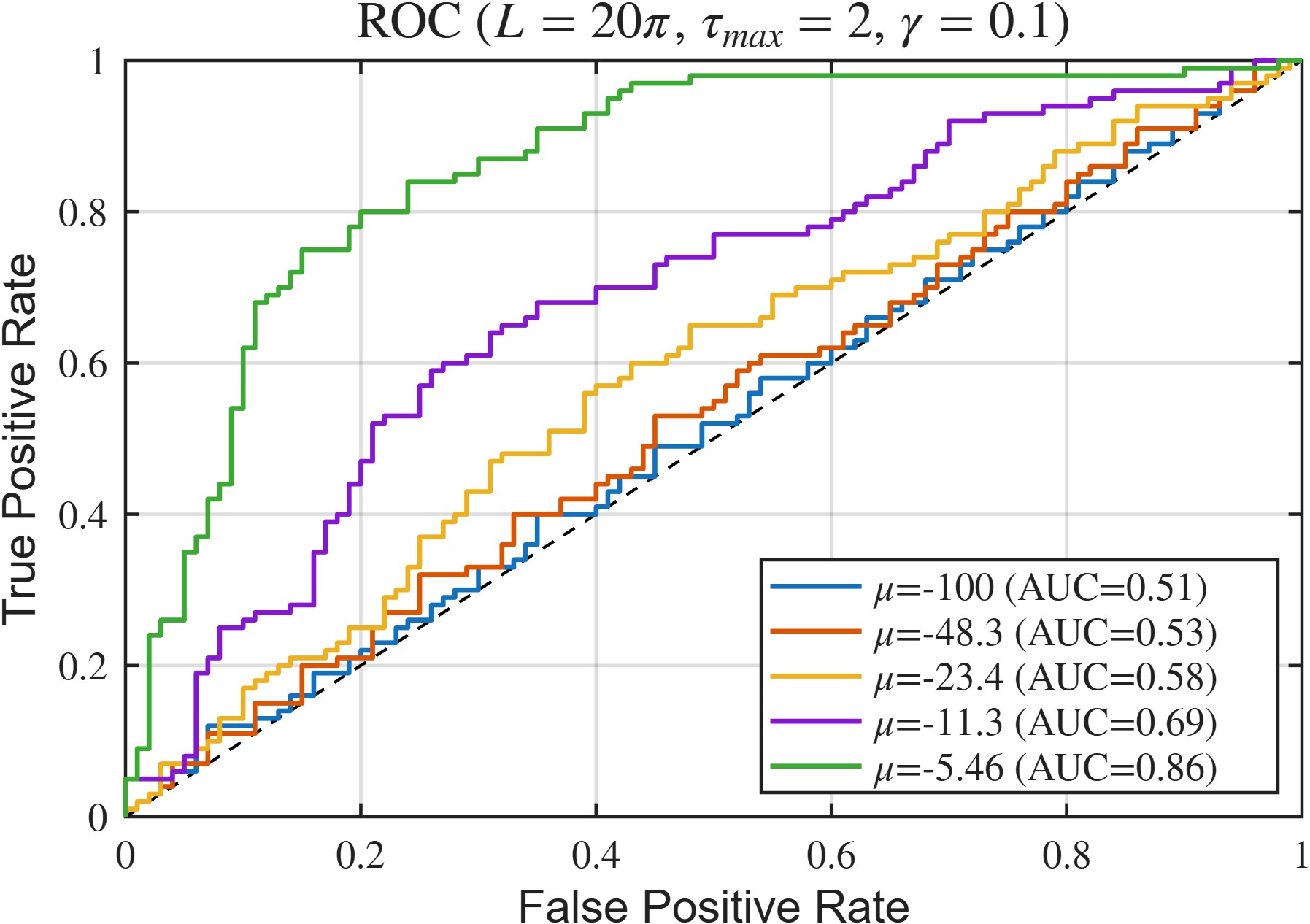}
\includegraphics[width=0.48\textwidth]{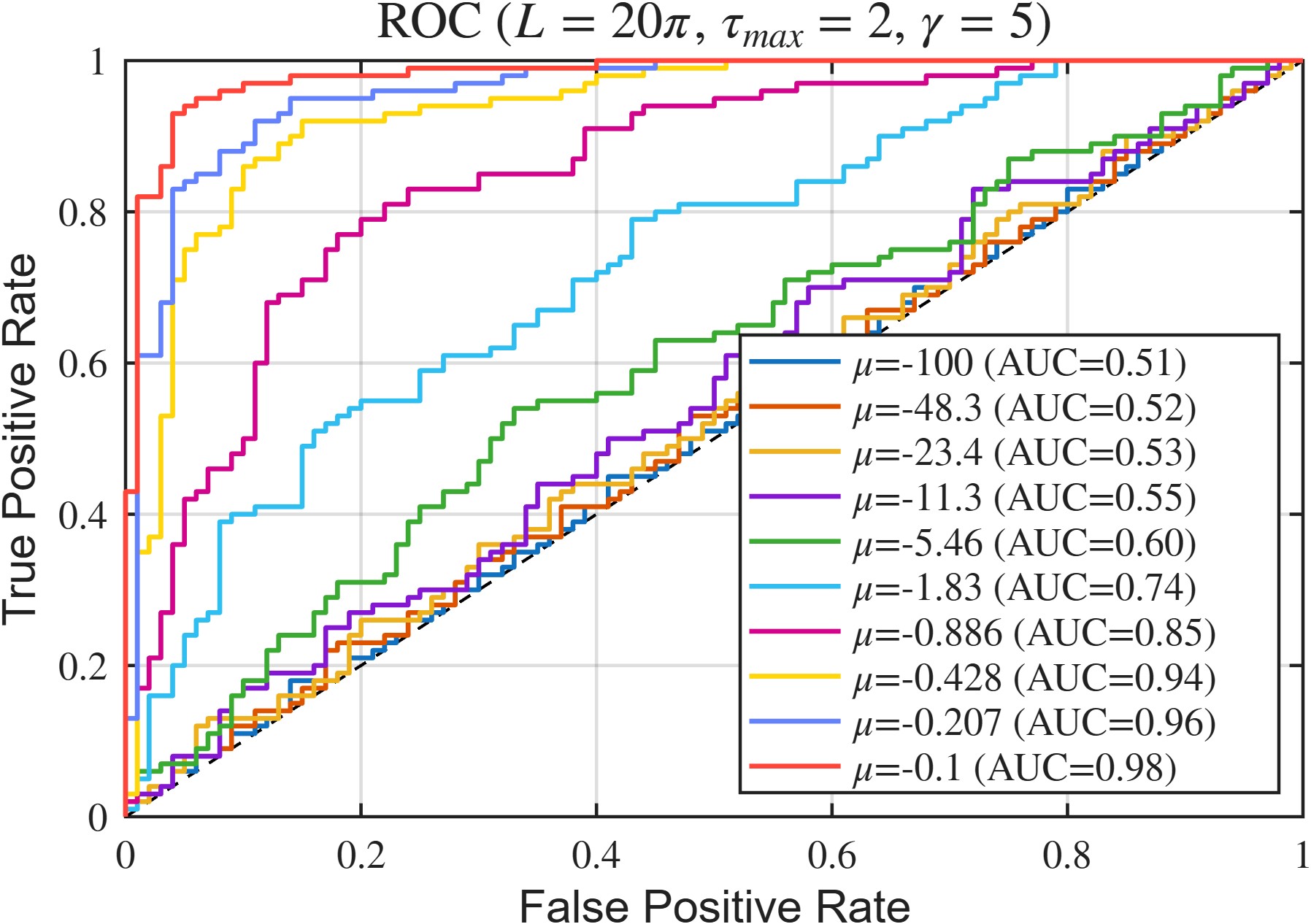}
\caption{
Representative ROC curves for $\gamma=0.1$ (left) and $\gamma=5$ (right) at fixed observation horizon $\tau_{\max}=2$. Each curve corresponds to a different value of $\mu$ and is constructed from the fitted values of $\hat{s}$. The axis and class conventions are as in Figure~\ref{fig:pitchfork_combined_roc}.
}
\label{fig:SH_combined_roc_gamma}
\end{figure}

The representative ROC curves in
Figure~\ref{fig:SH_combined_roc_gamma} illustrate the same change in
classification performance ($\gamma = 0.1$ and $\gamma = 5$). For $\gamma=0.1$, the supercritical and
subcritical estimates can be separated over a larger part of the stable
parameter range, whereas for $\gamma=5$ the ROC curves remain closer to the
diagonal. The corresponding parameter distributions are shown in
Appendix~\ref{app:SH_ridgefits}. To relate this change in classification
performance to the spatial dynamics, we also consider how the dispersion
relation and the observed spectrum change with $\gamma$.

The fitted dispersion relations in Appendix~\ref{app:SH_dispersion} show
that increasing $\gamma$ progressively narrows the band around the critical
wavenumber $k=1$. Consistently, the representative spectra in
Figure~\ref{fig:SH_gamma_spectra} show that the observed fluctuations become
concentrated over a narrower range of wavenumbers and that the effective
number of contributing modes decreases. This suggests that the change in
classification performance is related to the range of spatial modes
contributing to the observed dynamics, which we discuss further in
Section~\ref{sec:discussion_spatial}.

\section{Discussion}
\label{sec:discussion}

The results in Sections~\ref{sec:PF_results} and~\ref{sec:SH_results} show that in some cases the bifurcation criticality can be distinguished from noisy observations of a system before it crosses the associated bifurcation point or before noise-induced tipping. However, the classification performance depends
strongly on the observation horizon, the distance to the bifurcation, and,
in the spatial model, the distribution of fluctuations across spatial modes. At the same time, the linear stability properties (i.e., the closeness to the bifurcation point) remain generally
well recovered with the extended regression method introduced here. We first discuss the effects of the observation horizon and
the distance to the bifurcation, before considering the additional effects
introduced by spatial observations.

\subsection{Required observation time}
\label{sec:discussion_temporal}

Bifurcation criticality classification was better for longer observation times $\tau_{\max}$ and for $\mu$ closer to the bifurcation point (Figures~\ref{fig:pitchfork_auc_vs_r} and~\ref{fig:SH_auc_tau}); this was also observed in the direct parameter estimates (see Appendices~\ref{app:pitchfork} and~\ref{app:SH_ridgefits}). Longer
observation times provide more information for estimating the (cubic)
coefficient, while closer to the bifurcation the weaker linear restoring
dynamics allow the nonlinear contribution to become more pronounced in the
observed fluctuations. Longer observation horizons therefore not only
improve the classification close to the bifurcation, but also allow reliable
classification substantially farther away from the transition. In contrast, for very short observation times the two bifurcation criticality types may not be distinguishable at all; for example, the shortest observation
horizons in Figure~\ref{fig:SH_auc_tau} give AUC values close to $0.5$ over
most of the investigated range.

This distinction is particularly relevant when considering the method as an early-warning indicator. In the subcritical case, near-complete separation is generally reached only close to the bifurcation, where stochastic fluctuations also make noise-induced tipping more likely and limit how long trajectories can be observed. This trade-off is consistent with the first-passage-time analysis in Figure~\ref{fig:mfpt_kramers} and Appendix~\ref{app:hitting_time}, where the expected escape time decreases rapidly as $\mu$ approaches zero. Longer observation records are therefore valuable because they allow reliable classification at larger $|\mu|$, where noise-induced tipping is less likely.

In the experiments we have used rescaled parameters. It is therefore informative to understand these results in the context of the unscaled parameters. For the pitchfork model,
\begin{align}
\tau_{\max}^{\mathrm{PF}}
&=
t_{\max}\sigma\sqrt{|\alpha|},
&
|\mu^{\mathrm{PF}}|
&=
\frac{|r|}{\sigma\sqrt{|\alpha|}},
\end{align}
and for the Swift--Hohenberg model
\begin{align}
\tau_{\max}^{\mathrm{SH}}
&=
t_{\max}\sigma\sqrt{k_\ast|\alpha|},
&
|\mu^{\mathrm{SH}}|
&=
\frac{|r|}{\sigma\sqrt{k_\ast|\alpha|}}.
\end{align}
For fixed $k_\ast$, both models have the same qualitative dependence
on observation time, noise strength, nonlinearity, and distance from the
bifurcation. In particular, increasing $t_{\max}$ naturally increases the available scaled
observation horizon. Increasing either $\sigma$ or $|\alpha|$
simultaneously increases the scaled observation horizon and decreases
$|\mu|$. Both changes correspond to regimes in which the bifurcation type is more
easily distinguished in our numerical results. This is consistent with the
increase in fluctuation amplitude under stronger stochastic forcing or weaker
restoring dynamics near a critical transition
\cite{DitlevsenJohnsen2010,Kuehn2013}. Larger fluctuations increase the
relative contribution of the nonlinear term, making it easier to identify
in the model fit. However, in the subcritical case, these fluctuations also increase the possibility of noise-induced tipping, thus resulting in the trade-off between identifiability and escape
discussed above.

\subsection{Using spatial information}
\label{sec:discussion_spatial}

In the spatial Swift--Hohenberg equation, each observation contained data across a spatial field. Although these observations are spatially correlated, they provide variation across space in addition to that over time, which can help constrain the linear and nonlinear coefficients in the regression. Comparing the non-spatial pitchfork model (Figure~\ref{fig:pitchfork_auc_vs_r}) and the spatial Swift--Hohenberg equation (Figure~\ref{fig:SH_auc_tau}), this is reflected in the Swift--Hohenberg model reaching high classification performance at much shorter observation horizons. However, the spatial setting also appears more prone to noise-induced tipping than the non-spatial model, particularly close to the bifurcation point. One possible reason is that farther from the bifurcation, stronger restoring dynamics make sufficiently large fluctuations increasingly unlikely, so the presence of many possible escape locations has only a limited effect. This is consistent with the first-passage-time analysis in Figure~\ref{fig:mfpt_kramers}, where the expected escape time eventually increases exponentially with $\left|\mu\right|^2$. Closer to the bifurcation, the restoring dynamics weaken and the expected escape time decreases rapidly, making it more likely that a sufficiently large fluctuation somewhere in the spatial domain initiates departure from the local regime.

The strength of the spatial effects also influences the classification performance. Figure~\ref{fig:SH_auc_gamma} shows that the
bifurcation types become more difficult to distinguish as $\gamma$
increases. Increasing $\gamma$
narrows the band of weakly damped modes around the critical wavenumber. For stochastic
Swift--Hohenberg dynamics, the variance of individual Fourier modes is
directly related to their linear eigenvalues
\cite[Sec.~3.2]{GowdaKuehn2015}. Increasing $\gamma$ therefore
concentrates the observed fluctuations into a smaller range of spatial
modes, leaving less spatial variation available to distinguish the
different terms in the fitted dynamics. This provides a possible
explanation for the reduced ability to identify the nonlinear coefficient
and hence the bifurcation type.

These spatial effects can be related back to the unscaled parameters as
\begin{align}
\gamma
=
\frac{\Gamma k_\ast^4}
{\sigma\sqrt{k_\ast|\alpha|}}.
\label{eq:gamma_physical}
\end{align}
We see that for fixed $k_\ast$, an increase in either the noise strength $\sigma$ or the
magnitude of the nonlinear coefficient $|\alpha|$ leads to a decrease of
$\gamma$, an increase of $\tau_{\max}$, and a decrease of $|\mu|$. All three changes
move the system towards regimes in which the bifurcation type was easier to distinguish. In the spatial model, stronger
noise or nonlinearity therefore not only makes the nonlinear dynamics more
pronounced, but also increases the range of weakly damped spatial modes that
can contribute to distinguishing the linear and nonlinear dynamics.

The critical physical wavenumber $k_\ast$ introduces an additional spatial
dependence that cannot be determined from the separate effects of
$\mu$, $\tau_{\max}$, and $\gamma$ alone. We therefore investigated the
direct effect of increasing $k_\ast$ separately in
Appendix~\ref{app:SH_spatial_scaling}. The results show that increasing
$k_\ast$ improves the bifurcation criticality classification. A larger
$k_\ast$ corresponds to a shorter critical wavelength, so that more critical
wavelengths are contained within a fixed physical observation domain.
Equivalently, more spatial modes lie close to the critical wavenumber and
exhibit substantial variance, as also observed for the stochastic
Swift--Hohenberg equation on larger domains
\cite[Sec.~4.1]{GowdaKuehn2015}. This additional contribution from the other modes may
help constrain the fitted dynamics, and thereby improve estimation of $\hat{s}$.

\subsection{Applicability to observational data and limitations}
\label{sec:discussion_applicability}

The results of this study were obtained for controlled stochastic models,
for which the relevant linear and nonlinear terms are known and the data
can be observed at a prescribed temporal and spatial resolution. For
real systems, the underlying model structure will be
less certain. In this case, the regression
may need to include a larger library of candidate terms. Sparse
model-discovery methods, such as SINDy and its extensions, provide one
approach for identifying which terms are supported by the observations
\cite{BruntonEtAl2025,Song2024}. However, in all cases, the relevant dynamics
must be represented in the candidate library and sufficiently constrained by
the available data. Based on the findings in this paper, the main prerequisite for the use of a bifurcation criticality classification method is the availability of observations that are long enough and contain sufficiently strong fluctuations to identify the nonlinear dynamics before a transition occurs.

\section{Conclusion}
\label{sec:conclusion}

In this paper, we have introduced a bifurcation criticality classification method to distinguish imminent subcritical bifurcations from supercritical bifurcations using noisy observations. This extends linear stability estimates, such as those based on critical slowing down metrics \cite{Scheffer2009, van2007slow, wissel1984universal, dakos2011slowing}, by including the leading cubic nonlinear term in a regression framework. We further tested this method on a non-spatial pitchfork bifurcation and a spatial Swift--Hohenberg equation, to identify the situations in which it seems possible to determine a bifurcation's criticality. Our results indicate that bifurcation criticality classification is indeed possible in certain situations. In particular, the prediction improves with longer observation times, when the system is closer to the bifurcation, and in systems with stronger noise or more dominant nonlinear effects. However, in the subcritical case, stronger noise or nonlinearity also increases the risk of an earlier noise-induced transition, which limits the window of opportunity for criticality classification. Further, the inclusion of spatial information can lead to improved predictions, thus decreasing the observation time needed for good criticality prediction, although the precise usefulness depends on the range of spatial modes contributing to the observations, with broader spectral fluctuations providing better separation of the fitted nonlinear dynamics.

These results show that noisy observations of a system before a bifurcation do contain information not only about an approaching
instability, but also about the nonlinear character of the transition that
follows. Combining estimates of the linear and nonlinear dynamics
therefore provides a route towards early-warning systems that can predict
both the approach to a transition and its expected criticality. Extending this
framework to other bifurcation mechanisms, more complex systems, and
observational data provides a natural next step towards predicting the
character of critical transitions in real systems.

\section*{Acknowledgments}
B.H.\ was funded by the Deutsche Forschungsgemeinschaft (DFG, German Research
Foundation) -- Project number 543917644. P.A.S.\ thanks the Centre for Complex
Systems Studies (CCSS) at Utrecht University for supporting a research visit
to the Technical University of Munich through the Swaantje Mondt Fund.

\section*{Author Declarations}

\subsection*{Conflict of Interest}
The authors have no conflicts to disclose.

\subsection*{Use of generative AI}
During the preparation of this work, the authors used ChatGPT (OpenAI), Claude (Anthropic)
and Gemini (Google) to assist with
writing and improving MATLAB code for the simulations and analysis, to brainstorm ideas,
and to improve the clarity of the text and the presentation of the results. All
AI-assisted code was checked and tested by the authors, and all text and ideas were
critically reviewed and revised by the authors, who take full responsibility for the
content of this publication. No generative AI was used to create or alter research data
or figures.

\subsection*{Author Contributions}
\textbf{Paul A. Sanders:} Conceptualization (lead); Methodology (lead);
Software (lead); Formal analysis (lead); Investigation (lead);
Visualization (lead); Writing -- original draft (lead);
Writing -- review \& editing (equal).
\textbf{Christian Kuehn:} Conceptualization (supporting); Methodology
(supporting); Supervision (equal);
Writing -- review \& editing (supporting).
\textbf{Bastian Hilder:} Conceptualization (supporting); Methodology
(supporting); Supervision (equal); 
Writing -- review \& editing (supporting).
\textbf{Robbin Bastiaansen:} Conceptualization (supporting); Methodology
(supporting); Formal analysis (supporting); Supervision (equal);
Writing -- review \& editing (equal).

\section*{Data Availability}
The MATLAB code used to generate and analyze the data in this study, together
with the processed results and the figures of the paper runs, is openly
available on GitHub at
\url{https://github.com/Paul-A-Sanders/Criticality_classification}. The raw simulated
trajectories are not stored but can be regenerated with this code from the
recorded random seeds.

\clearpage
\bibliographystyle{unsrturl}
\bibliography{References}

\clearpage
\appendix
\numberwithin{figure}{section}
\numberwithin{table}{section}

\section{Mean first passage time for the subcritical pitchfork model}\label{app:hitting_time}

This appendix quantifies the time scale on which trajectories of the subcritical pitchfork model can leave the neighborhood of the stable equilibrium before the deterministic bifurcation is reached. This provides a reference for the noise-induced tipping discussed in Section~\ref{sec:PF_setup} and Figure~\ref{fig:mfpt_kramers}. We follow the standard backward-equation approach for first-passage problems \cite{Gardiner1985} and compare the resulting mean first passage time with the Eyring--Kramers approximation \cite{Kramers1940,berglund2013kramers}.

For the subcritical case $s=1$, the reduced pitchfork equation in Equation~\eqref{eq:reduced_sde} has drift
\begin{align}
f(q)=\mu q+q^3.
\label{eq:pitchfork_drift}
\end{align}
For $\mu<0$, the equilibrium $q=0$ is stable and is bounded by the two unstable equilibria $q=\pm\sqrt{-\mu}$. For an initial condition $q(0)=q$ with $|q|<\sqrt{-\mu}$, we define the first exit time from this interval as
\begin{align}
\tau_e=\inf\left\{\tau>0:\lvert q(\tau)\rvert=\sqrt{-\mu}\right\}.
\label{eq:exit_time}
\end{align}
The mean first passage time from the starting point $q$ is
\begin{align}
\mathcal{T}(q)=\mathbb{E}[\tau_e|q(0)=q].
\end{align}

The infinitesimal generator of the reduced stochastic process is
\begin{align}
\mathcal{L}\phi(q)
=
(\mu q+q^3)\phi'(q)+\frac{1}{2}\phi''(q).
\label{eq:pf_generator}
\end{align}
The mean first passage time satisfies the standard backward equation \cite{Gardiner1985}
\begin{align}
f(q)\mathcal{T}'(q)+\frac{1}{2}\mathcal{T}''(q)=-1,
\qquad
\mathcal{T}\!\left(-\sqrt{-\mu}\right)
=
\mathcal{T}\!\left(\sqrt{-\mu}\right)
=
0.
\label{eq:mfpt_ode}
\end{align}
The boundary conditions represent absorption once either unstable equilibrium is reached. Since $\mathcal{T}(q)$ is a function of the initial state rather than time, Equation~\eqref{eq:mfpt_ode} is a boundary-value problem and does not require a separate initial condition. By symmetry, we can also change one of the boundary conditions to $\mathcal{T}'(0)=0$.

Writing $f(q)=-V'(q)$ with
\begin{align}
V(q)=-\frac{\mu}{2}q^2-\frac{1}{4}q^4
\label{eq:Potential_PF}
\end{align}
 and applying an integrating factor to solve the backward equation~\eqref{eq:mfpt_ode} gives the exact integral representation from our equilibrium at $q=0$
\begin{align}
\mathcal{T}(0)
=
2\int_0^a e^{2V(z)}
\left(\int_0^z e^{-2V(y)}\,dy\right)dz,
\qquad
a=\sqrt{-\mu}.
\label{eq:T0_final}
\end{align}

For $|\mu|\gg 1$, the Eyring--Kramers approximation gives
\cite{Kramers1940,berglund2013kramers}
\begin{align}
\mathcal{T}_{\mathrm{K}}(0)
\sim
\frac{\pi}{\sqrt{2}|\mu|}
\exp\!\left(\frac{\mu^2}{2}\right).
\label{eq:mfpt_kramers}
\end{align}
The comparison in Figure~\ref{fig:mfpt_kramers} uses $C\mathcal{T}_{\mathrm{K}}(0)$ with $C=0.6$. This empirical factor improves the prefactor over the finite parameter range shown in the figure.

\subsection{Gamma approximation of the exit-time distribution}

\begin{figure}[tbp]
\centering
\includegraphics[width=0.48\textwidth]{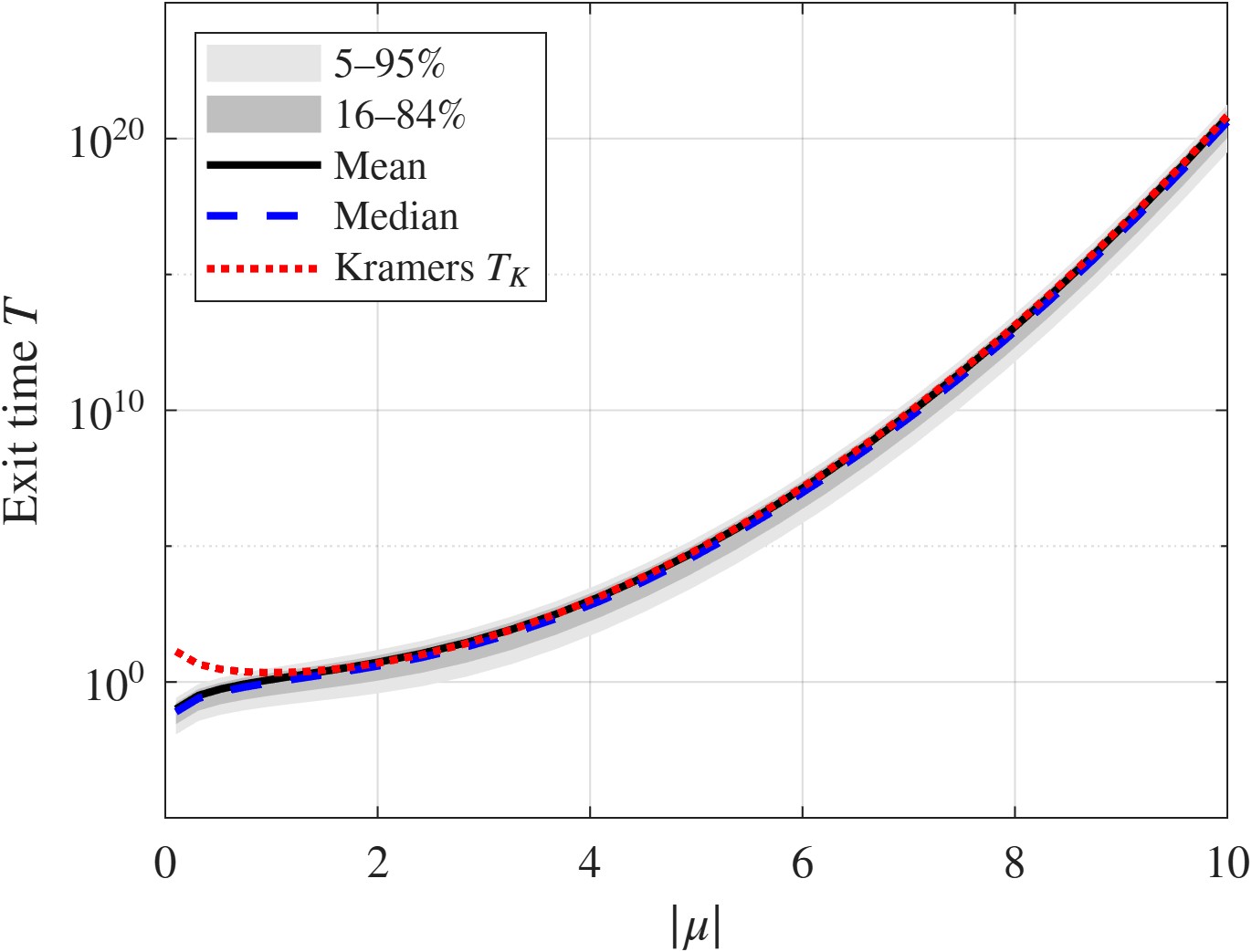}
\includegraphics[width=0.48\textwidth]{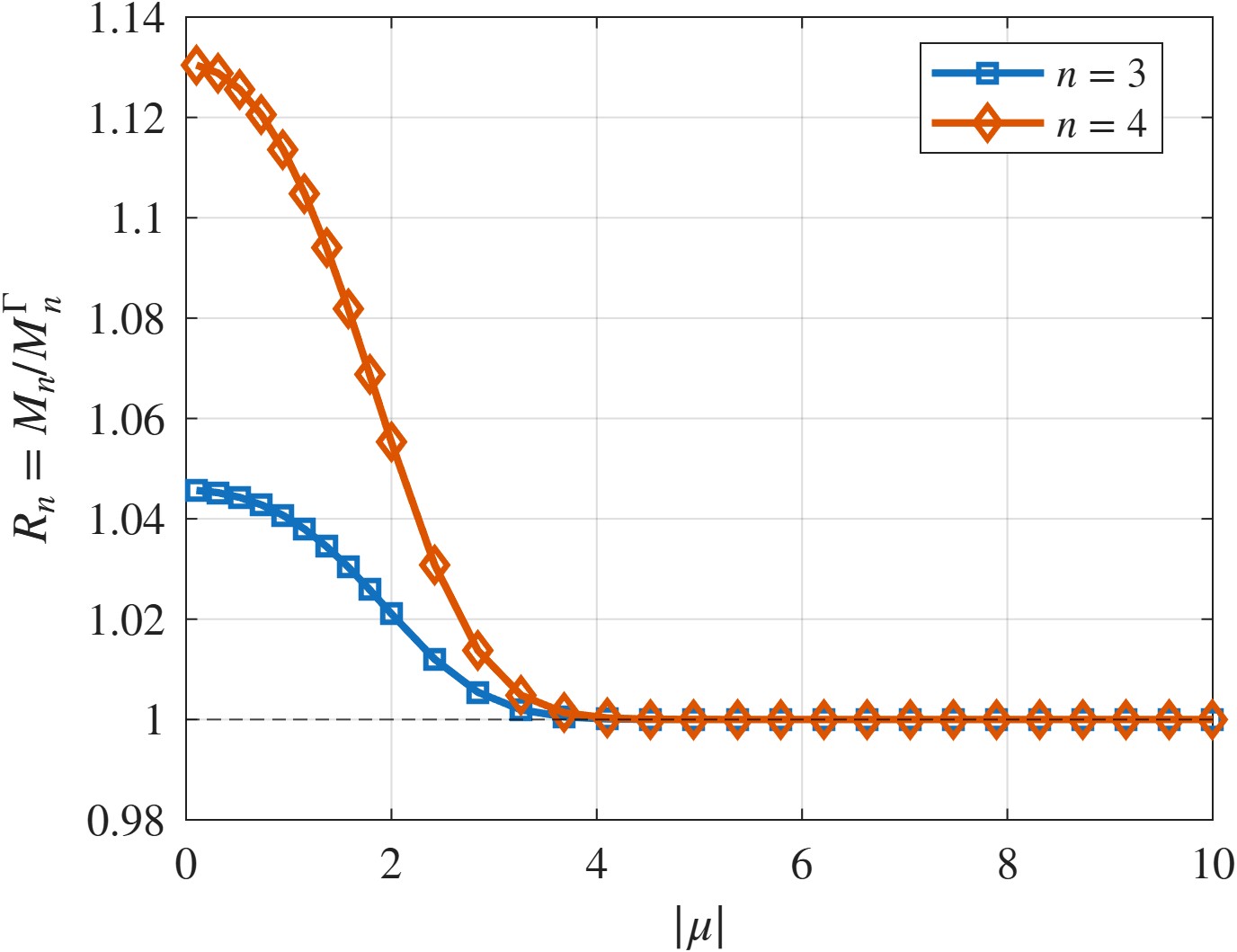}
\caption{
(a) Mean first passage time and quantile intervals obtained from the two-moment Gamma approximation, shown together with the Eyring--Kramers scaling. (b) Ratios between the third and fourth moments of the exit time and the corresponding moments of the fitted Gamma distribution; the reference value $R_n=1$ corresponds to equal moments.
}
\label{fig:mfpt_gamma}
\end{figure}

The moments $M_n(q)=\mathbb{E}_q[\tau_e^n]$ satisfy the backward hierarchy \cite{Gardiner1985}
\begin{align}
\mathcal{L}M_n(q)=-nM_{n-1}(q),
\qquad
M_n(\pm a)=0,
\end{align}
with $M_0(q)=1$. We approximate the exit-time distribution by a Gamma distribution, choosing its shape and scale parameters to match the mean and variance of the exit time at $q=0$. This follows the standard moment-matching procedure for the Gamma distribution \cite{Siegrist2017}. Similar Gamma approximations have been used for stochastic exit-time distributions \cite{Simpson2021Gamma}.
\begin{align}
m&=M_1(0),
&
v&=M_2(0)-M_1^2(0),
\\
\kappa&=\frac{m^2}{v},
&
\theta&=\frac{v}{m}.
\end{align}
The corresponding quantiles of the exit-time distribution are then approximated by
\begin{align}
Q_p \approx Q_p^\Gamma = F_\Gamma^{-1}(p;\kappa,\theta),
\end{align}
where $F_\Gamma^{-1}$ is the inverse cumulative distribution function of the fitted Gamma distribution.
To test this approximation, we compare the higher moments with the Gamma moments and its ratio
\begin{align}
M_n^{\Gamma}=\theta^n\frac{\Gamma(\kappa+n)}{\Gamma(\kappa)},
\qquad
R_n=\frac{M_n(0)}{M_n^{\Gamma}}.
\end{align}
Values of $R_n$ close to one indicate that the two-moment Gamma approximation also describes the higher moments. In Figure \ref{fig:mfpt_gamma} we show the ratios $R_3$ and $R_4$.

\clearpage
\section{Parameter distributions for the pitchfork model}
\label{app:pitchfork}
To provide additional insight into the fitting procedure and the ROC analysis presented in Section~\ref{sec:PF_results}, we show the distributions of the estimated linear and nonlinear coefficients obtained from the stochastic simulations. The distributions are visualized using kernel density estimates computed with the MATLAB \texttt{ksdensity} function using a Gaussian kernel. The estimates of $\hat{\mu}$ and $\hat{s}$ are shown separately to make the accuracy and spread of both fitted coefficients easier to assess.

\begin{figure}[tbp]
    \centering
    \includegraphics[width=0.8\textwidth]{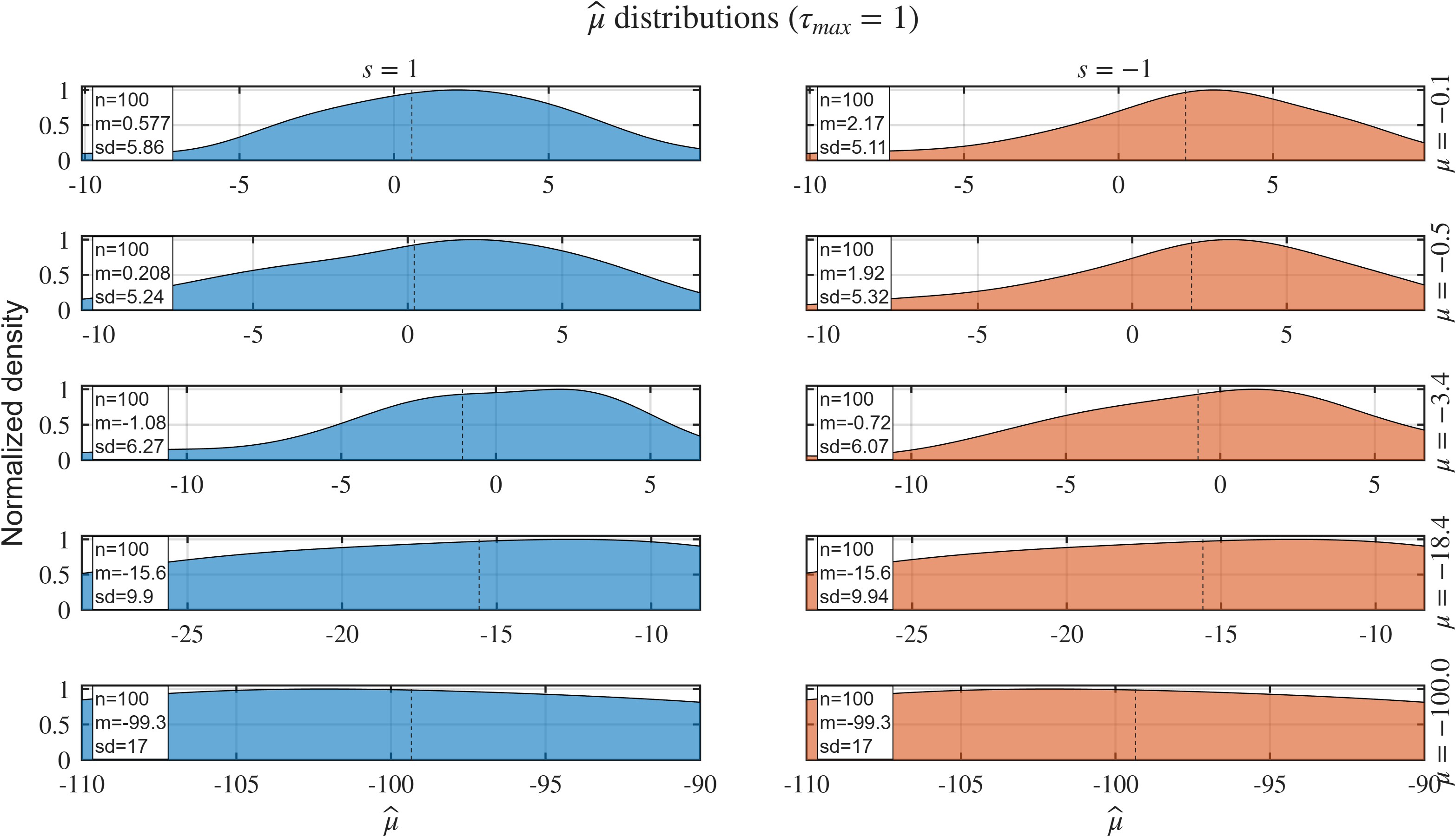}
    \caption{
Distributions of the estimated linear coefficient $\hat{\mu}$ for an observation horizon $\tau_{\max}=1$. Each panel corresponds to a fixed combination of the true parameters $\mu$ and $s$, based on 100 independent stochastic realizations when sufficiently many realizations can be included in the analysis. The filled curves show the estimated probability densities, with the fill colors distinguishing the supercritical ($s=-1$) and subcritical ($s=1$) models. The black line outlines each density and the vertical dashed line marks the sample mean of the fitted estimates. The text box gives the number of available realizations, the sample mean, and the sample standard deviation. Parameter combinations for which too few realizations can be included in the analysis are left empty.
}
    \label{fig:PC_ridge_mu_t=1}
\end{figure}

\begin{figure}[tbp]
    \centering
    \includegraphics[width=0.8\textwidth]{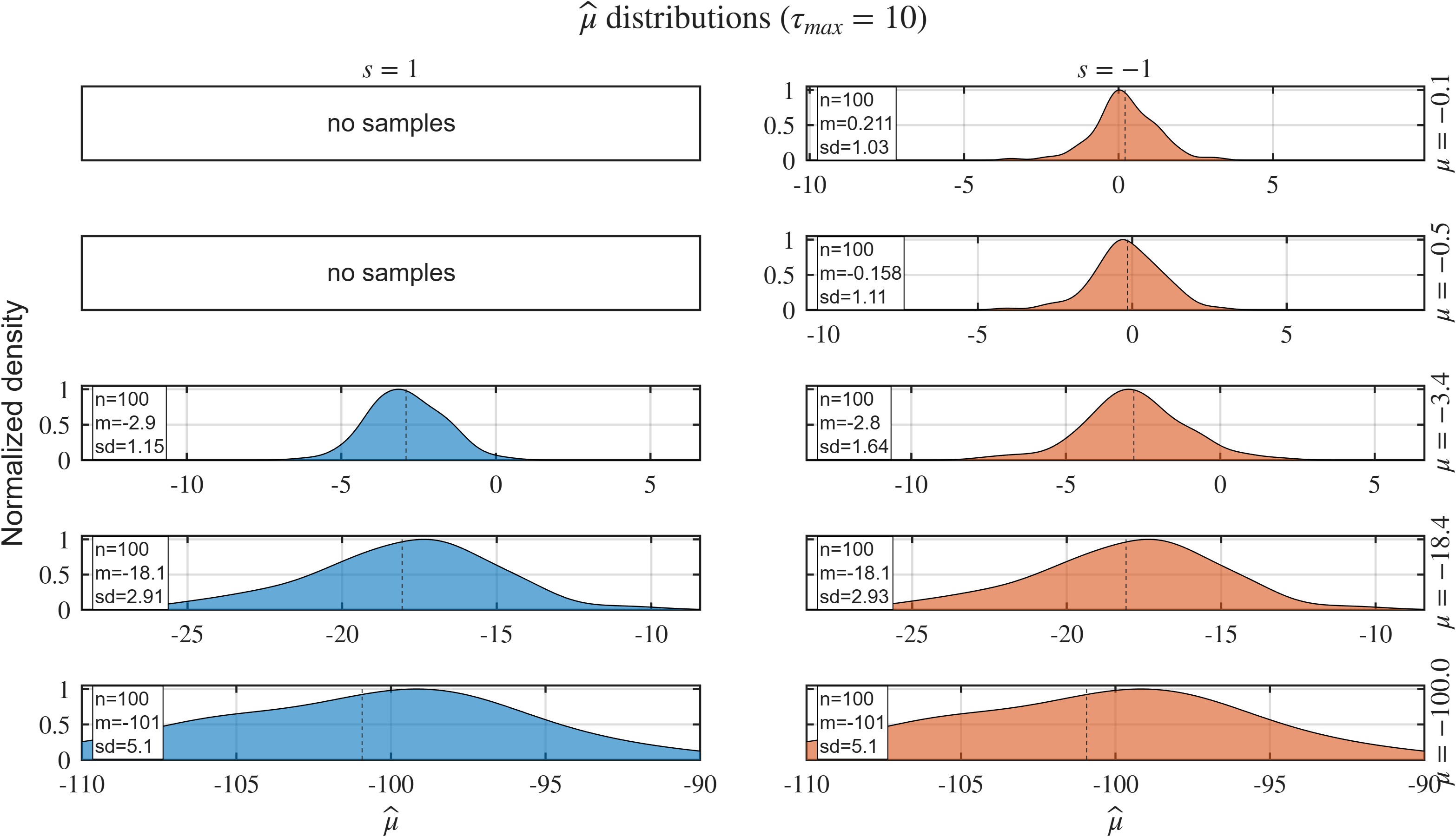}
    \caption{
Distributions of the estimated linear coefficient $\hat{\mu}$ for $\tau_{\max}=10$. Plotting conventions are as in Figure~\ref{fig:PC_ridge_mu_t=1}.
}
    \label{fig:PC_ridge_mu_t=10}
\end{figure}

\begin{figure}[tbp]
    \centering
    \includegraphics[width=0.8\textwidth]{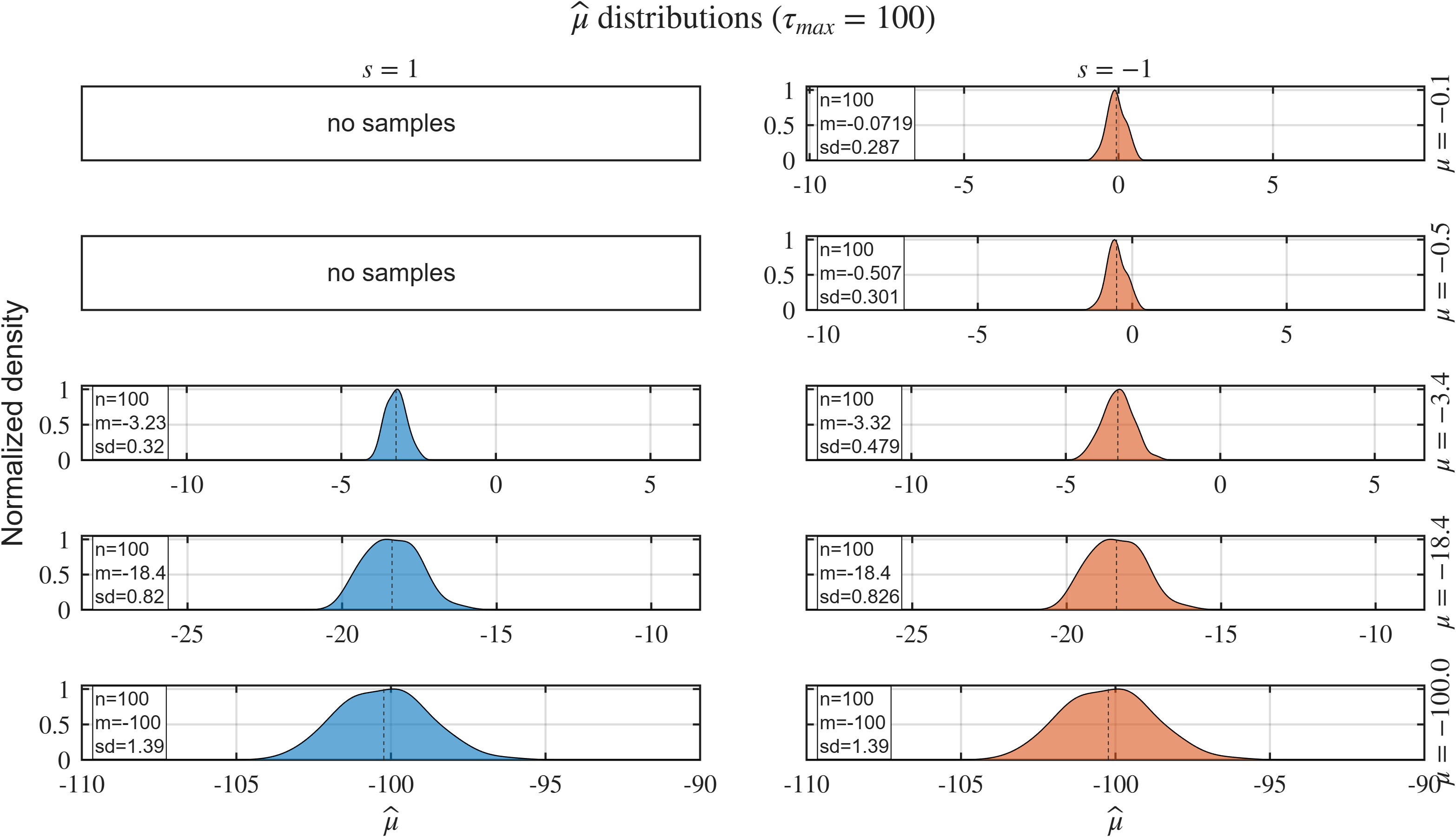}
    \caption{
Distributions of the estimated linear coefficient $\hat{\mu}$ for $\tau_{\max}=100$. Plotting conventions are as in Figure~\ref{fig:PC_ridge_mu_t=1}.
}
    \label{fig:PC_ridge_mu_t=100}
\end{figure}

\begin{figure}[tbp]
    \centering
    \includegraphics[width=0.8\textwidth]{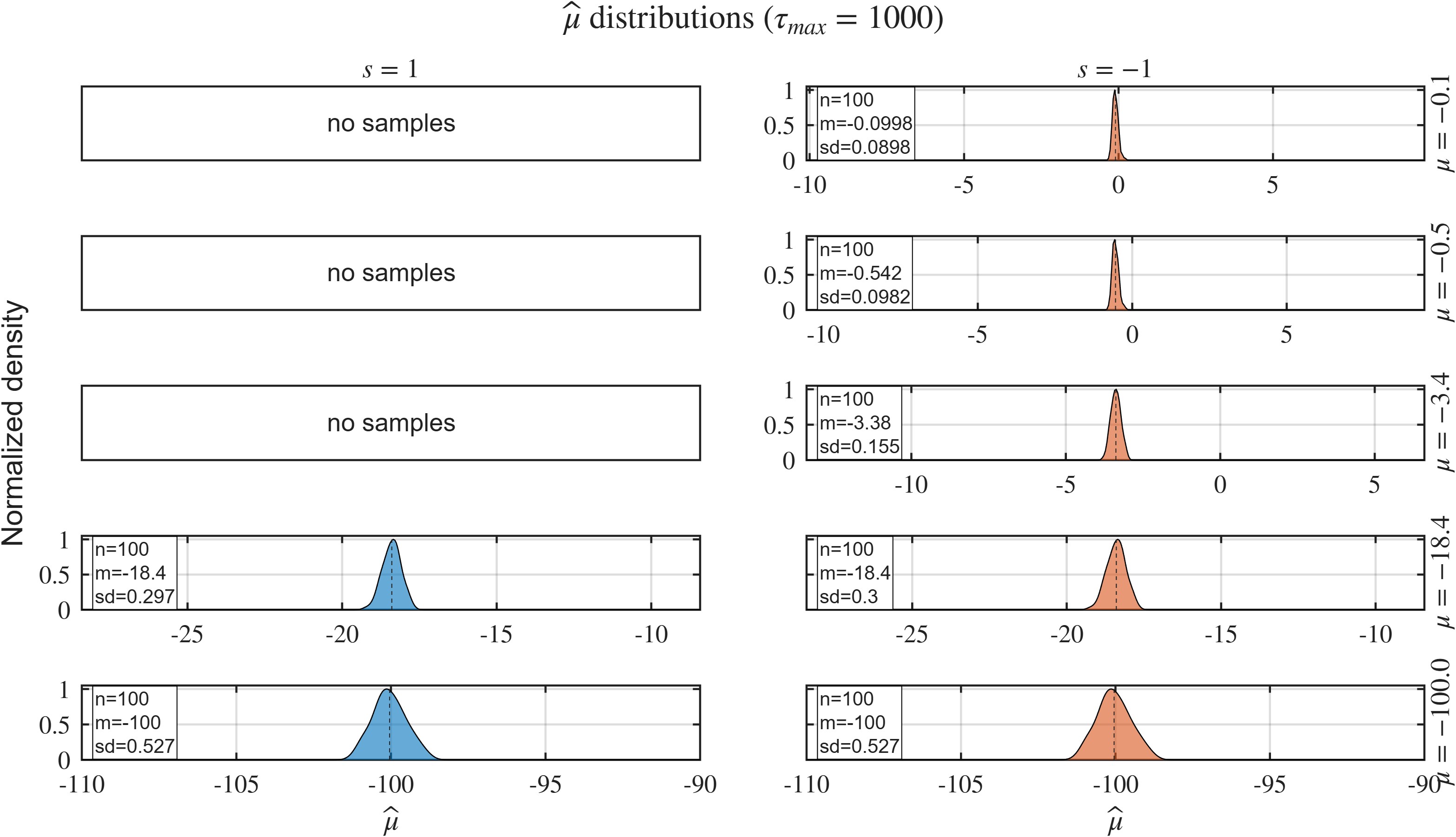}
    \caption{
Distributions of the estimated linear coefficient $\hat{\mu}$ for $\tau_{\max}=1000$. Plotting conventions are as in Figure~\ref{fig:PC_ridge_mu_t=1}.
}
    \label{fig:PC_ridge_mu_t=1000}
\end{figure}

The distributions of $\hat{\mu}$ for representative observation horizons are shown in Figures~\ref{fig:PC_ridge_mu_t=1}--\ref{fig:PC_ridge_mu_t=1000}. In general, the linear coefficient is recovered more accurately for longer observation horizons and closer to the bifurcation point.

\begin{figure}[tbp]
\centering
\includegraphics[width=0.8\textwidth]{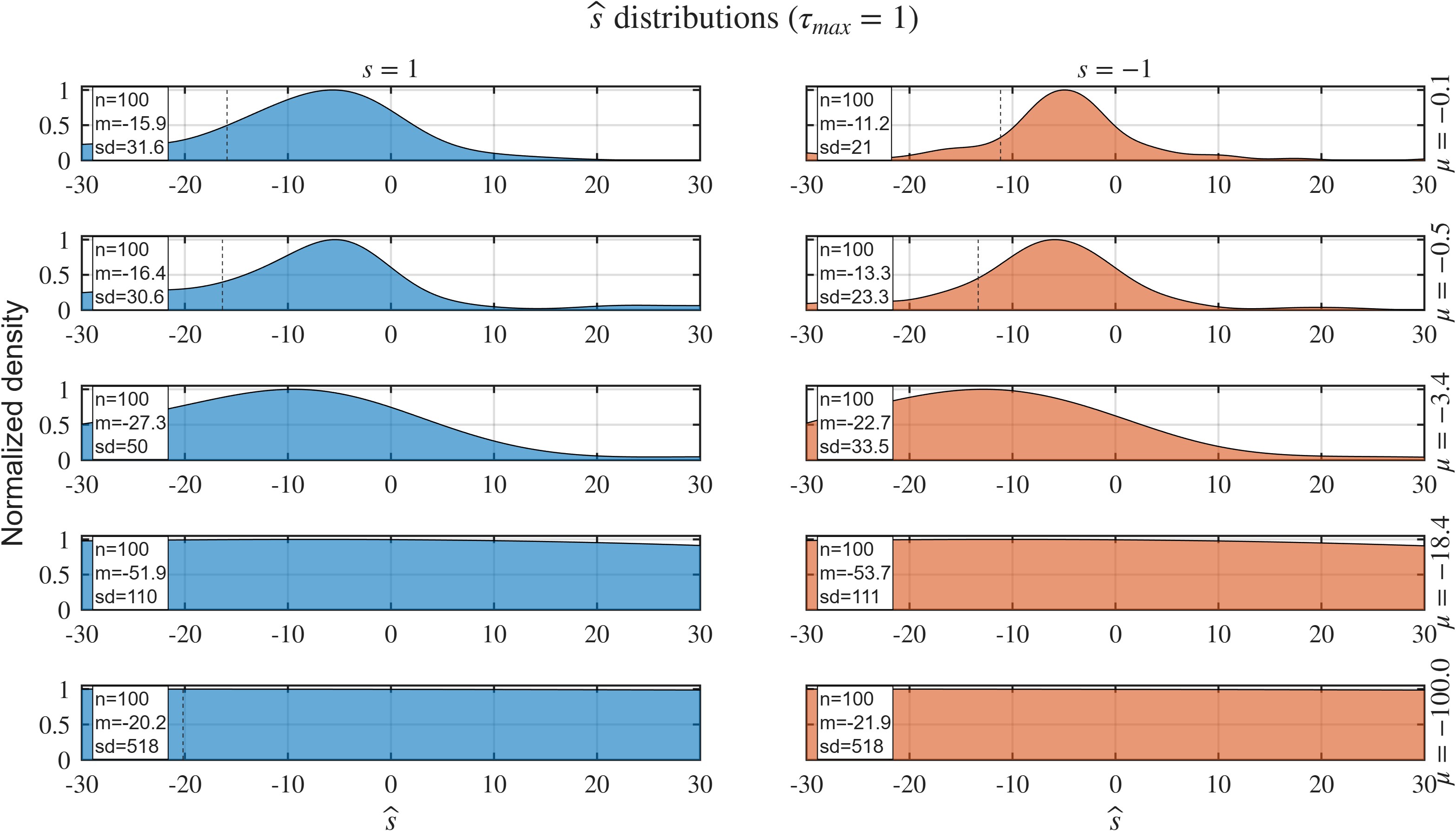}
    \caption{
Distributions of the estimated nonlinear coefficient $\hat{s}$ for $\tau_{\max}=1$. The true values $s=-1$ and $s=1$ correspond to the supercritical and subcritical models, respectively. Plotting conventions are as in Figure~\ref{fig:PC_ridge_mu_t=1}.
}
    \label{fig:PC_ridge_s_t=1}
\end{figure}

\begin{figure}[tbp]
    \centering
    \includegraphics[width=0.8\textwidth]{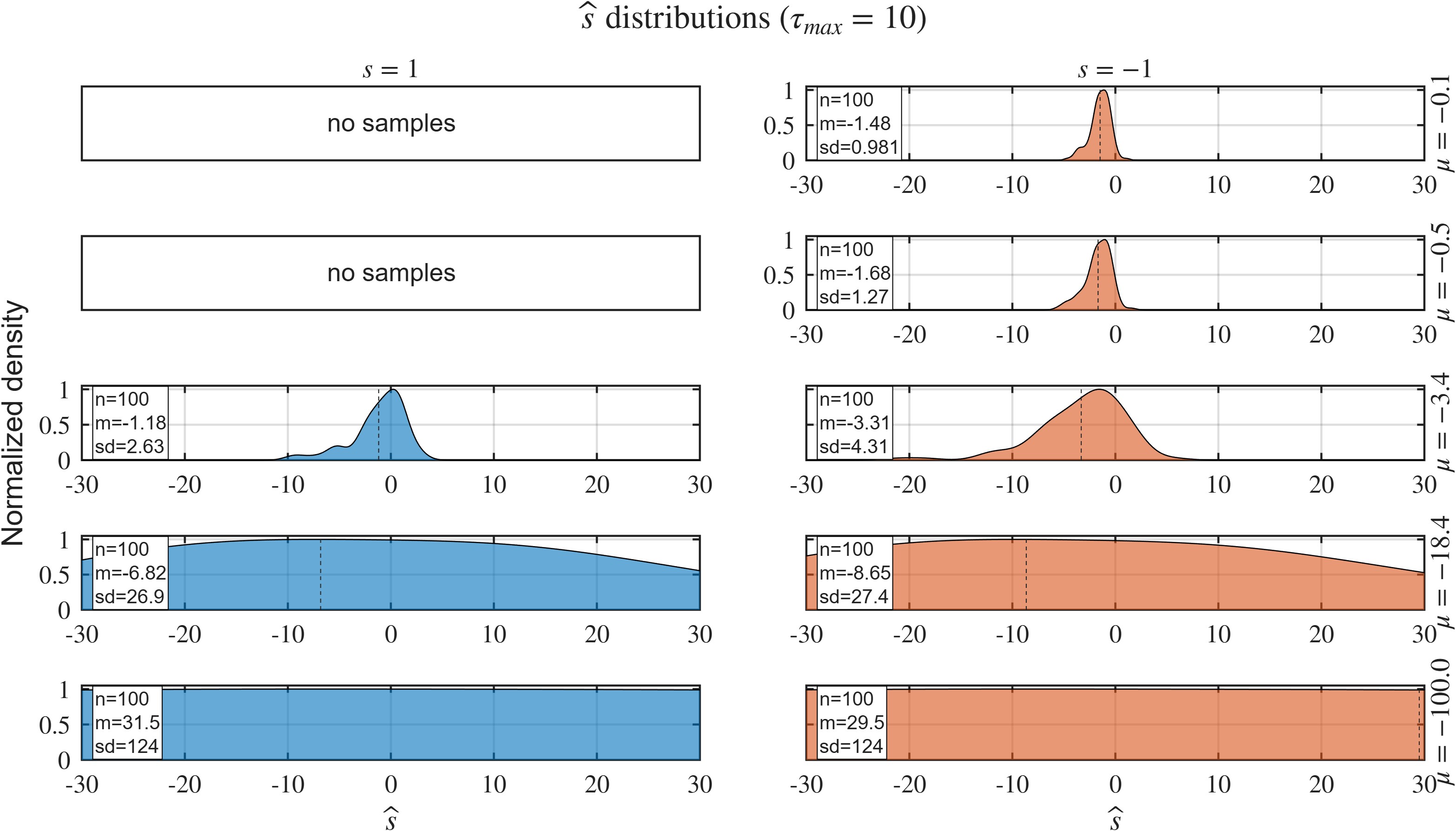}
    \caption{
Distributions of the estimated nonlinear coefficient $\hat{s}$ for $\tau_{\max}=10$. Plotting conventions are as in Figure~\ref{fig:PC_ridge_mu_t=1}.
}
    \label{fig:PC_ridge_s_t=10}
\end{figure}

\begin{figure}[tbp]
    \centering
    \includegraphics[width=0.8\textwidth]{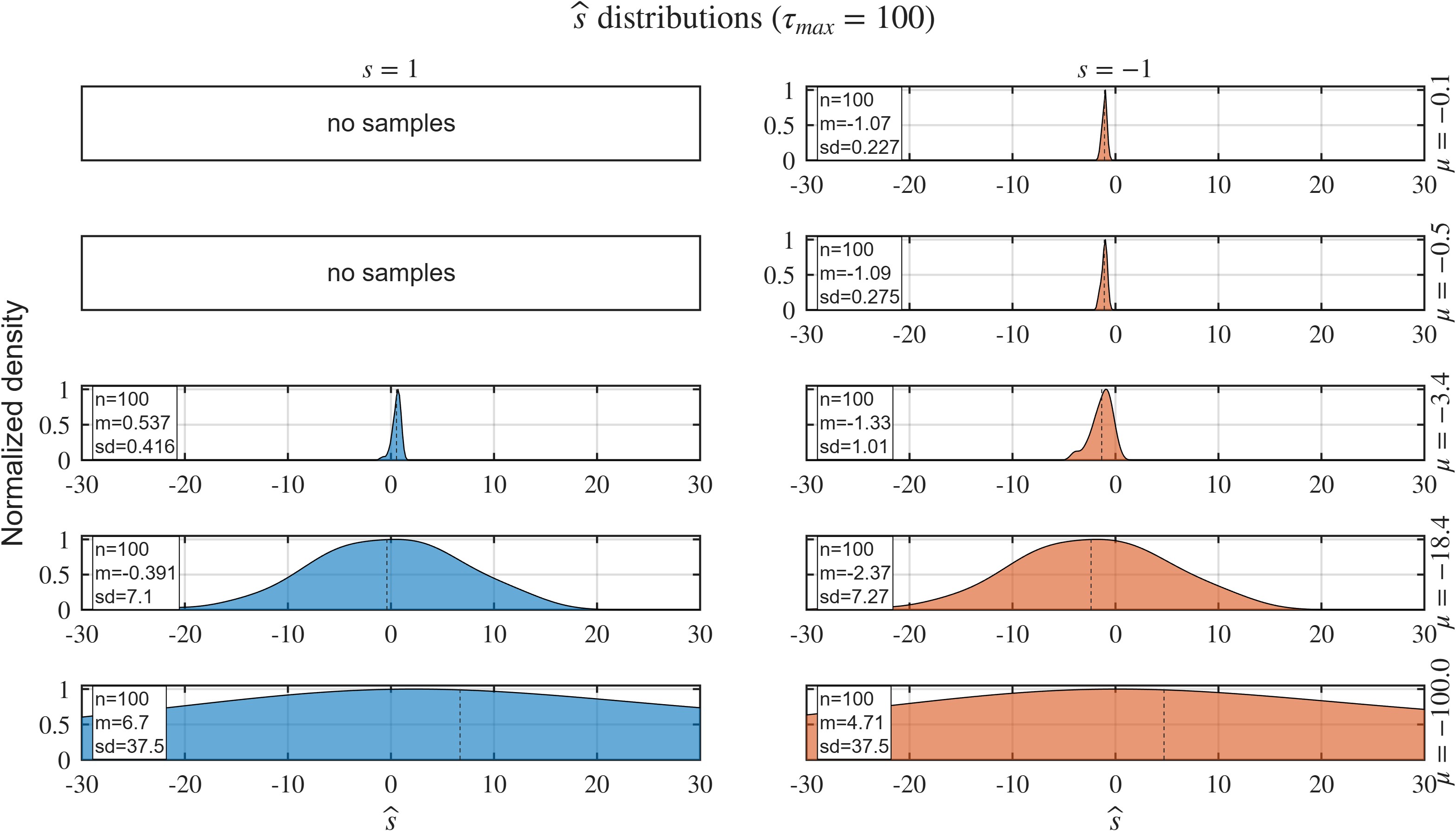}
    \caption{
Distributions of the estimated nonlinear coefficient $\hat{s}$ for $\tau_{\max}=100$. Plotting conventions are as in Figure~\ref{fig:PC_ridge_mu_t=1}.
}
    \label{fig:PC_ridge_s_t=100}
\end{figure}

\begin{figure}[tbp]
    \centering
    \includegraphics[width=0.8\textwidth]{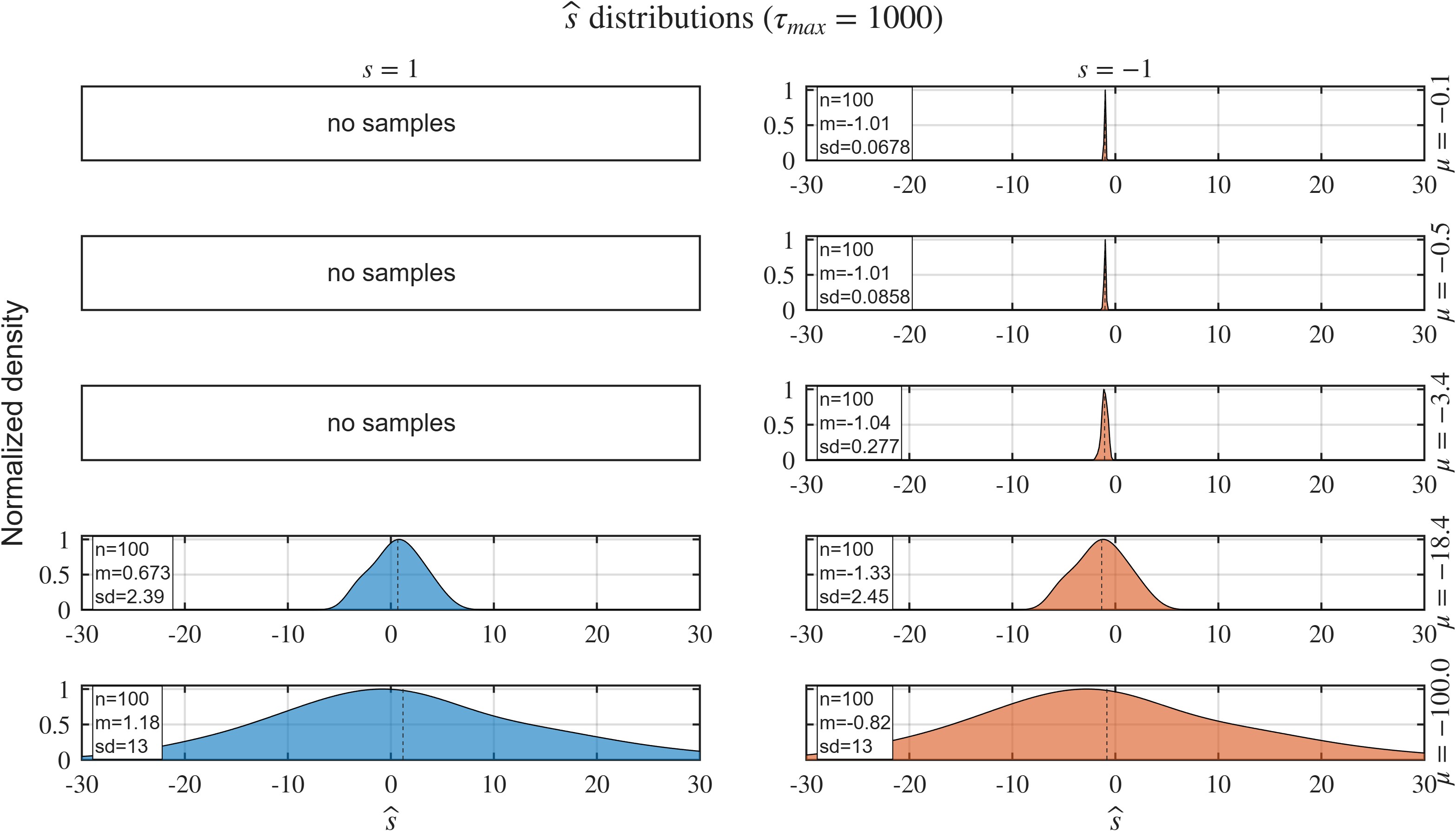}
    \caption{
Distributions of the estimated nonlinear coefficient $\hat{s}$ for $\tau_{\max}=1000$. Plotting conventions are as in Figure~\ref{fig:PC_ridge_mu_t=1}.
}
    \label{fig:PC_ridge_s_t=1000}
\end{figure}

The distributions of $\hat{s}$, shown in Figures~\ref{fig:PC_ridge_s_t=1}--\ref{fig:PC_ridge_s_t=1000}, show how well the two bifurcation types can be distinguished directly from the fitted nonlinear coefficient. The true values $s=-1$ and $s=1$ correspond to the supercritical and subcritical models, respectively. Their separation therefore determines the ROC and AUC performance discussed in Section~\ref{sec:Stochastic_pitchfork}.

\clearpage
\section{Validation of the Swift--Hohenberg fits}
\label{app:SH_ridgefits}

To provide additional insight into the fitting procedure and the ROC
analyses presented in Section~\ref{sec:Swift-Hohenberg}, we first assess
the recovery of the complete linear dispersion relation and subsequently
examine the fitted linear and nonlinear coefficients separately. For the
varying-$\gamma$ experiment, we additionally show how the spatial spectrum
changes as the linear spatial filtering becomes stronger. Together, these
results provide additional context for both the accuracy of the recovered
linear stability and the ability to distinguish the supercritical and
subcritical models.

\subsection{Fitted dispersion relations}
\label{app:SH_dispersion}

\begin{figure}[tbp]
    \centering
    \includegraphics[width=0.9\textwidth]{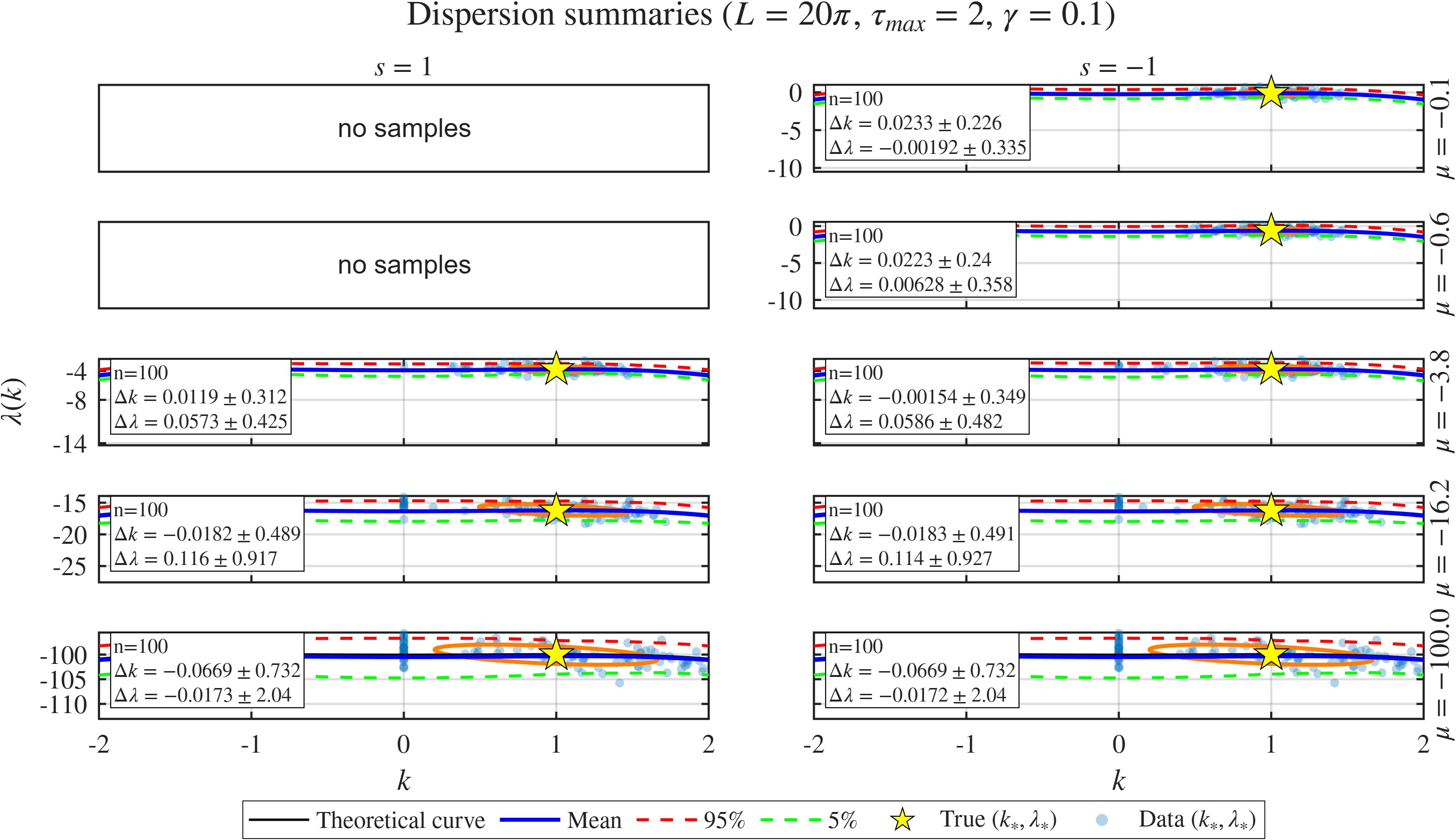}
    \caption{
Fitted dispersion relations for $\gamma=0.1$ and $\tau_{\max}=2$. Each panel corresponds to a fixed combination of $\mu$ and $s$. Plotting conventions are as in Figure~\ref{fig:SH_dispersion_representative}; the text boxes report the number of available realizations and the mean and standard deviation of $\Delta k$ and $\Delta\lambda$.
}
    \label{fig:SH_disp_gamma01}
\end{figure}

\begin{figure}[tbp]
    \centering
    \includegraphics[width=0.9\textwidth]{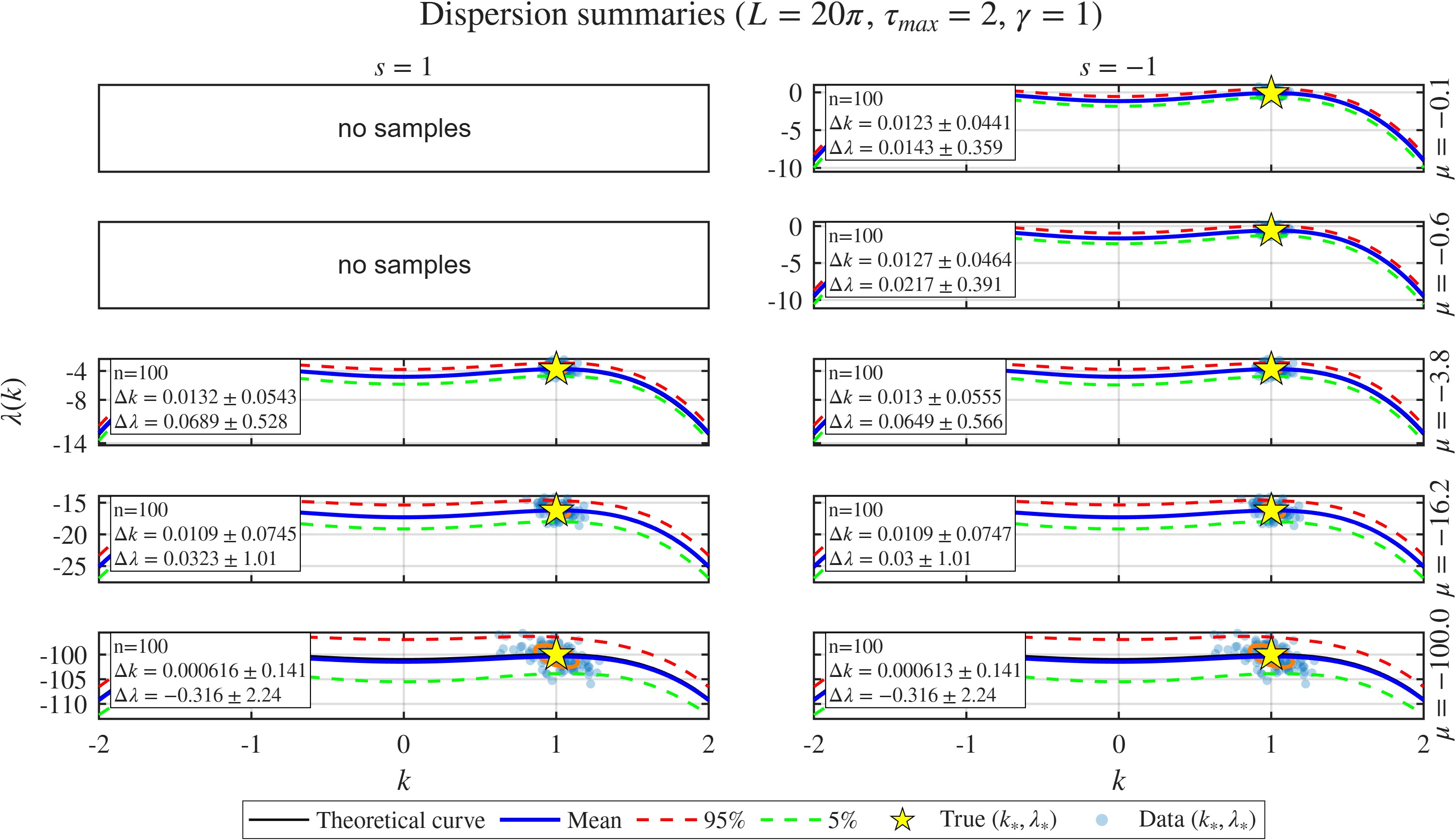}
    \caption{
Fitted dispersion relations for $\gamma=1$ and $\tau_{\max}=2$. Plotting conventions and error definitions are as in Figure~\ref{fig:SH_disp_gamma01}.
}
    \label{fig:SH_disp_gamma1}
\end{figure}

\begin{figure}[tbp]
    \centering
    \includegraphics[width=0.9\textwidth]{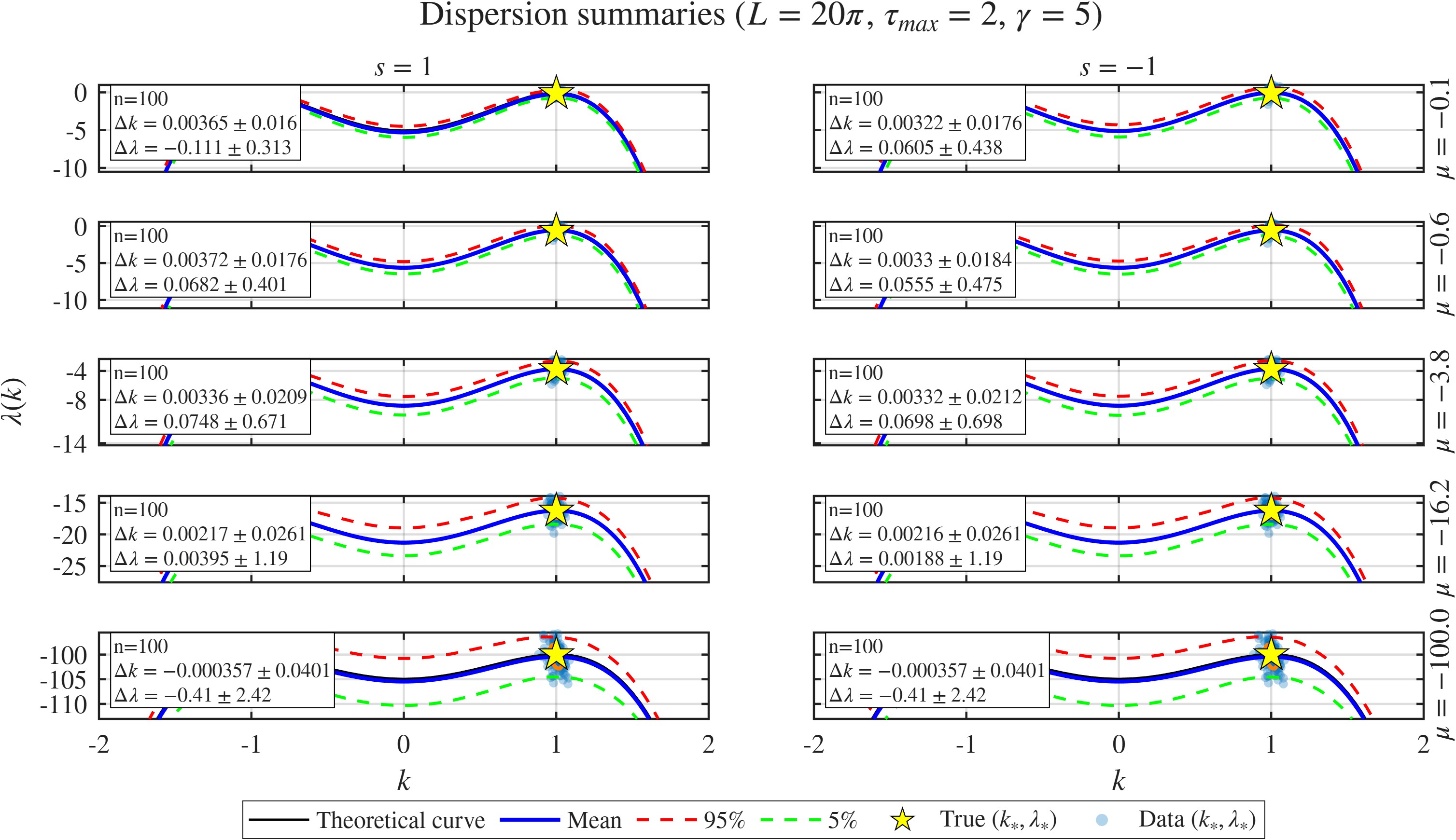}
    \caption{
Fitted dispersion relations for $\gamma=5$ and $\tau_{\max}=2$. Plotting conventions and error definitions are as in Figure~\ref{fig:SH_disp_gamma01}.
}
    \label{fig:SH_disp_gamma5}
\end{figure}

We first assess the combined recovery of the linear spatial dynamics by
comparing the fitted dispersion relations with the theoretical relation in
Equation~\eqref{eq:SH_dispersion}. Representative results are
shown for $\gamma=0.1$, $\gamma=1$, and $\gamma=5$ at the fixed observation
horizon $\tau_{\max}=2$ (Figures~\ref{fig:SH_disp_gamma01}--\ref{fig:SH_disp_gamma5}). These cases span the transition from a broad,
relatively flat maximum to a narrow and steep maximum around the critical
mode $k_*=1$.

For each realization, the fitted coefficients define
$\hat{\lambda}(k)$ through Equation~\eqref{eq:SH_estimated_dispersion}.
The errors in the estimated maximum are summarized by
\begin{align}
\Delta k
&=
\hat{k}_*-k_*,
&
\Delta\lambda
&=
\hat{\lambda}_*-\lambda_*,
\end{align}
with $k_*=1$ and $\lambda_*=\mu$. Across the three values of $\gamma$,
the mean fitted dispersion relations remain close to the theoretical
curves, showing that the linear spatial dynamics are generally recovered
well alongside the nonlinear coefficient.

The two error measures reflect different aspects of this reconstruction.
For small $\gamma$, the dispersion relation is relatively flat near its
maximum, so uncertainty in the location $\hat{k}_*$ produces only a small
change in the corresponding growth rate. As $\gamma$ increases, the
maximum becomes increasingly sharp, constraining $\hat{k}_*$ more tightly
and reducing the spread in $\Delta k$. At the same time, deviations in the
fitted dispersion curve can produce larger errors in
$\Delta\lambda$. Hence, improved localization of the dominant wavenumber
does not necessarily imply a more accurate estimate of its growth rate.

\subsection{Fitting of the linear coefficient \texorpdfstring{$\hat{a}_0$}{a0}}

We next examine the fitted coefficients individually. For the linear fit,
the displayed coefficient is $\hat{a}_0$, whose theoretical value is
\begin{align}
a_0=\mu-\gamma.
\end{align}

\subsubsection*{Varying the observation horizon}

\begin{figure}[tbp]
    \centering
    \includegraphics[width=0.8\textwidth]{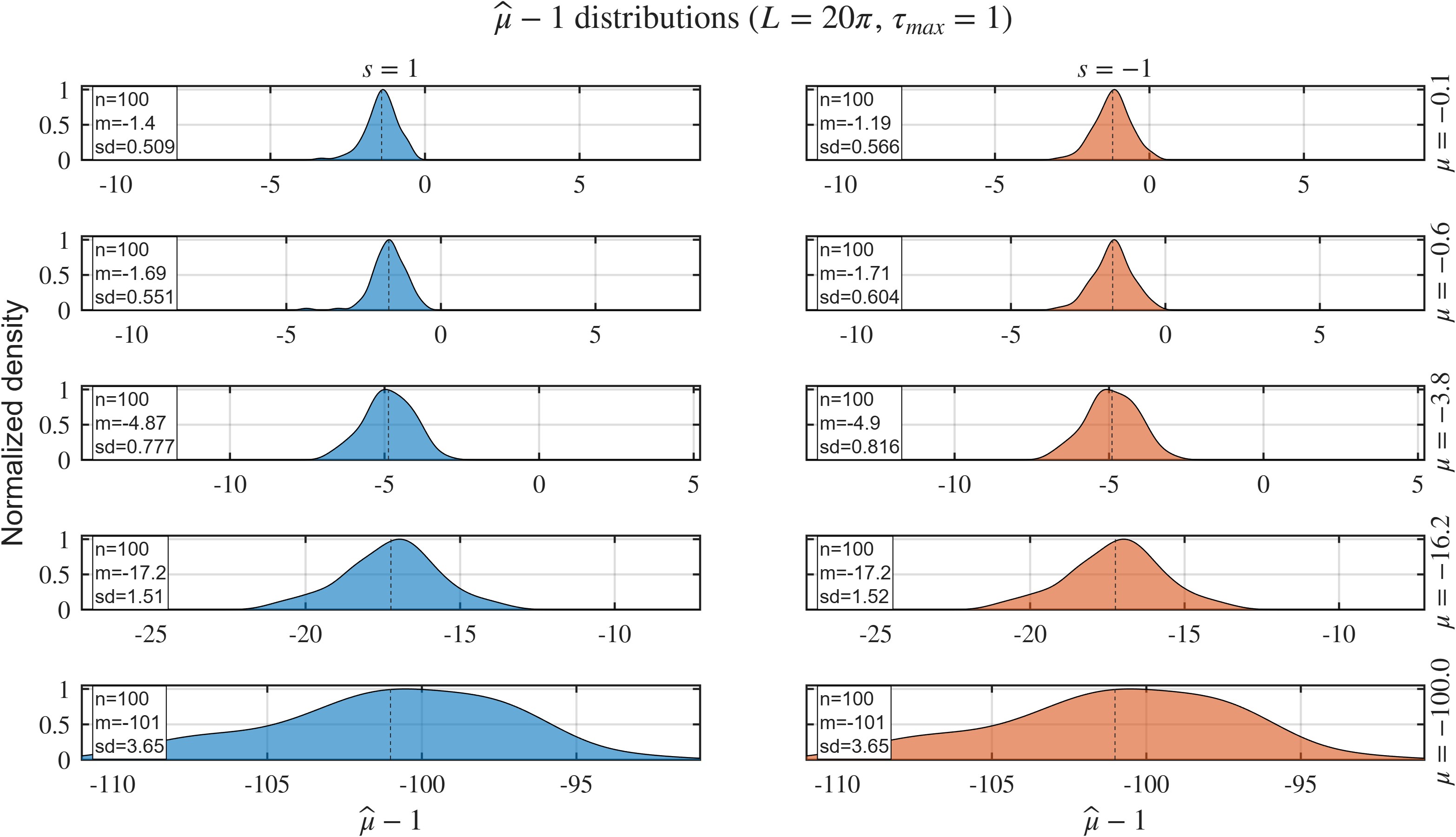}
    \caption{
Distributions of the estimated linear coefficient $\hat{a}_0$ for $\tau_{\max}=1$ and $\gamma=1$, for which the theoretical value is $a_0=\mu-1$. The panel layout and plotting conventions are as in Figure~\ref{fig:PC_ridge_mu_t=1}.
}
    \label{fig:SH_ridge_a0_tau=1}
\end{figure}

\begin{figure}[tbp]
    \centering
    \includegraphics[width=0.8\textwidth]{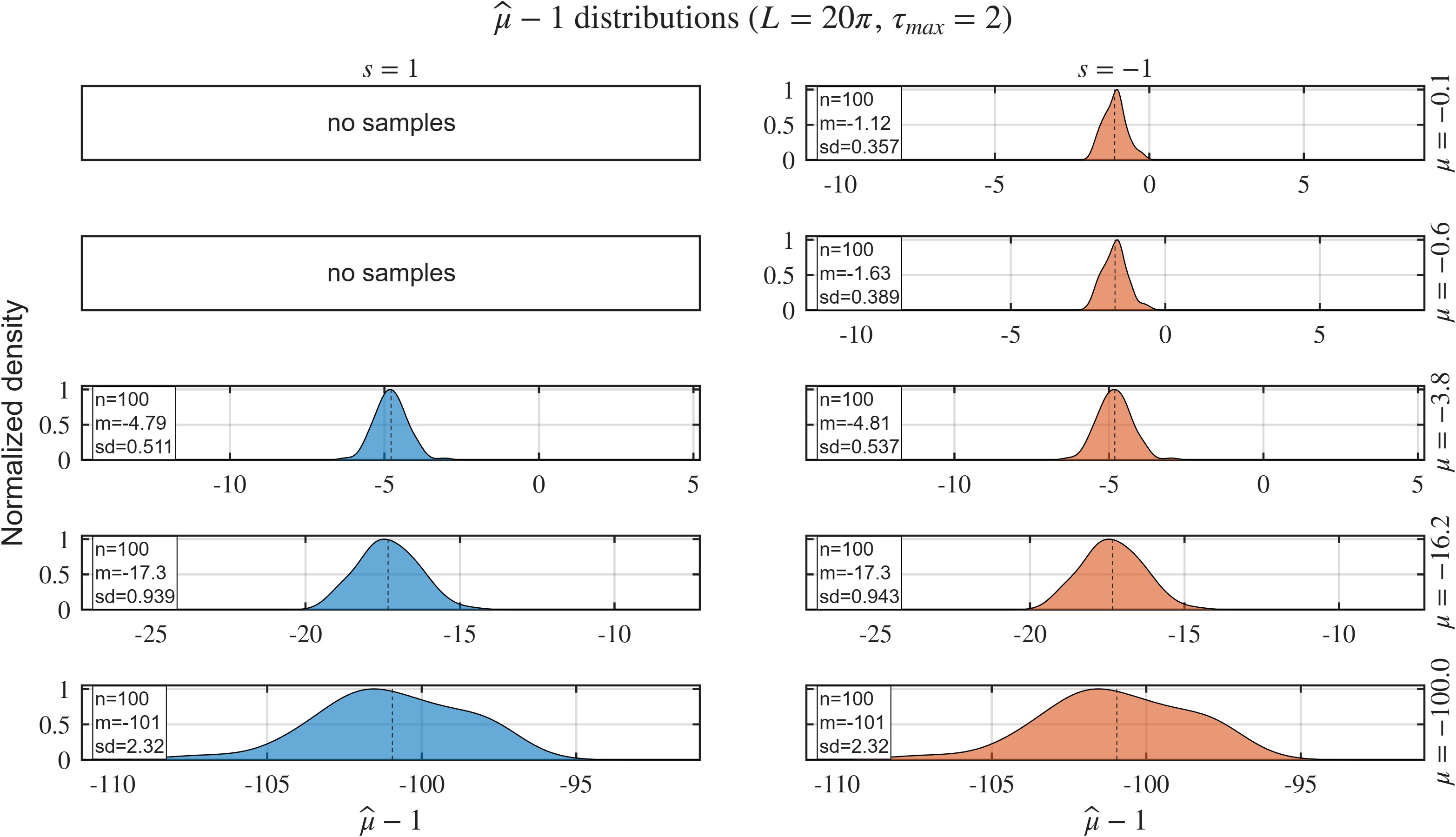}
    \caption{
Distributions of the estimated linear coefficient $\hat{a}_0$ for $\tau_{\max}=2$ and $\gamma=1$. Plotting conventions are as in Figure~\ref{fig:SH_ridge_a0_tau=1}.
}
    \label{fig:SH_ridge_a0_tau=2}
\end{figure}

\begin{figure}[tbp]
    \centering
    \includegraphics[width=0.8\textwidth]{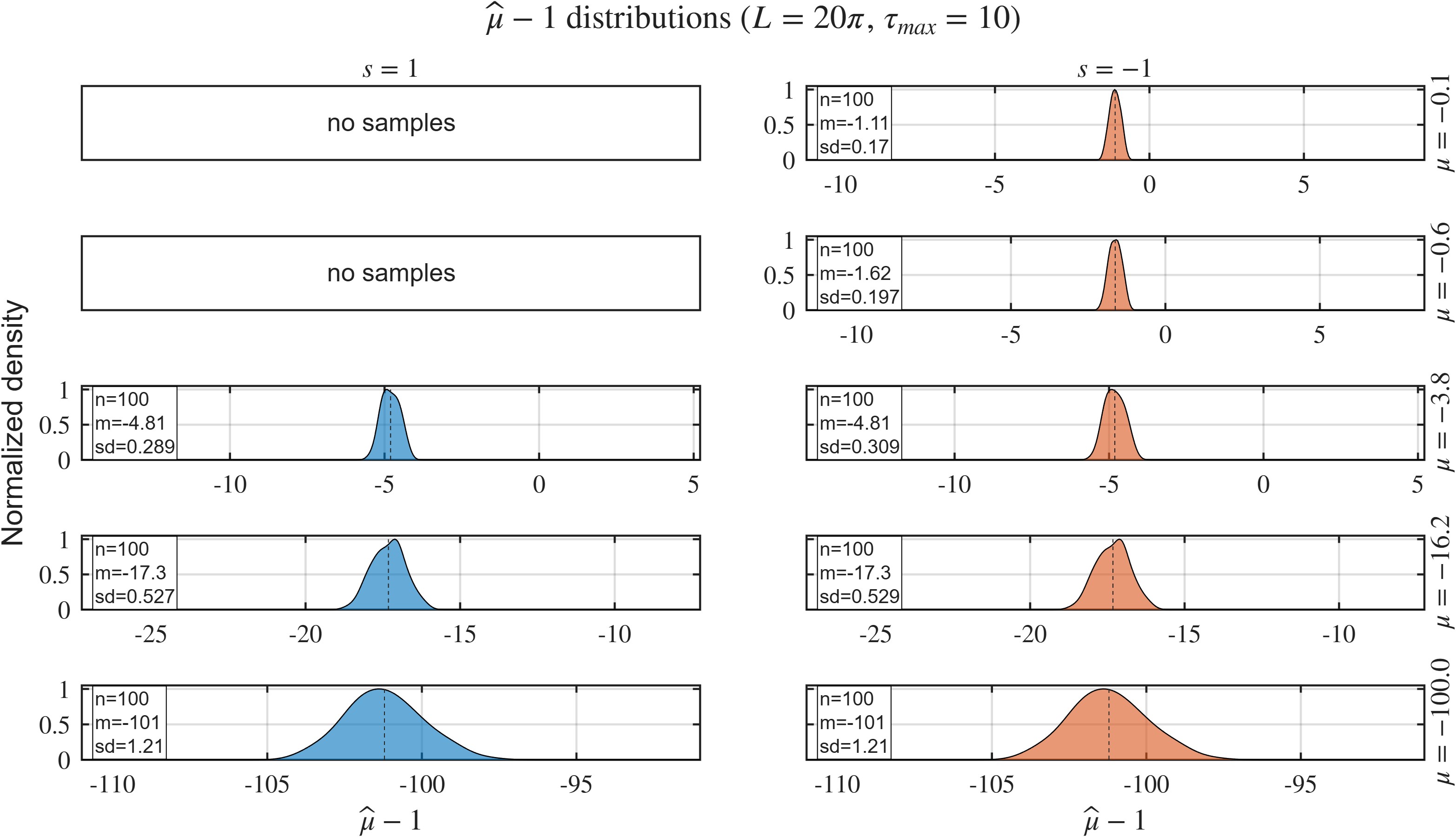}
    \caption{
Distributions of the estimated linear coefficient $\hat{a}_0$ for $\tau_{\max}=10$ and $\gamma=1$. Plotting conventions are as in Figure~\ref{fig:SH_ridge_a0_tau=1}.
}
    \label{fig:SH_ridge_a0_tau=10}
\end{figure}

\begin{figure}[tbp]
    \centering
    \includegraphics[width=0.8\textwidth]{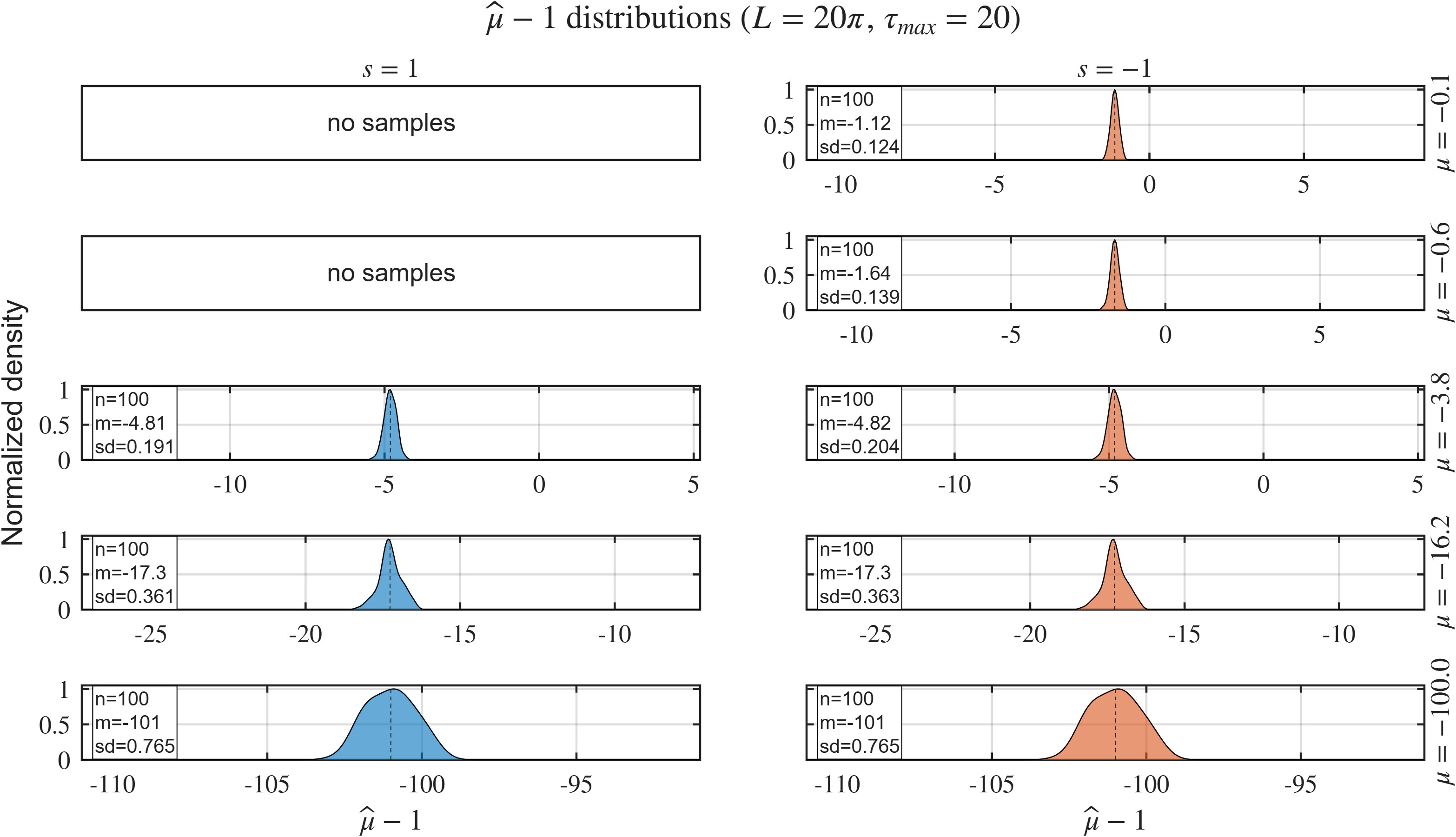}
    \caption{
Distributions of the estimated linear coefficient $\hat{a}_0$ for $\tau_{\max}=20$ and $\gamma=1$. Plotting conventions are as in Figure~\ref{fig:SH_ridge_a0_tau=1}.
}
    \label{fig:SH_ridge_a0_tau=20}
\end{figure}

Figures~\ref{fig:SH_ridge_a0_tau=1}--\ref{fig:SH_ridge_a0_tau=20} show the distributions of the fitted linear coefficient for varying observation horizons at fixed $\gamma=1$, for which the theoretical value is $a_0=\mu-1$.

\subsubsection*{Varying the spatial parameter \texorpdfstring{$\gamma$}{gamma}}

\begin{figure}[tbp]
    \centering
    \includegraphics[width=0.8\textwidth]{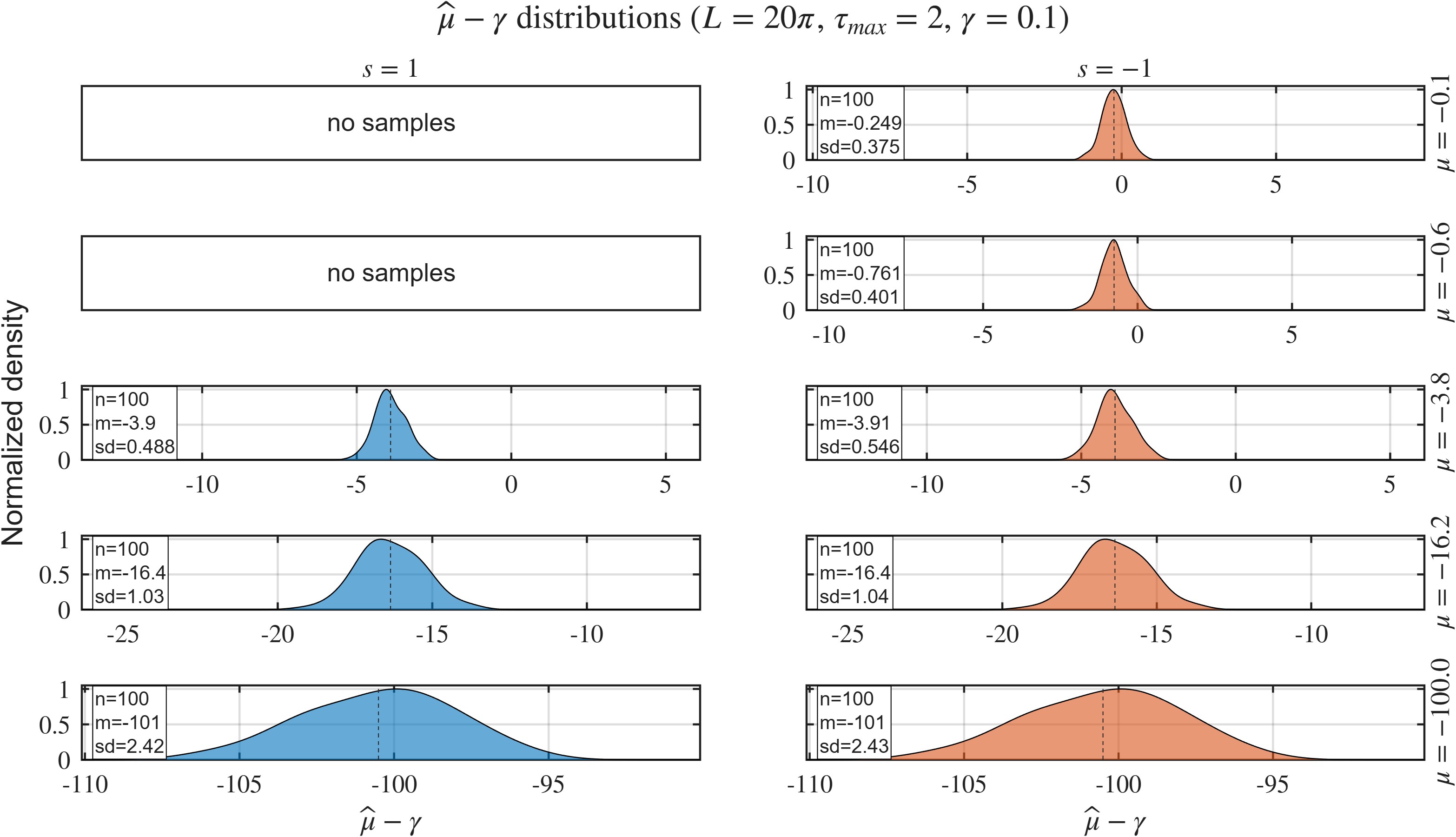}
    \caption{
Distributions of the estimated linear coefficient $\hat{a}_0$ for $\gamma=0.1$ and $\tau_{\max}=2$, with theoretical value $a_0=\mu-0.1$. Plotting conventions are as in Figure~\ref{fig:SH_ridge_a0_tau=1}.
}
    \label{fig:SH_ridge_a0_gamma=0.1}
\end{figure}

\begin{figure}[tbp]
    \centering
    \includegraphics[width=0.8\textwidth]{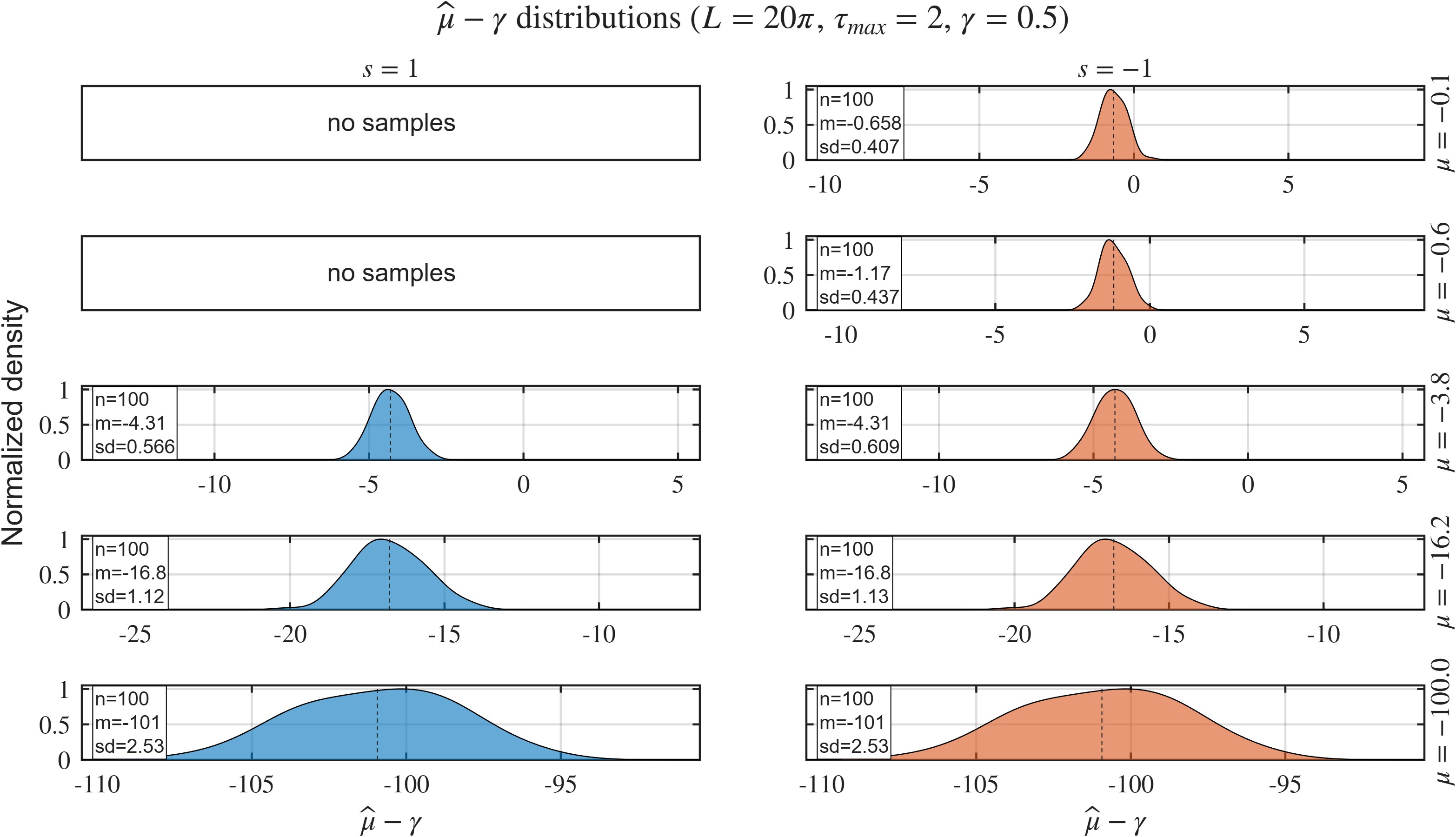}
    \caption{
Distributions of the estimated linear coefficient $\hat{a}_0$ for $\gamma=0.5$ and $\tau_{\max}=2$, with theoretical value $a_0=\mu-0.5$. Plotting conventions are as in Figure~\ref{fig:SH_ridge_a0_tau=1}.
}
    \label{fig:SH_ridge_a0_gamma=0.5}
\end{figure}

\begin{figure}[tbp]
    \centering
    \includegraphics[width=0.8\textwidth]{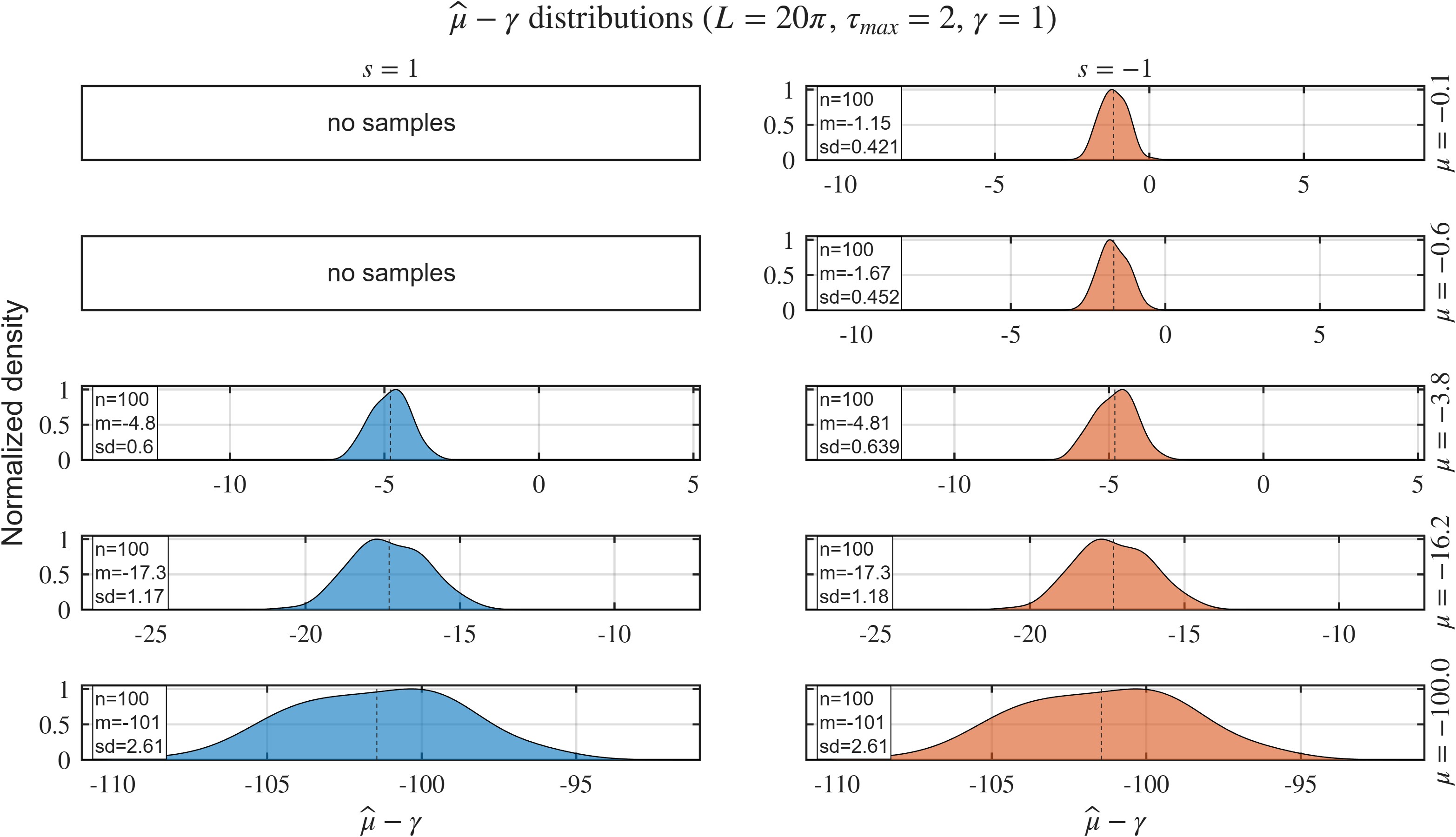}
    \caption{
Distributions of the estimated linear coefficient $\hat{a}_0$ for $\gamma=1$ and $\tau_{\max}=2$, with theoretical value $a_0=\mu-1$. Plotting conventions are as in Figure~\ref{fig:SH_ridge_a0_tau=1}.
}
    \label{fig:SH_ridge_a0_gamma=1}
\end{figure}

\begin{figure}[tbp]
    \centering
    \includegraphics[width=0.8\textwidth]{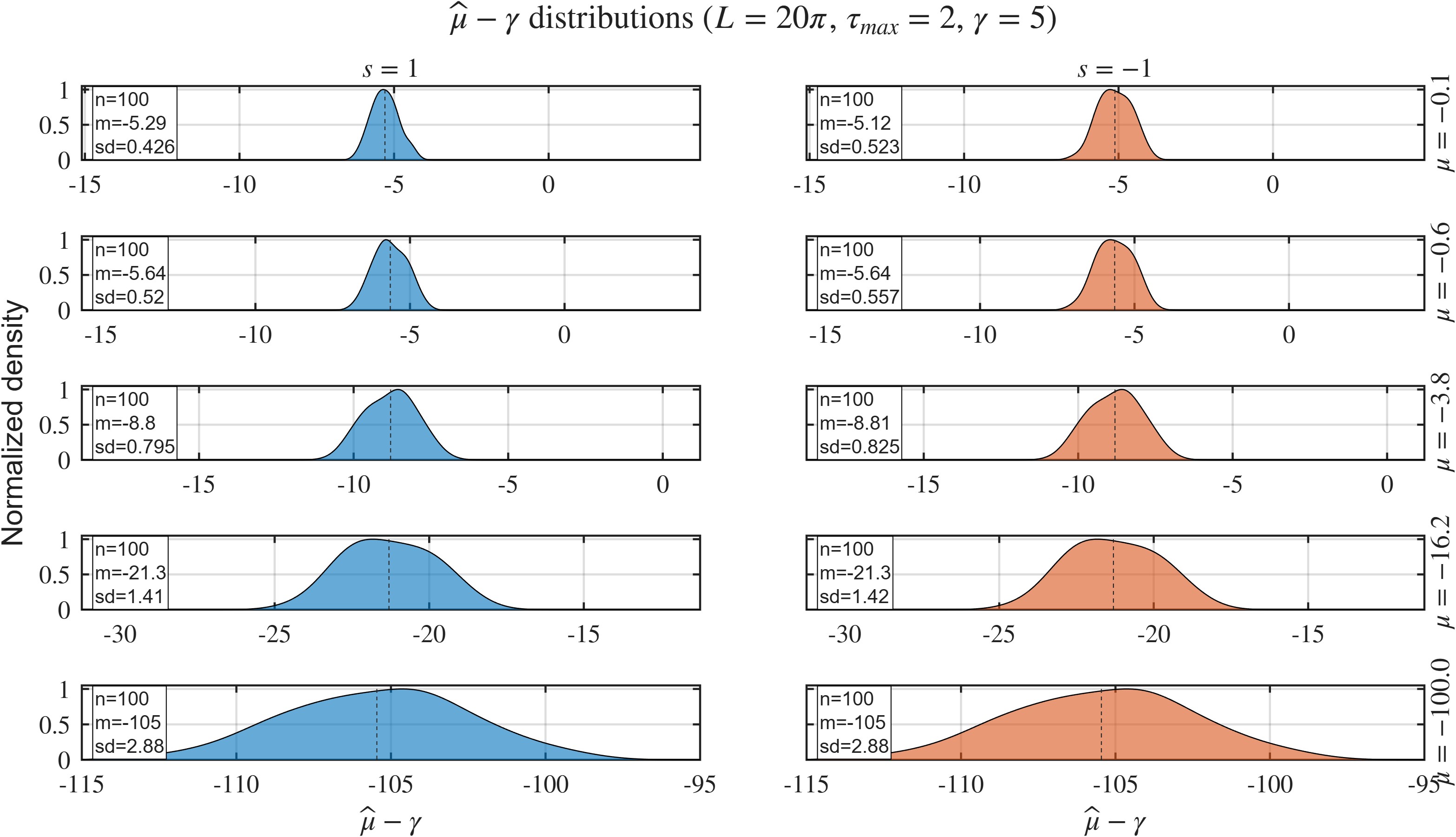}
    \caption{
Distributions of the estimated linear coefficient $\hat{a}_0$ for $\gamma=5$ and $\tau_{\max}=2$, with theoretical value $a_0=\mu-5$. Plotting conventions are as in Figure~\ref{fig:SH_ridge_a0_tau=1}.
}
    \label{fig:SH_ridge_a0_gamma=5}
\end{figure}

Figures~\ref{fig:SH_ridge_a0_gamma=0.1}--\ref{fig:SH_ridge_a0_gamma=5} show the corresponding distributions when $\gamma$ is varied at fixed $\tau_{\max}=2$. In this case the theoretical linear coefficient changes with both $\mu$ and $\gamma$ according to $a_0=\mu-\gamma$.

\subsection{Fitting of the nonlinear coefficient \texorpdfstring{$\hat{s}$}{s}}

The nonlinear coefficient determines the predicted bifurcation type, with $s=-1$ corresponding to the supercritical model and $s=1$ to the subcritical model. The overlap between the fitted distributions of $\hat{s}$ therefore determines how well the supercritical and subcritical models can be distinguished, as quantified by the ROC and AUC analyses in Section~\ref{sec:Swift-Hohenberg}.

\subsubsection*{Varying the observation horizon}
\begin{figure}[tbp]
    \centering
    \includegraphics[width=0.8\textwidth]{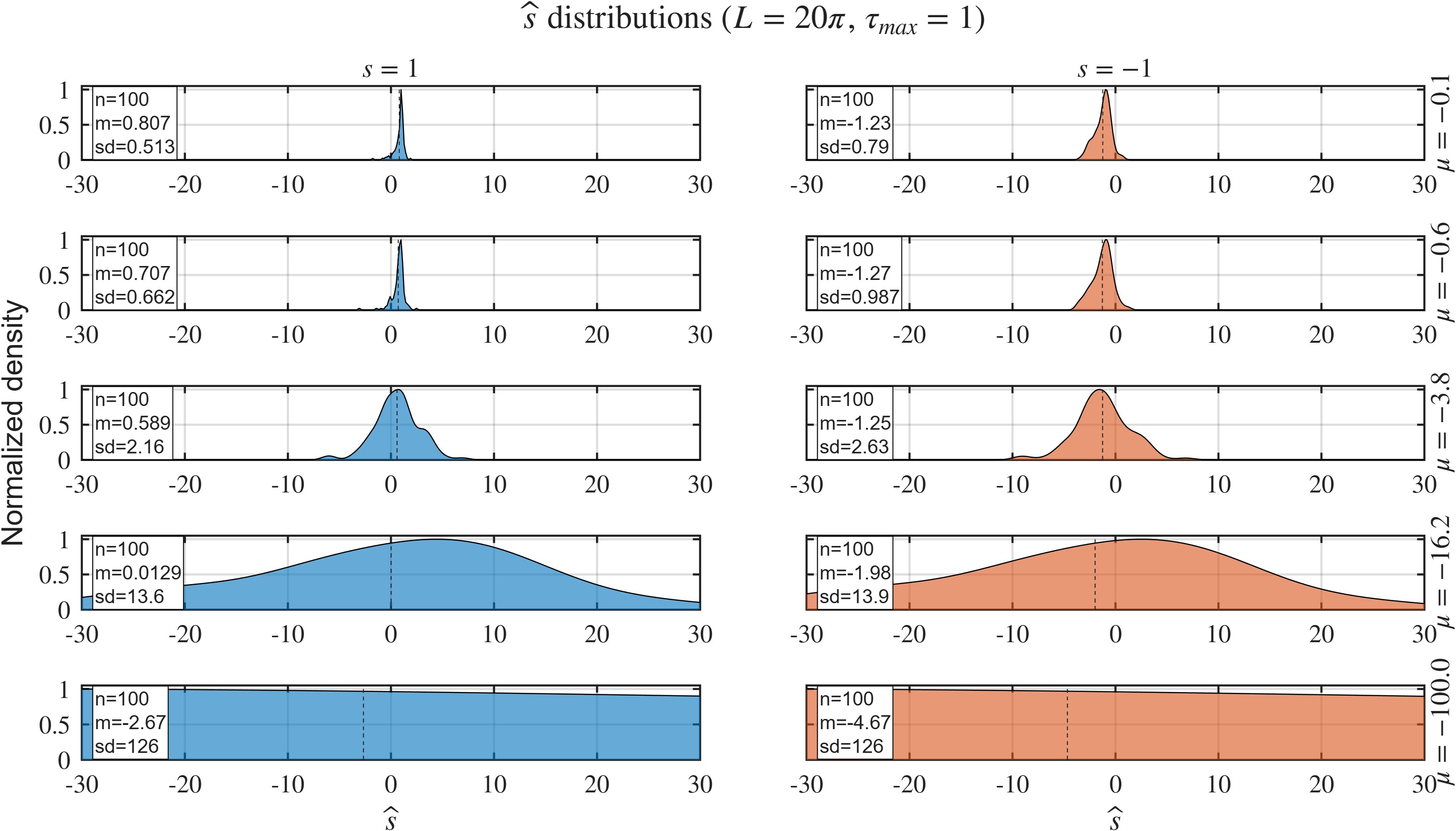}
    \caption{
Distributions of the estimated nonlinear coefficient $\hat{s}$ for $\tau_{\max}=1$ and $\gamma=1$. Plotting conventions are as in Figure~\ref{fig:SH_ridge_a0_tau=1}.
}
    \label{fig:SH_ridge_s_tau=1}
\end{figure}

\begin{figure}[tbp]
    \centering
    \includegraphics[width=0.8\textwidth]{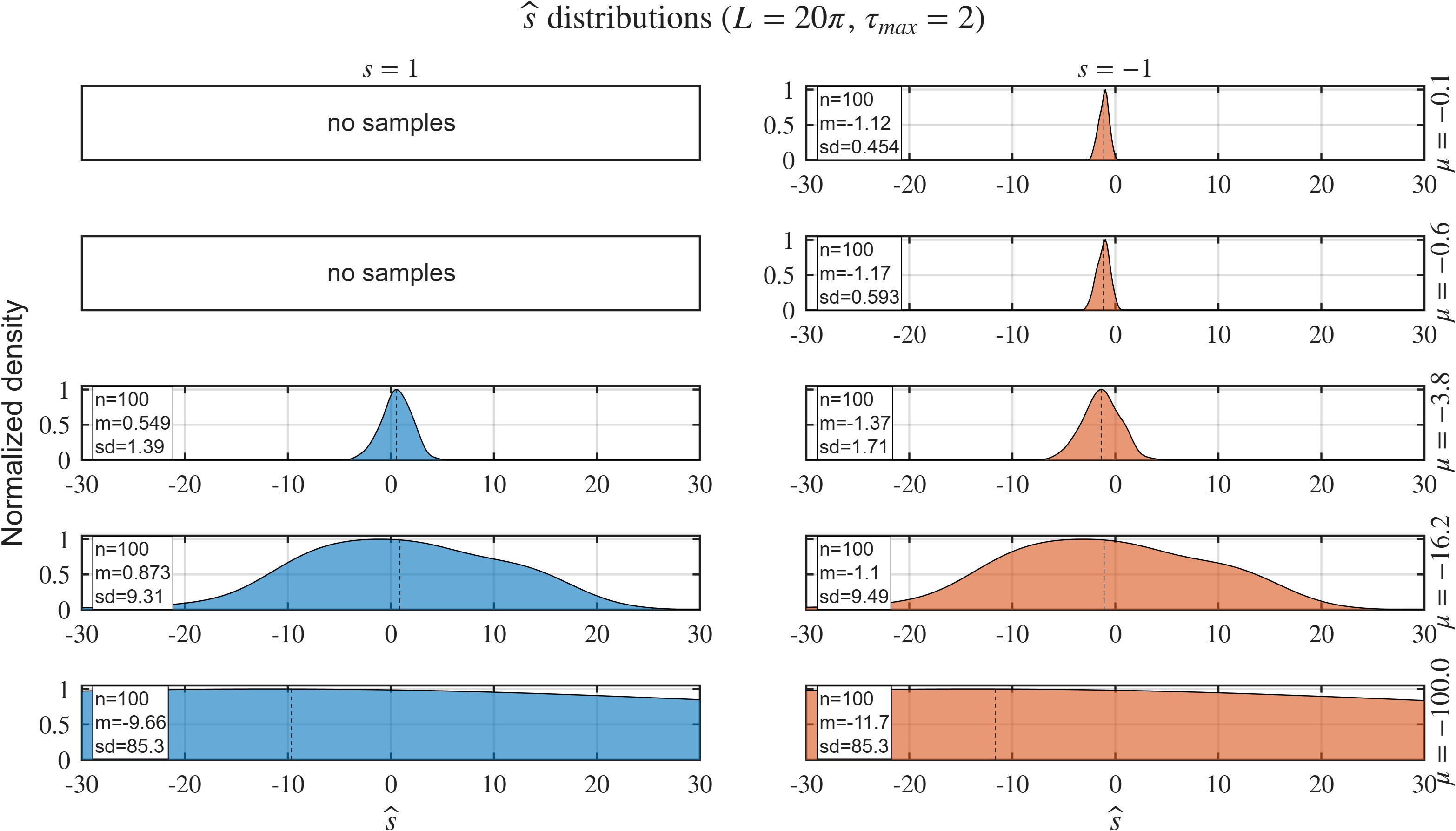}
    \caption{
Distributions of the estimated nonlinear coefficient $\hat{s}$ for $\tau_{\max}=2$ and $\gamma=1$. Plotting conventions are as in Figure~\ref{fig:SH_ridge_a0_tau=1}.
}
    \label{fig:SH_ridge_s_tau=2}
\end{figure}

\begin{figure}[tbp]
    \centering
    \includegraphics[width=0.8\textwidth]{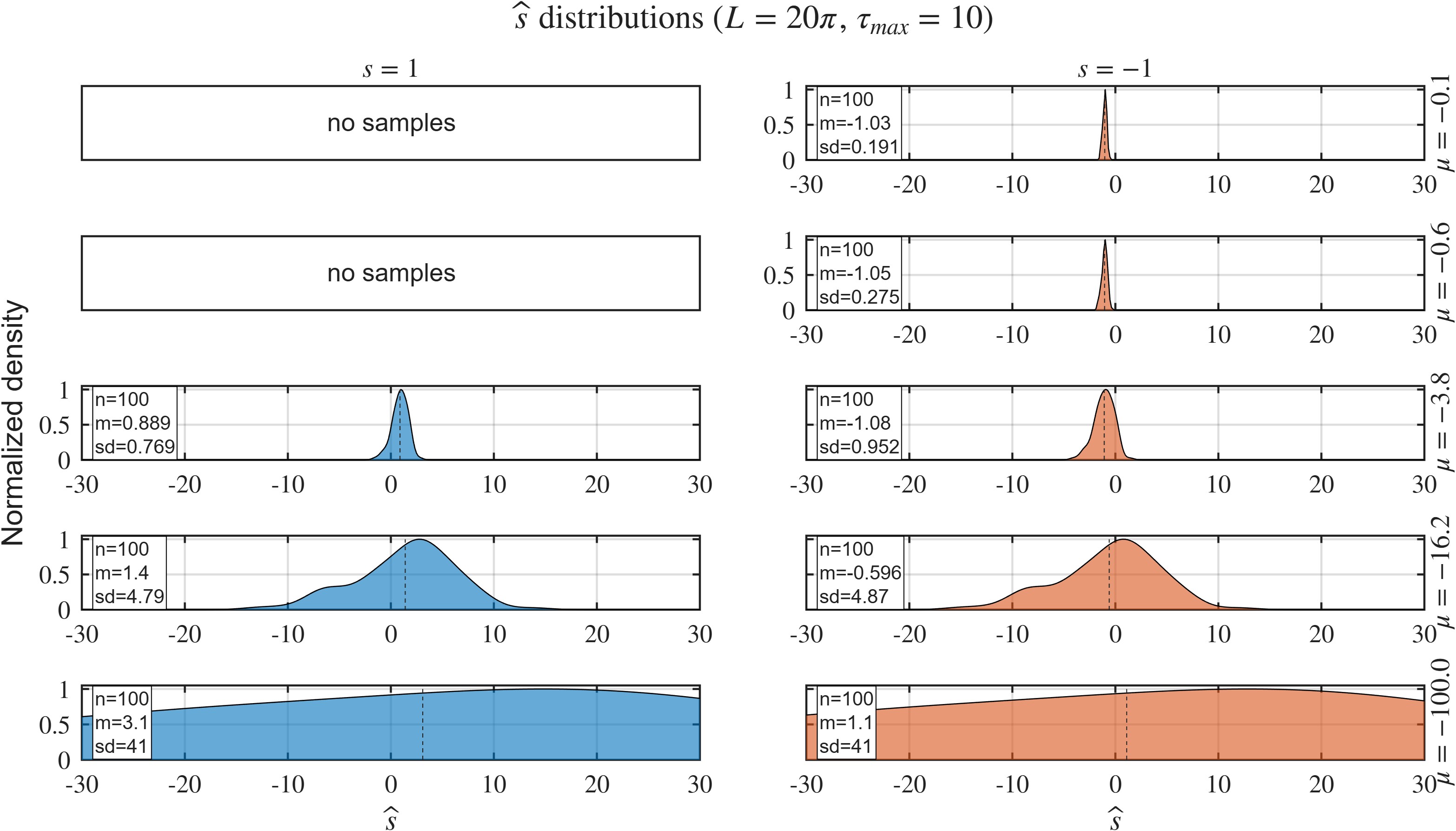}
    \caption{
Distributions of the estimated nonlinear coefficient $\hat{s}$ for $\tau_{\max}=10$ and $\gamma=1$. Plotting conventions are as in Figure~\ref{fig:SH_ridge_a0_tau=1}.
}
    \label{fig:SH_ridge_s_tau=10}
\end{figure}

\begin{figure}[tbp]
    \centering
    \includegraphics[width=0.8\textwidth]{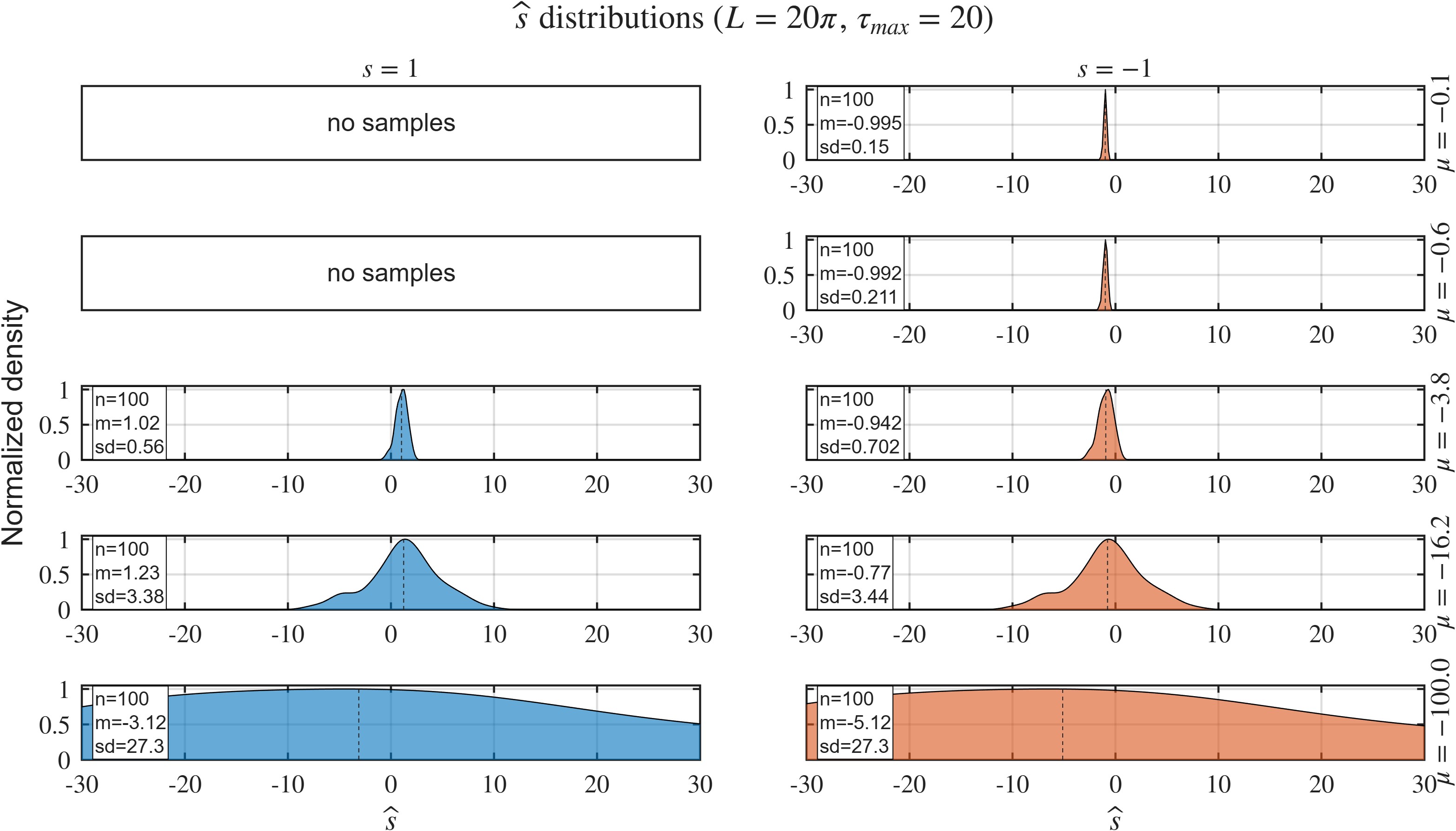}
    \caption{
Distributions of the estimated nonlinear coefficient $\hat{s}$ for $\tau_{\max}=20$ and $\gamma=1$. Plotting conventions are as in Figure~\ref{fig:SH_ridge_a0_tau=1}.
}
    \label{fig:SH_ridge_s_tau=20}
\end{figure}

Figures~\ref{fig:SH_ridge_s_tau=1}--\ref{fig:SH_ridge_s_tau=20} show the distributions of $\hat{s}$ for representative observation horizons at fixed $\gamma=1$. With increasing observation horizon, the fitted distributions generally become narrower and the supercritical and subcritical cases become more clearly separated, consistent with the improved classification performance reported in Section~\ref{sec:SH_results}.

\subsubsection*{Varying the spatial parameter \texorpdfstring{$\gamma$}{gamma}}

\begin{figure}[tbp]
    \centering
    \includegraphics[width=\textwidth]
    {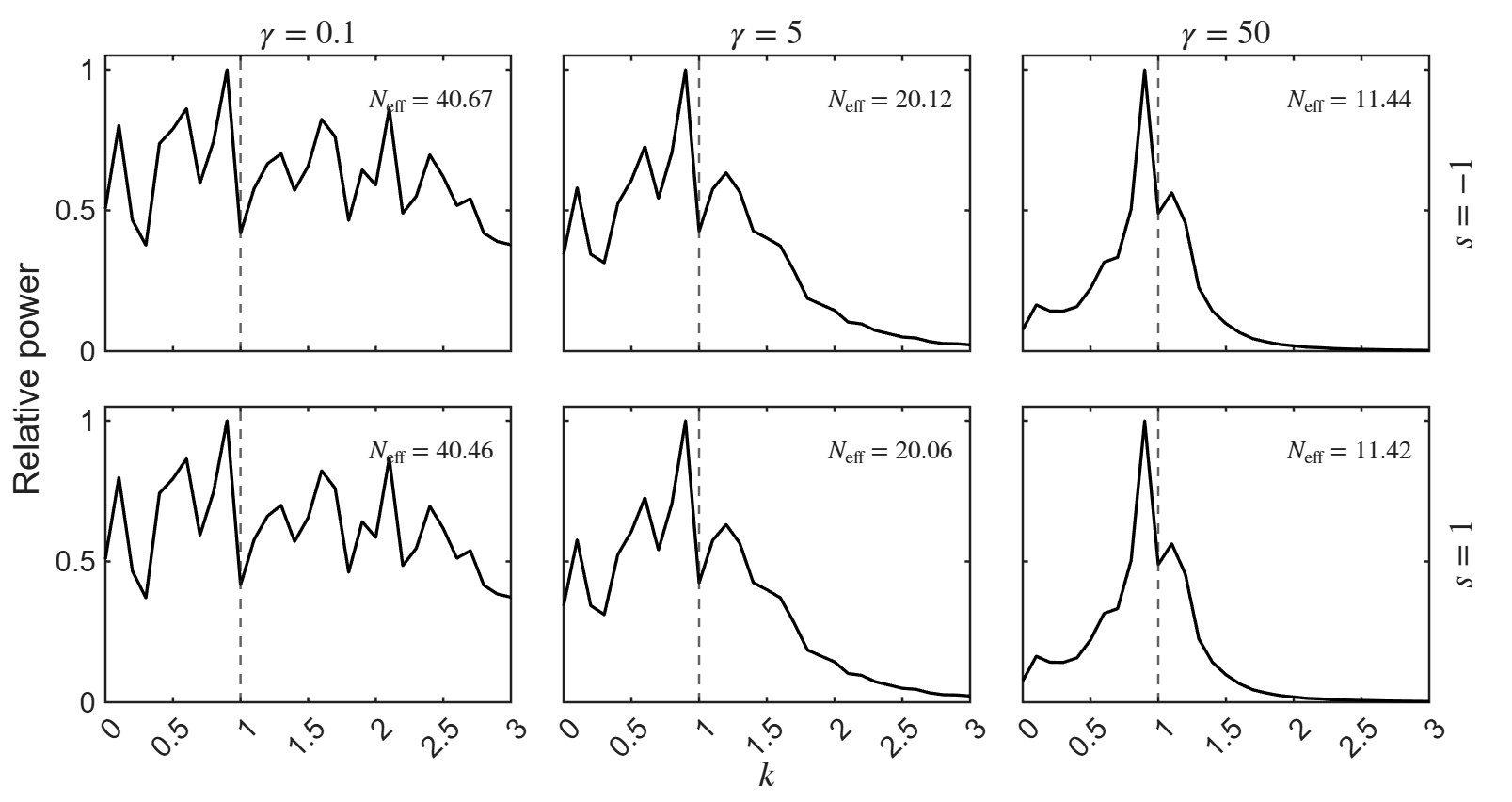}
    \caption{
Time-averaged one-sided spatial power spectra for the reduced stochastic Swift--Hohenberg model at $\mu=-10$ and $\tau_{\max}=2$. Columns correspond to $\gamma=0.1$, $5$, and $50$, and rows to the supercritical ($s=-1$) and subcritical ($s=1$) models. Each spectrum is normalized by its maximum value for visualization, and the vertical dashed line marks the critical wavenumber $k=1$. The displayed effective spectral mode count is $N_{\mathrm{eff}}=(\sum_k p_k^2)^{-1}$, where $p_k=P_k/\sum_j P_j$ is the fraction of the total spectral power in mode $k$. The first $25\%$ of each simulation is discarded before the time average is computed.
}
    \label{fig:SH_gamma_spectra}
\end{figure}

\begin{figure}[tbp]
    \centering
    \includegraphics[width=0.8\textwidth]{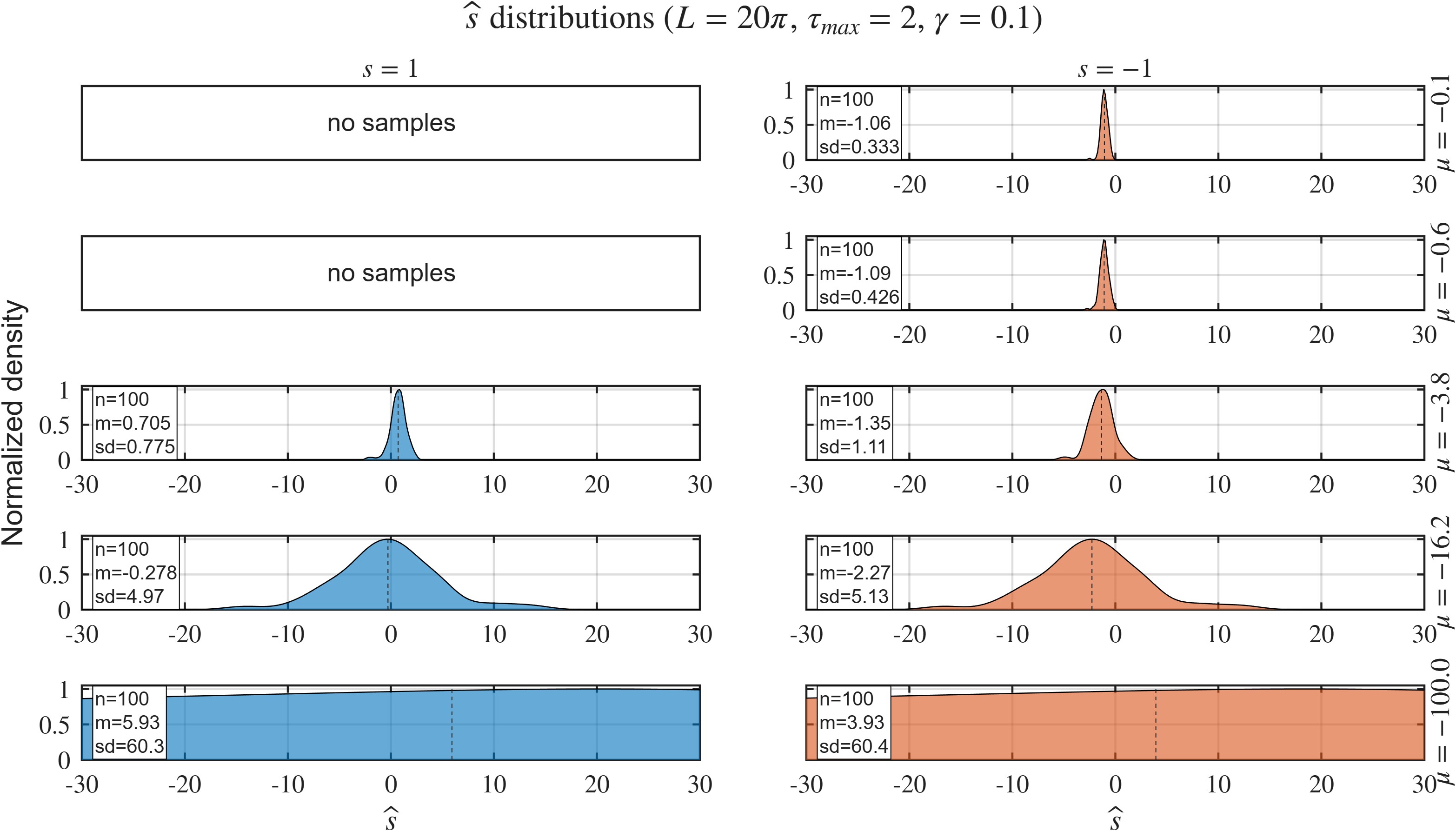}
    \caption{
Distributions of the estimated nonlinear coefficient $\hat{s}$ for $\gamma=0.1$ and $\tau_{\max}=2$. Plotting conventions are as in Figure~\ref{fig:SH_ridge_a0_tau=1}.
}
    \label{fig:SH_ridge_s_gamma=0.1}
\end{figure}

\begin{figure}[tbp]
    \centering
    \includegraphics[width=0.8\textwidth]{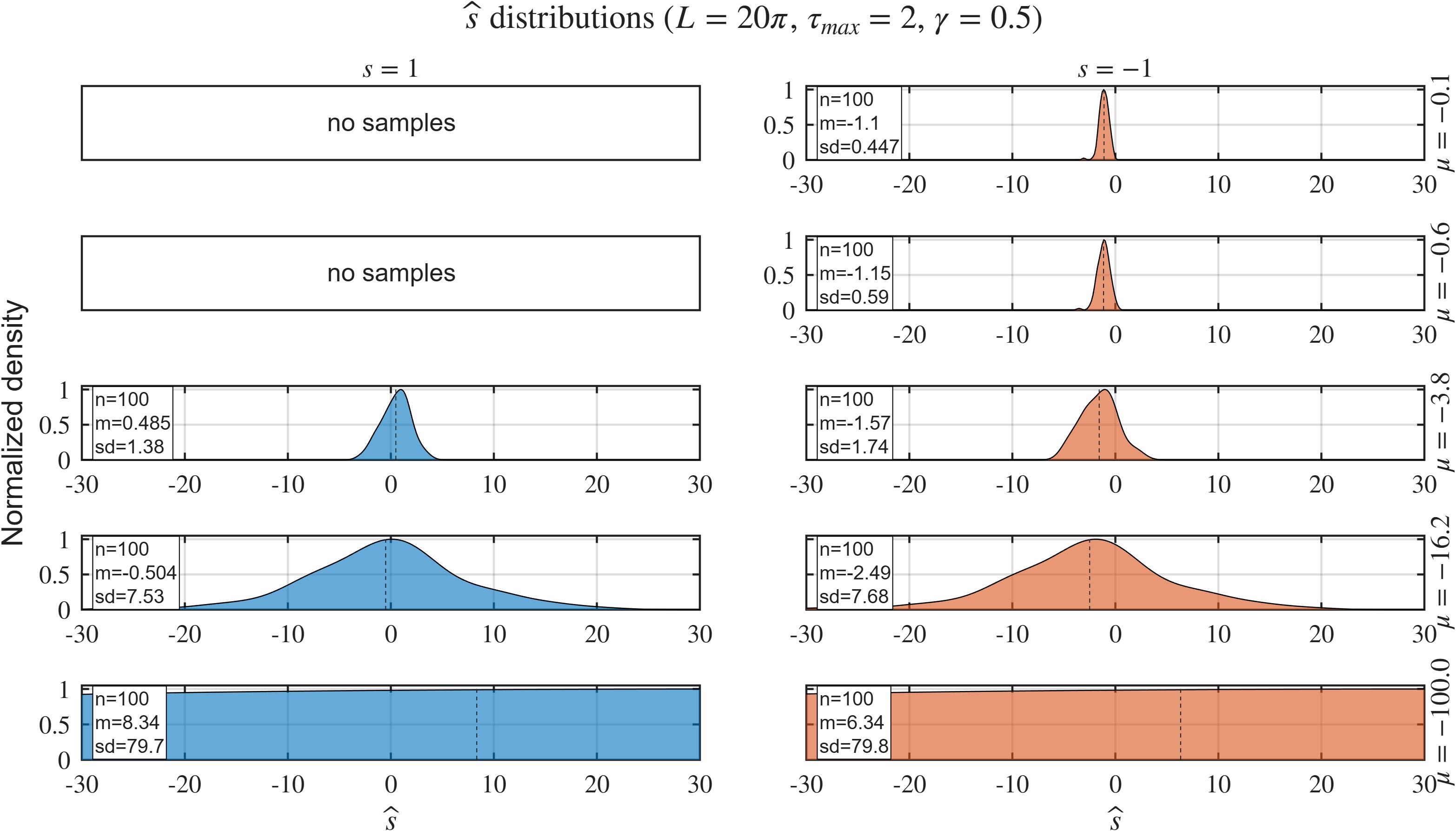}
    \caption{
Distributions of the estimated nonlinear coefficient $\hat{s}$ for $\gamma=0.5$ and $\tau_{\max}=2$. Plotting conventions are as in Figure~\ref{fig:SH_ridge_a0_tau=1}.
}
    \label{fig:SH_ridge_s_gamma=0.5}
\end{figure}

\begin{figure}[tbp]
    \centering
    \includegraphics[width=0.8\textwidth]{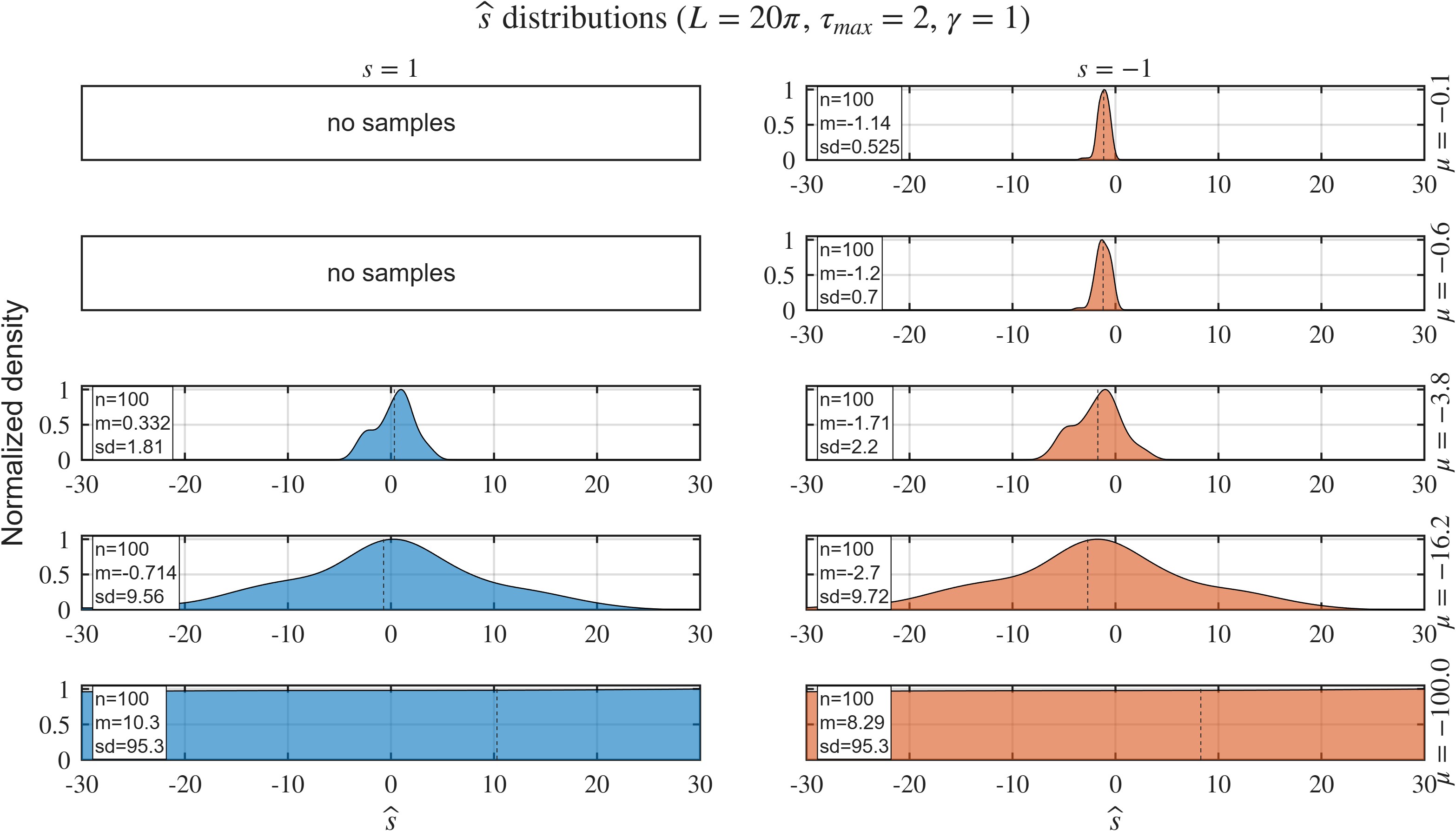}
    \caption{
Distributions of the estimated nonlinear coefficient $\hat{s}$ for $\gamma=1$ and $\tau_{\max}=2$. Plotting conventions are as in Figure~\ref{fig:SH_ridge_a0_tau=1}.
}
    \label{fig:SH_ridge_s_gamma=1}
\end{figure}

\begin{figure}[tbp]
    \centering
    \includegraphics[width=0.8\textwidth]{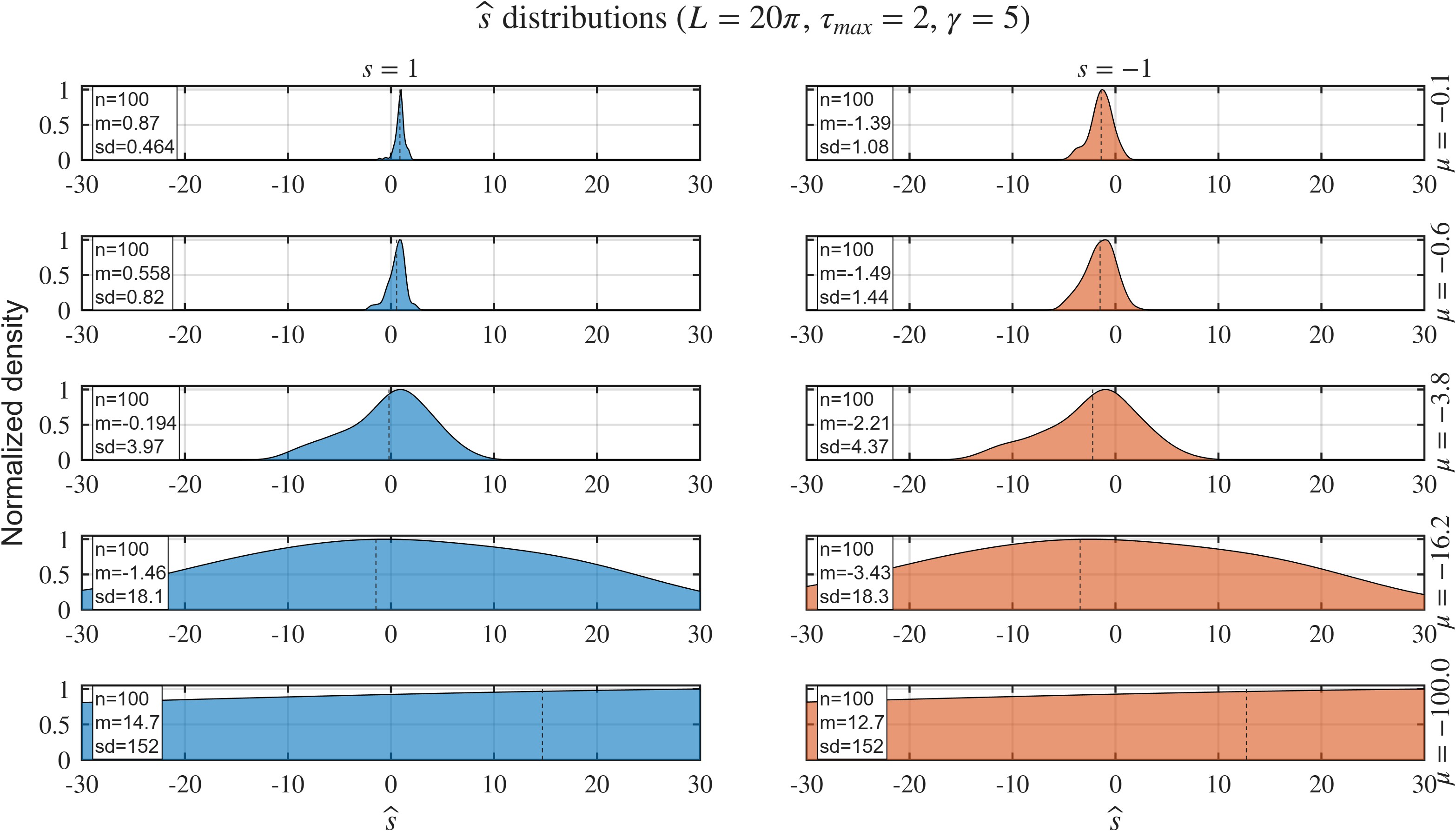}
    \caption{
Distributions of the estimated nonlinear coefficient $\hat{s}$ for $\gamma=5$ and $\tau_{\max}=2$. Plotting conventions are as in Figure~\ref{fig:SH_ridge_a0_tau=1}.
}
    \label{fig:SH_ridge_s_gamma=5}
\end{figure}

Figure~\ref{fig:SH_gamma_spectra} shows representative power spectra over
a wider illustrative range of $\gamma$. Increasing $\gamma$ concentrates
the fluctuations more strongly around the critical wavenumber $k=1$ and
decreases the effective number of substantially contributing modes. The
same spectral narrowing is visible for both bifurcation criticality types. Its possible
implications for the identifiability of the cubic term are discussed in
Section~\ref{sec:discussion_spatial}. The corresponding distributions of $\hat{s}$ are shown in Figures~\ref{fig:SH_ridge_s_gamma=0.1}--\ref{fig:SH_ridge_s_gamma=5}.

\clearpage
\section{Effect of the critical spatial scale}
\label{app:SH_spatial_scaling}

The physical critical wavenumber $k_\ast$ enters the reduced spatial
coordinate through
\begin{align}
X=k_\ast x.
\end{align}
Consequently, for a fixed physical observation domain of length $L_x$ and
a fixed physical spatial grid, increasing $k_\ast$ increases both the
reduced domain length $L$ and the reduced grid spacing $\Delta X$,
\begin{align}
L=k_\ast L_x,
\qquad
\Delta X=k_\ast\Delta x,
\end{align}
while the number of spatial observations remains unchanged. Since the
critical physical wavelength is $2\pi/k_\ast$, increasing $k_\ast$ also
means that more critical wavelengths are contained within the same physical
observation domain. To isolate this direct spatial effect from the
simultaneous dependence of $\mu$, $\gamma$, and $\tau_{\max}$ on
$k_\ast$, we vary $L$ and $\Delta X$ proportionally while keeping the
reduced dynamical parameters and the number of spatial observations fixed.

We use $N_x=100$, $\gamma=2$, $\tau_{\max}=1$, and
$N_{\mathrm{sim}}=100$ independent realizations for each parameter
combination. The parameter $\mu$ is varied over $20$ logarithmically spaced
values in $[-10^2,-10^{-1}]$, with $s\in\{-1,1\}$, and the numerical time step is
$\Delta\tau=2.5\times10^{-4}$. The reduced domain length is varied as
\begin{align}
L=\{14\pi,16\pi,\ldots,30\pi\},
\qquad
\Delta X=\frac{L}{N_x}.
\end{align}
Increasing $L$ and $\Delta X$ in this way corresponds directly to the
spatial change produced by increasing $k_\ast$ for a fixed physical
observation setup.

\begin{figure}[tbp]
    \centering
    \includegraphics[width=0.75\textwidth]{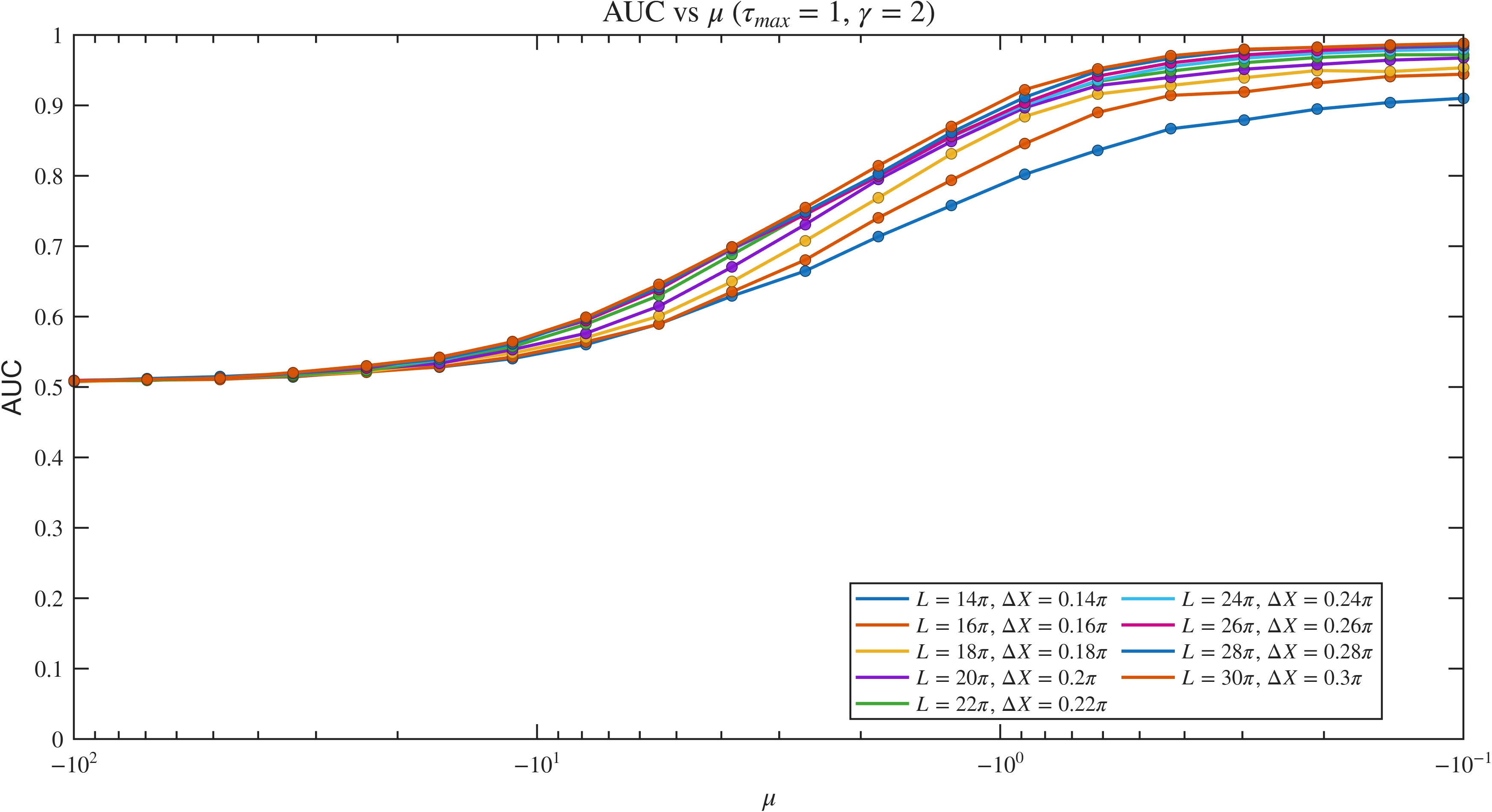}
    \caption{
    AUC for bifurcation criticality classification as a function of $\mu$ for
    different spatial scalings corresponding to increasing $k_\ast$. The
    number of spatial observations is fixed at $N_x=100$, with
    $\Delta X=L/N_x$, while $\gamma=2$ and $\tau_{\max}=1$ are fixed.
    The AUC and class conventions are as in
    Figure~\ref{fig:pitchfork_auc_vs_r}.
    }
    \label{fig:SH_L_auc}
\end{figure}

Figure~\ref{fig:SH_L_auc} shows the resulting AUC as a function of
$\mu$. Over the investigated range, the AUC generally increases with $L$,
showing that the direct spatial effect associated with increasing $k_\ast$
improves the distinction between the supercritical and subcritical
dynamics.

\end{document}